\documentclass[a4paper,Haag duality]{mathscan}
\newtheorem{thm}{Theorem}[section] 
\newtheorem{pro}[thm]{Proposition}

\theoremstyle{definition}           
\newtheorem{rem}[thm]{Remark}       
\newtheorem{defn}[thm]{Definition}  
\newcommand{\NI}{\noindent}

\newcommand{\bea}{\begin{eqnarray}}
\newcommand{\eea}{\end{eqnarray}}
\newcommand{\dsp}{\displaystyle}
\def \b #1 {\bf #1}
\newcommand{\IR}{I\!\!R}

\newcommand{\IC}{I\!\!C}

\newcommand{\IT}{I\!\!T}
\newcommand{\IZ}{Z\!\!\!Z}
\newcommand{\cal}{\mathcal}
\newcommand{\clk}{{\cal K}}

\newcommand{\cla}{{\cal A}}
\newcommand{\clz}{{\cal Z}}
\newcommand{\cli}{{\cal I}}

\newcommand{\clf}{{\cal F}}

\newcommand{\clh}{{\cal H}}
\newcommand{\clp}{{\cal P}}
\newcommand{\clo}{{\cal O}}
\newcommand{\clb}{{\cal B}}

\newcommand{\cle}{{\cal E}}
\newcommand{\clj}{{\cal J}}
\newcommand{\cln}{{\cal N}}
\newcommand{\cld}{{\cal D}}

\newcommand{\clm}{{\cal M}}

\newcommand{\raro}{\rightarrow}

\newcommand{\vsp}{\vskip 1em}

\def \qed {\hfill \vrule height6pt width 6pt depth 0pt}
\newcommand{\be}{\begin{equation}}
\newcommand{\ee}{\end{equation}}
\newcommand{\ben}{\begin{eqnarray*}}
\newcommand{\een}{\end{eqnarray*}}
\begin{document} 

\title{Translation invariant pure state on $\otimes_{\!Z}\!M_d(\!C)$ and it's split property }

\author{ Anilesh Mohari }
\thanks{...}

\address{ The Institute of Mathematical Sciences, CIT Campus, Taramani, Chennai-600113 }

\email{anilesh@imsc.res.in}



\thanks{ This paper has grown over the years starting with initial work in the 
middle of 2005. The author gratefully acknowledge discussion with Ola Bratteli and Palle E. T. Jorgensen 
for inspiring participation in sharing the intricacy of the present problem. Finally the author is indebted to 
Taku Matsui for valuable comments on an earlier draft of the present problem where the author made an attempt to prove 
Haag duality. }

\begin{abstract}

We prove Haag duality property of any translation invariant pure state on $\clb = \otimes_{\IZ}M_d(\!C), \;d \ge 2$, where $M_d(\!C)$ is the set of $d \times d$ dimensional matrices over field of complex numbers. We also prove a necessary and sufficient condition for a translation invariant factor state to be pure on $\clb$. This result makes it possible to study such a pure state with additional symmetry. We prove that exponentially decaying two point spacial correlation function of a real lattice symmetric reflection positive translation invariant pure state is a split state. Further there exists no translation invariant pure state on $\clb$ that is real, lattice symmetric, refection positive and $su(2)$ invariant when $d$ is an even integer. This in particular says that Heisenberg iso-spin anti-ferromagnets model for $1/2$-odd integer spin degrees of freedom admits
spontaneous symmetry breaking at it's ground states

\end{abstract}

\maketitle 

\section{ Introduction }

\vsp
Let $\omega$ be a translation invariant factor state on $\clb=\otimes_{\IZ}M_d(\!C)$ and $(\clh,\pi_{\omega},\Omega)$
be the GNS space associated with $(\clb,\omega)$. Let $\clb_R=\otimes_{n \ge 1}M_d(\!C)$ ( respectively $\clb_{L} =\otimes_{n < 1} M_d(\!C)$ be the right (left) sided $C^*$-sub-algebra of $\clb$ and 
$e_0$ ( respectively $\tilde{e}_0$ ) be the support projection of $\omega$ in 
$\pi_{\omega}(\clb_R)''$ (respectively in $\pi_{\omega}(\clb_L)''$. Thus we have two commuting projections 
$e_0,\tilde{e}_0$ defined by $e_0=[\pi_{\omega}(\clb_R)'\Omega]$ and $\tilde{e}_0=[\pi_{\omega}(\clb_L)'\Omega]$ as
$\pi_{\omega}(\clb_R)''$ commutes with $\pi_{\omega}(\clb_L)''$. We set $Q_0=e_0\tilde{e}_0$ and take $\clk_0$ 
to be the Hilbert subspace of $\clh$ determined by the projection $Q_0$. Also set von-Neumann algebras $\clm^1_0=Q_0\pi_{\omega}(\clb_R)''Q_0$ and $\tilde{\clm}^1_0=Q_0\pi_{\omega}(\clb_L)''Q_0$. So by our construction we
have $\tilde{\clm}^1_0 \subseteq (\clm^1_0)'$. Further the vector state $\phi_0(x)=<\Omega,x \Omega>$ on $\clb(\clk_0)$ is faithful and normal on $\clm_0^1$ and $\tilde{\clm}_0^1$. Further $\omega$ being a factor state, we will also have 
factor property of $\clm_0^1$ and $\tilde{\clm}_0^1$ ( Theorem 2.4 ). What is less obvious when we can expect 
cyclic property i.e. $[\clm^1_0\Omega]=[\tilde{\clm}^1_0\Omega]=I_{\clk_0}$, identity of $\clk_0$. This question is 
far from obvious to prove that the following statements are equivalent:

\NI (a) $[\clm^1_0\Omega]=I_{\clk_0},\;[\tilde{\clm}^1_0\Omega]=I_{\clk_0}$

\NI (b) $(\clm^1)'_0=\tilde{\clm}^1_0$;

\NI (c) $\pi_{\omega}(\clb_R)'=\pi_{\omega}(\clb_L)''$;

\NI (d) $[\pi_{\omega}(\clb_R)'\Omega]=[\pi_{\omega}(\clb_L)''\Omega]$;

\NI (e) $\omega$ is pure. 

\vsp 
Before we elaborate further on equivalence of above statements we briefly recall results on translation invariant pure state on $\clb =\otimes_{\IZ}M_d(\!C)$ that finds it's relevance while proving Haag duality (c). There is a one to one affine map between translation invariant states on $\clb$ and translation invariant states on $\clb_R =\otimes_{\IZ_+}M_d(\!C)$ by $\omega \raro \omega_R=\omega_{|\clb_R}$. The inverse map is the inductive limit 
state of $(\clb_R,\psi_R) \raro^{\lambda_n} (\clb_R,\psi_R)$ where $(\lambda_n:n \ge 0)$ is the canonical semi-group 
of right shifts on $\clb_R$. Pure states on a $\mbox{UHF}$ algebra are studied in the general framework of [Pow1]. Such a situation has been investigated also in details at various degrees of generality in [BJP] and [BJKW] with primary motivation to develop a $C^*$ algebraic method to study iterative function systems and its associated wavelet theory. One interesting result in [BJP] says that any translation invariant pure state on $\clb_R$ is also a product state and the canonical endomorphism associated with two such states are unitary equivalent. However such a statement is not true for two translation invariant pure states on $\clb$ as their restriction to $\clb_R$ need not be isomorphic. Theorem 3.4 in [Mo3] says that $\omega_R$ is either a type-I or a type-III factor state on $\clb_R$. 
Both type of factors are known to exists in the literature arises in statistical mechanics [BR2,Si,Ma1]. Thus the classification problem of translation dynamics on $\clb$ with invariant pure states on $\clb$ up to unitary isomorphism is a delicate one. However in [Mo3] we have got partial success when such states admits Kolmogorov's property introduces in [AM]. Since a $\theta$ invariant states $\omega$ on $\clb$ is completely determined by its restriction $\omega_R$ to $\clb_R$, in principle it is possible to describe various property of $\omega$ including pureness by studying certain properties of their restriction $\omega_R$. Since pureness of $\omega_R$ is not necessary for pureness of $\omega$, a valid fundamental question that arises here: What are the parameters that determines both $\omega$ and $\omega_R$ uniquely and how these parameters determines properties of $\omega$ and $\omega_R$?   

\vsp 
Theorem 7.1 in [BJKW] has aimed towards a sufficient condition on the associated minimal data Popescu elements for 
purity of the translation invariant state. However the proof is faulty as certain argument used in the proof is not time reversal symmetric. Thus Lemma 7.6 in [BJKW] needs additional assumption related to the support projection of the dual Cuntz's state. Besides this additional structure proof of Lemma 7.8 in [BJKW] is also not complete unless we find a proof for  $\tilde{\clm}=\clm'$ ( we retained same notations here in the text) for such a factor state $\omega$. Such a problem could have been solved if there were any method which shows directly that Takesaki's conditional expectation exists from $\clm'$ onto $\tilde{\clm}$. Thus main body of the proof for Theorem 7.1 in [BJKW] is incomplete. 

\vsp 
Translation invariant states on $\clb$ as an inductive limit states are investigated in a series of papers [Mo1],[Mo2],[Mo3]. Though we will be working here within the framework of amalgamated representation of $\tilde{\clo}_d 
\otimes \clo_d$ introduced in [BJKW], criteria on asymptotic behavior of translation invariant pure states is used to prove first that (d) indeed implies (e). The statement (d) implies (e) has originated from section 7 of [BJKW]. 
Converse problem is directly related to the main result of section 3 i.e. Haag duality $\pi_{\omega}(\clb_R)'=\pi_{\omega}(\clb_L)''$ for a pure translation invariant state (Theorem 3.6) $\omega$ on $\clb$. 

\vsp 
We explore the set of representation of $\clb$ quasi-equivalent to $\pi_{\omega}$ and equip it with 
a partial ordering to prove Haag duality. Mackey's system of imprimitivity plays a crucial role even though pure state not necessarily give rises a Mackey's system of imprimitivity generated by support projection $E_0$ with respect shift. Though we have worked here with amalgamated representation of $\tilde{\clo}_d \otimes \clo_d$, it seems that just for Haag duality one can avoid doing so. It seems that the underlining group $\IZ$ can easily be replaced by $\IZ^k$ for some $k \ge 2$ and wedge duality for a pointed cone can be proved following the same ideas. We defer this line of analysis leaving it for a possible future direction of work as it's relation with problems in quantum spin chain in higher dimensional lattice needs some additional structure. Haag duality property plays an important role in studying the factor state $\omega_R$ or $\omega_L$ when $\omega$ admits some additional symmetry apart from translation symmetry [Mo5] to determine split property [Ma1,Ma2] of a pure translation invariant state on $\clb$ and it's relation with decaying property of special correlation functions. 

\vsp
We briefly set the standard notations and known relations in the following text. The quantum spin chain that
we consider here is described by a {\bf UHF} $C^*$-algebra denoted by $\clb=\otimes_{\IZ}M_d(\!C)$. Here $\clb$
is the $C^*$ -completion of the infinite tensor product of the algebra ${\bf M}_d(\!C)$ of $d$ by $d$ complex matrices,
each component of the tensor product element is indexed by an integer $j$. Let $Q$ be a matrix in ${\bf M}_d(\!C)$. By
$Q^{(j)}$ we denote the element $...\otimes 1 \otimes 1 ... 1 \otimes Q \otimes 1 \otimes ... 1\otimes ,,. $, where $Q$
appears in the $j$-th component. Given a subset $\Lambda$ of $\!Z$, $\clb_{\Lambda}$ is defined as the $C^*$-sub-algebra
of $\clb$ generated by all $Q^{(j)}$ with $Q \in {\bf M}_d(\!C)$, $j \in \Lambda$. We also set
$$\clb_{loc}= \bigcup_{\Lambda:|\Lambda| < \infty } \clb_{\Lambda}$$
where $|\Lambda|$ is the cardinality of $\Lambda$. Let $\omega$ be a state on $\clb$. The restriction of $\omega$
to $\clb_{\Lambda}$ is denoted by $\omega_{\Lambda}$. We also set $\omega_{R}=\omega_{[1,\infty)}$ and $\omega_{L}=
\omega_{(-\infty,0]}$. The translation $\theta_k$ is an automorphism of $\clb$ defined by $\theta_k(Q^{(j)})=Q^{(j+k)}$.
Thus $\theta_1,\theta_{-1}$ are unital $*$-endomorphism on $\clb_R$ and $\clb_L$ respectively. We say $\omega$ is
translation invariant if $\omega \circ \theta_k = \omega$ on $\clb$ ( $\omega \circ \theta_1 = \omega$ on $\clb$ ).
In such a case $(\clb_R,\theta_1,\omega_{R})$ and $(\clb_L,\theta_{-1},\omega_{L})$ are two unital $*$-endomorphisms with
invariant states. Main result obtained in this paper given below.  

\begin{thm}
Let $\omega$ be a translation invariant factor state on $\clb$. Then the following holds:

\NI (a) $\omega$ is pure if and only if $[\pi_{\omega}(\clb_L)''\Omega]=[\pi_{\omega}(\clb_R)'\Omega]$ where $\pi_{\omega}$ 
is the $GNS$ representation of $\clb$ associated with $\omega$; 

\NI (b) $\omega$ is pure if and only if it admits Haag duality property i.e. $\pi_{\omega}(\clb_R)'=\pi_{\omega}(\clb_L)''$ (Theorem 3.6);
\end{thm}

\vsp 
A general mathematical question that is central here now: when and how can we guarantee that 
$\omega_R (\omega_L)$ are type-I factors or type-III factor by studying additional symmetry 
of the state? To that end we first recall [Ma2] a standard definition of a state 
to be split in the following.

\vsp
\begin{defn} 
Let $\omega$ be a translation invariant state on $\clb$. We say that $\omega$ is {\it split}
if the following condition is valid: Given any $\epsilon > 0$ there exists a $k \ge 1$ so that

\be
\mbox{sup}_{||Q|| < 1}|\omega(Q)-\omega_L \otimes \omega_R(Q)| \le \epsilon
\ee
where the above supremum is taken over all local elements $Q \in \clb_{(-\infty,-k] \cup [k, \infty)}$
with the norm less than $1$.
\end{defn}

\vsp
Here we recall few well known facts from [Pow1,BR,Ma1,Ma2]. The uniform cluster condition is valid if and only if the state $\omega$ is quasi-equivalent to the product state $\psi_L \otimes \psi_R$ of a state $\psi_L$ of $\clb_L$ and another state $\psi_R$ of $\clb_R$. Thus a Gibbs state of a finite range interaction is split. If $\omega$ is a pure translation invariant state, then $\omega_R( \omega_L )$ is type-I if and only if $\omega$ is also a split state. The canonical trace is a non-pure split state and unique ground state of XY model [AMa,Ma2] is a non-split pure state. One central aim is to find a criteria for a pure translation invariant state to be split. To that end we present a precise definition for exponential decay.

\vsp
\begin{defn} 
Let $\omega$ be a translation invariant state on one dimensional spin chain $\clb$. We say the two
point spacial correlation functions for $\omega$ {\it decay exponentially} if there exists a $\delta > 0$ so that
\be
e^{\delta k} |\omega(Q_1 \theta_k(Q_2))-\omega(Q_1)\omega(Q_2)| \raro 0
\ee
as $|k| \raro \infty$ for any local elements $Q_1,Q_2 \in \clb.$
\end{defn}

\vsp
For any translation invariant state $\omega$ on $\clb$ we set translation invariant state $\tilde{\omega}$ by reflecting around the point ${1 \over 2}$ on $\clb$ by $$\tilde{\omega}(Q_{-l}^{(-l)} \otimes Q_{-l+1}^{(-l+1)} \otimes ... \otimes Q_{-1}^{(-1)} \otimes Q_0^{(0)} \otimes Q_1^{(1)} ... \otimes Q_n^{(n)})$$
\be 
= \omega(Q_n^{(-n+1)}... \otimes Q_1^{(0)} \otimes Q_0^{(1)} \otimes Q_{-1}^{(2)} 
\otimes ... Q_{-l+1}^{(l)} \otimes Q_{-l}^{(l+1)})
\ee
for all $n,l \ge 1$ and $Q_{-l},..Q_{-1},Q_0,Q_1,..,Q_n \in 
M_n(\IC)$ where $Q^{(k)}$ is the matrix $Q$ at lattice point $k$. We define $\tilde{\omega}$ 
on $\clb$ by extending linearly to any $Q \in \clb_{loc}$. Thus $\omega \raro \tilde{\omega}$ is an
affine one to one onto map on the set of translation invariant states on $\clb$. 
Thus the state $\tilde{\omega}$ is translation invariant, ergodic, factor state if and only if $\omega$ is
translation invariant, ergodic, factor state respectively. We say $\omega$ is {\it lattice 
reflection symmetric} if $\omega=\tilde{\omega}$.   

\vsp 
If $Q= Q^{(l)}_0 \otimes Q^{(l+1)}_1 \otimes ....\otimes Q^{(l+m)}_m$ we set
$Q^t={Q^t}^{(l)}_0 \otimes {Q^t}^{(l+1)}_1 \otimes ..\otimes {Q^t}^{(l+m)}_m$
where $Q_0,Q_1,...,Q_m$ are arbitrary elements in $M_d$ and $Q_0^t,Q^t_1,..$ stands for transpose
with respect to an orthonormal basis $(e_i)$ for $\IC^d$ (not complex conjugate) of $Q_0,Q_1,..$ 
respectively. We define $Q^t$ by extending linearly for any
$Q \in \clb_{loc}$. For a state $\omega$ on UHF$_d$ $C^*$ algebra $\otimes_{\IZ} M_d$ we define
a state $\bar{\omega}$ on $\otimes_{\IZ} M_d$ by the following prescription
\be
\bar{\omega}(Q) = \omega(Q^t)
\ee
Thus the state $\bar{\omega}$ is translation invariant, ergodic, factor state if and only if $\omega$ is
translation invariant, ergodic, factor state respectively. We say $\omega$ is {\it real } if
$\bar{\omega}=\omega$.

\vsp 
A translation invariant state $\omega$ is said to be in {\it detailed balance } if $\omega$ is {\it lattice reflection symmetric } and {\it real} (for further details see section 3 ). The canonical trace on $\clb$ is both real and lattice symmetric. This notion of detailed balance state is introduced as an reminiscence of Onsager's reciprocal relations explored in recent articles [AM,Mo1,Mo2,Mo3] on non-commutative probability theory. Here we also recall well known notion [DLS] that a state $\omega$ on $\clb$ is called reflection positive if 
\be 
\omega(\overline{\clj(x)} x) \ge 0
\ee for all $x \in \clb_R$ where $\clj$ is the reflection map with a twist $g_0 \in U_d(\!C)$ from $\clb_R$ onto $\clb_L$ and $\overline{x}$ stands for complex conjugation with respect to a basis $(e_i)$ for $\!C^d$ applies globally on each matrices of the lattice simultaneously. 

\vsp 
Let $G$ be a compact group and $g \raro v(g)$ be a $d-$dimensional unitary representation of $G$. By $\gamma_g$ we 
denote the product action of $G$ on the infinite tensor product $\clb$ induced by $v(g)$,
\be 
\gamma_g(Q)=(..\otimes v(g) \otimes v(g)\otimes v(g)...)Q(...\otimes v(g)^*\otimes v(g)^*\otimes v(g)^*...)
\ee
for any $Q \in \clb$. We say $\omega$ is $G$-invariant if $\omega(\gamma_g(Q))=\omega(Q)$ for all $Q \in \clb_{loc}$. 

\vsp 
Our main results on symmetry in section 4 uses Haag duality property crucially to study lattice reflection symmetry property of a translation invariant pure state $\omega$. Here we make a short list of end results obtained in 
this paper. 

\vsp
\begin{thm} 
Let $\omega$ be a pure lattice symmetric translation invariant real ( with respect 
to a basis $(e_i)$ of $\IC^d$ ) state on $\clb$. Then the following holds:

\NI (a) If $\omega$ is also reflection positive with a twist $g_0$ and two point spacial correlation function for $\omega$ decays exponentially then $\omega$ is a split state 
i.e. $\pi_{\omega}(\clb_R)''$ is a type-I factor (Theorem 5.4).

\NI (b) If $\omega$ is also reflection positive and $G$-invariant then $g  \raro u(g)$
is a real representation with respect to the orthonormal basis $(e_i)$ for $\IC^d$ provided 
the invariant subspace of the representation $g \raro \overline{u(g)} \otimes u(g)$ is one 
dimensional.     
\end{thm}

\vsp 
\begin{thm} 
Let $H$ be a translation invariant Hamiltonian of the form $H=\sum_{ k \in \IZ} 
\theta_k(h_0)$ with $h_0=h_0^* \in \clb_{loc}$ and $H$ be real, lattice reflection symmetric and 
reflection positive with a twist $g_0 \in U_d(\IC)$. Let $H$ be also $SU(2)$ invariant where
$g \raro u(g)$ is an irreducible representation of $SU(2)$.  Then ground state for $H$ is not 
unique when $d=2s+1$ and $s$ is a half-odd integer.   
\end{thm} 

\vsp 
We end this paper with a brief application of our results to antiffero-magnet Heisenberg models 
which includes prime examples such as $H_{XY}$ and $H_{XXX}$ models.    

\vsp
The paper is organized as follows. In section 2 we study Popescu's dilationassociated with a translation invariant state on Cuntz algebra $\clo_d$ and review `commutant lifting theorem' investigated in [BJKW]. The proof presented here remove the murky part of the proof of Theorem 5.1 in [BJKW]. In section 3 we explore both the notion of Kolmogorov's shift and it's intimate relation with Mackey's imprimitivity system to explore a duality argument introduced in [BJKW]. We find a useful necessary and sufficient condition (Theorem 1.1 (a) ) in terms of support projection of Cuntz's state for a translation invariant factor state $\omega$ on $\clb$ to be pure. The criteria on support projection is crucial to prove our main mathematical result Theorem 3.6. 

\vsp 
Section 4 studies discrete symmetry and section 5 gives the proof of the statement (a) of Theorem 1.4. Section 6 studies continuous symmetry and gives proof of the statement (b) of Theorem 1.4 and also proof of Theorem 1.5. 

\begin{rem} 
The paper tittled ``On Haag Duality for Pure States of Quantum Spin Chain'' 
by authors: M. Keyl, Taku Matsui, D. Schlingemann, R. F. Werner, Rev. Math. Phys. 20:707-724,2008 has an incomplete proof for Haag duality property as Lemma 4.3 in that paper has a faulty argument.   
\end{rem}

\section{ States on $\clo_d$ and the commutant lifting theorem }

\vsp
In this section we essentially recall results from [BJKW] and organize it with additional remarks and arguments as 
it demands to understand the present problem investigated in section 3. First we recall that the Cuntz algebra 
$\clo_d ( d \in \{2,3,.., \} )$ is the universal $C^*$-algebra generated by the elements $\{s_1,s_2,...,s_d \}$ subject 
to the relations:

$$s^*_is_j=\delta^i_j1$$
$$\sum_{1 \le i \le d } s_is^*_i=1.$$

\vsp
There is a canonical action of the group $U(d)$ of unitary $d \times d$ matrices on $\clo_d$ given by
$$\beta_g(s_i)=\sum_{1 \le j \le d}\overline{g^j_i}s_j$$
for $g=((g^i_j) \in U(d)$. In particular the gauge action is defined by
$$\beta_z(s_i)=zs_i,\;\;z \in \IT =S^1= \{z \in \!C: |z|=1 \}.$$
If UHF$_d$ is the fixed point sub-algebra under the gauge action, then UHF$_d$ is the closure of the
linear span of all wick ordered monomials of the form
$$s_{i_1}...s_{i_k}s^*_{j_k}...s^*_{j_1}$$
which is also isomorphic to the UHF$_d$ algebra
$$M_{d^\infty}=\otimes^{\infty}_1M_d$$
so that the isomorphism carries the wick ordered monomial above into the matrix element
$$e^{i_1}_{j_1}(1)\otimes e^{i_2}_{j_2}(2) \otimes....\otimes e^{i_k}_{j_k}(k) \otimes 1 \otimes 1 ....$$
and the restriction of $\beta_g$ to $UHF_d$ is then carried into action
$$Ad(g)\otimes Ad(g) \otimes Ad(g) \otimes ....$$

\vsp
We also define the canonical endomorphism $\lambda$ on $\clo_d$ by
$$\lambda(x)=\sum_{1 \le i \le d}s_ixs^*_i$$
and the isomorphism carries $\lambda$ restricted to UHF$_d$ into the one-sided shift
$$y_1 \otimes y_2 \otimes ... \raro 1 \otimes y_1 \otimes y_2 ....$$
on $\otimes^{\infty}_1 M_d$. Note that $\lambda \beta_g = \beta_g \lambda $ on UHF$_d$.

\vsp
Let $d \in \{2,3,..,,..\}$ and $\IZ_d$ be a set of $d$ elements.  $\cli$ be the set of finite sequences
$I=(i_1,i_2,...,i_m)$ where
$i_k \in \IZ_d$ and $m \ge 1$. We also include empty set $\emptyset \in \cli$ and set $s_{\emptyset }=1=s^*_{\emptyset}$,
$s_{I}=s_{i_1}......s_{i_m} \in \clo_d $ and $s^*_{I}=s^*_{i_m}...s^*_{i_1} \in \clo_d$. In the following we recall 
a commutant lifting theorem ( Theorem 5.1 in [BJKW] ), crucial for our purpose.

\begin{thm}
Let $v_1,v_2,...,v_d$ be a family of bounded operators on a Hilbert space
$\clk$ so that $\sum_{1 \le k \le d} v_kv_k^*=I$. Then there exists a unique   up to
isomorphism Hilbert space $\clh$, a projection $P$ on $\clk$ and a family of isometries
$\{S_k:,\;1 \le k \le d,\;P\}$ satisfying Cuntz's relation so that
\be
PS^{*}_{k}P=S_k^*P=v^*_k
\ee
for all $1 \le k \le d$ and $\clk$ is cyclic for the representation i.e. the vectors
$\{ S_I\clk: |I| < \infty \}$ are total in $\clh$.

Moreover the following holds:

\NI (a) $\Lambda_n(P) \uparrow I$ as $n \uparrow \infty$;

\NI (b) For any $D \in \clb_{\tau}(\clk)$, $\Lambda_n(D) \raro X'$ weakly as $n \raro \infty$
for some $X'$ in the commutant $\{S_k,S^*_k: 1 \le k \le d \}'$ so that $PX'P=D$. Moreover
the self adjoint elements in the commutant $\{S_k,S^*_k: 1 \le k \le d \}'$ is isometrically
order isomorphic with the self adjoint elements in $\clb_{\tau}(\clk)$ via the
surjective map $X' \raro PX'P$, where $\clb_{\tau}(\clk)=\{x \in \clb(\clk): \sum_{1 \le k \le d }
v_kxv^*_k=x \}.$

\NI (c) $\{v_k,v^*_k,\;1 \le k \le d \}' \subseteq \clb_{\tau}(\clk)$ and equality holds if and only if
$P \in \{S_k,S_k,\;1 \le k \le d \}''$.

\vsp 
If $(w_i)$ be another such an Popescu elements on a Hilbert space $\clk'$ such that there exists an operator 
$u:\clk \raro \clk'$ so that $\sum_k w_kuv_k^*=u$ then there exists an operator $U:\clh_v \raro \clh_w$ so 
that $\pi'(x)U=U\pi(x)$ where $(\clh_w,\pi',S'_i)$ are Popescu dilation of $(w_i)$ and $\pi'$ is the associated
minimal representation of $\clo_d$. In particular $U$ is isometry, unitary if $u$ is so respectively. If $u$ is unitary 
and $\clk=\clk'$ then we can as well take $\clh_v=\clh_w$.  

\end{thm}

\vsp
\NI {\bf PROOF: } Following Popescu [Po] we define a completely positive map
$R: \clo_d \raro \clb(\clk)$ by
\be
R(s_Is^*_J)=v_Iv^*_J
\ee
for all $|I|,|J| < \infty$. The representation $S_1,..,S_d$ of $\clo_d$ on $\clh$ thus may be taken to be the Stinespring
dilation of $R$ [BR, vol-2] and uniqueness up to unitary equivalence follows from uniqueness of the Stinespring
representation. That $\clk$ is cyclic for the representation follows from the minimality property of the Stinespring
dilation. For (a) let $Q$ be the limiting projection. Then we have $\Lambda(Q)=Q$, hence $Q \in \{S_k,S^*_k \}'$
and $Q \ge P$. In particular $QS_If=S_If$ for all $f \in \clk$ and $|I| < \infty$. Hence $Q=I$ by the
cyclicity of $\clk$. For (b) essentially we deffer from the argument used in Theorem 5.1 in [BJKW]. We fix any $D
\in \clb_{\tau}(\clk)$ and note that $P\Lambda_k(D)P=\tau_k(D)=D$ for any $k \ge 1$. Thus for any integers $n > m$
we have
$$\Lambda_m(P)\Lambda_n(D)\Lambda_m(P)=\Lambda_m(P\Lambda_{n-m}(D)P)=\Lambda_m(D)$$
Hence for any fix $m \ge 1$ limit $<f,\Lambda_n(D)g>$ as $n \raro \infty$ exists for all
$f,g \in \Lambda_m(P)$. Since the family of operators $\Lambda_n(D)$ is uniformly bounded and $\Lambda_m(P) \uparrow I$
as $m \raro \infty$, a standard density argument guarantees that the weak operator limit of $\Lambda_n(D)$ exists as
$n \raro \infty$. Let $X'$ be the limit. So $\Lambda(X')=X'$, by Cuntz's relation, $X' \in \{S_k,S^*_k:1 \le k \le k \}'$.
Since $P\Lambda_n(D)P=D$ for all $n \ge 1$, we also conclude that $PX'P=D$ by taking limit $n \raro \infty$. Conversely
it is obvious that $P\{ S_k,S^*_k:\; k \ge 1 \}'P \subseteq \clb_{\tau}(\clk)$. Hence
we can identify $P\{S_k,S^*_k:\; k \ge 1 \}'P$ with $\clb_{\tau}(\clk)$.

\vsp
Further it is obvious that $X'$ is self-adjoint if and only if $D=PX'P$ is self-adjoint. Now fix any self-adjoint
element $D \in \clb_{\tau}(\clk)$. Since identity operator on $\clk$ is an element in $\clb_{\tau}(\clk)$ for any
$\alpha \ge 0$ for which $- \alpha P  \le D \le \alpha P$, we have $\alpha \Lambda_n(P) \le \Lambda_n(D)
\le \alpha \Lambda_n(P)$ for all $n \ge 1$. By taking limit $n \raro \infty$ we
conclude that $- \alpha I \le X' \le \alpha I$, where $PX'P=D$. Since operator norm of a self-adjoint element $A$
in a Hilbert space is given by
$$||A||= \mbox{inf}_{\alpha \ge 0}\{\alpha: - \alpha I  \le A \le \alpha I \}$$
we conclude that $||X'|| \le ||D||$. That $||D||=||PX'P|| \le ||X'||$ is obvious, $P$ being a
projection. Thus the map is isometrically order isomorphic taking self-adjoint elements of the commutant
to self-adjoint elements of $\clb_{\tau}(\clk)$.

\vsp
We are left to prove (c). Inclusion is trivial.  For the last part note that for any invariant element $D$
in $\clb(\clk)$ there exists an element $X'$ in $\{S_k,S^*_k,\;1 \le k \le d\}'$ so that $PX'P=D$. In such a
case we verify that $Dv^*_k= PX'PS^*_kP=PX'S^*_kP=PS^*_kX'P=PS^*_kPX'P=v^*_kD$. We also have $D^* \in \clb_
{\tau}(\clk)$ and thus $D^*v^*_k=v^*_kD^*$. Hence $D \in \{v_k,v^*_k: 1 \le k \le d \}'$. Since $P\pi_{\hat{\omega}}
(\clo_d)'P = \clb(\clk)_{\tau}$, we conclude that $\clb(\clk)_{\tau} \subseteq \clm'$. Thus equality holds whenever
$P \in \{S_k,S_k^*,\;1 \le k \le d \}''$. For converse note that by commutant lifting property self-adjoint elements
of the commutant $\{S_k,S^*_k,1 \le k \le d \}'$ is order isometric with the algebra $\clm'$ via the map $X' \raro PX'P$.
Hence $P \in \{S_k,S^*_k,1 \le k \le d \}''$ by Proposition 4.2 in [BJKW]. 

\vsp 
For the proof of intertwining relation and their property we refer to main body of the proof of 
Theorem 5.1 in [BJKW] \qed

\vsp
A family $(v_k,1 \le k \le d)$ of contractive operators on a Hilbert space $\clk$ is called Popescu's elements and 
dilation $(\clh,P,\clk,S_k,1 \le k \le d)$ in Theorem 2.1 is called Popescu's dilation to Cuntz elements. In the 
following proposition we deal with a family of minimal Popescu elements for a state on $\clo_d$. 

\begin{pro}
There exists a canonical one-one correspondence between the following objects:

\vsp
\NI (a) States $\psi$ on $\clo_d$

\vsp
\NI (b) Function $C: \cli \times \cli \raro \!C$ with the following properties:

\NI (i) $C(\emptyset, \emptyset)=1$;

\NI (ii) for any function $\lambda:\cli \raro \!C$ with finite support we have
         $$\sum_{I,J \in \cli} \overline{\lambda(I)}C(I,J)\lambda(J) \ge 0$$

\NI (iii) $\sum_{i \in \IZ_d} C(Ii,Ji) = C(I,J)$ for all $I,J \in \cli$.

\vsp
\NI (c) Unitary equivalence class of objects $(\clk,\Omega,v_1,..,v_d)$ where

\NI (i) $\clk$ is a Hilbert space and $\Omega$ is an unit vector in $\clk$;

\NI (ii) $v_1,..,v_d \in \clb(\clk)$ so that $\sum_{i \in \IZ_d} v_iv^*_i=1$;

\NI (iii) the linear span of the vectors of the form $v^*_I\Omega$, where $I \in \cli$, is dense in $\clk$.

\vsp
Where the correspondence is given by a unique   completely positive map $R: \clo_d \raro \clb(\clk)$ so that

\NI (i) $R(s_Is^*_J)=v_Iv^*_J;$

\NI (ii) $\psi(x)=<\Omega,R(x)\Omega>;$

\NI (iii) $\psi(s_Is^*_J)=C(I,J)=<v^*_I \Omega,v^*_J\Omega>;$

\NI (iv) For any fix $g \in U_d$ and the completely positive map $R_g:\clo_d \raro \clb(\clk)$ defined by $R_g= R
\circ \beta_g$ give rises to a Popescu system given by $(\clk,\Omega,\beta_g(v_i),..,\beta_g(v_d))$ where
$\beta_g(v_i)=\sum_{1 \le j \le d} \overline{g^i_j} v_j.$

\end{pro}

\vsp
\NI{\bf PROOF: } For a proof we simply refer to Proposition 2.1 in [BJKW]. \qed

\vsp
The following is a simple consequence of Theorem 2.1 valid for a $\lambda$-invariant state $\psi$ on $\clo_d$. 
This proposition will have very little application in main body of this paper but this gives a clear picture explaining the
delicacy of the present problems. 

\begin{pro}

Let $\psi$ be a state on $\clo_d$ and $(\clh,\pi,\Omega)$ be the GNS space associated with 
$(\clo_d,\psi)$. We set $S_i=\pi(s_i)$ and normal state $\psi_{\Omega}$ on $\pi(\clo_d)''$ defined by 
$$\psi_{\Omega}(X) = <\Omega, X \Omega>$$
Let $P$ be the projection on the closed subspace $\clk$ generated by the vectors $\{S^*_I\Omega: |I| < \infty \}$ and 
\be
v_k=PS_kP
\ee
for $1 \le k \le d$. Then following holds:

\vsp
\NI (a) $\{v^*_I\Omega: |I| < \infty \}$ is total in $\clk$.

\NI (b) $\sum_{1 \le k \le d} v_kv_k^*=I$;

\NI (c) $S_k^*P=PS_k^*P$ for all $1 \le k \le d$; 

\NI (d) For any $I=(i_1,i_2,...,i_k),J=(j_1,j_2,...,j_l)$ with $|I|,|J| < \infty$ we have 
\be 
\psi(s_Is^*_J)= <\Omega,v_Iv^*_J\Omega>
\ee 
and the vectors $\{ S_If: f \in \clk,\;|I| < \infty \}$ are total in
the GNS Hilbert space associated with $(\clo_d,\psi)$. Further such a family $(\clk,\;v_k,\;1 \le k \le d,\;\omega)$ satisfying
(a) to (d) are determined uniquely up to isomorphism. 

\vsp
Conversely given a Popescu system $(\clk,v_k,\;1 \le k \le d,\Omega)$ satisfying (a) and (b) there exists a unique   state 
$\psi$ on $\clo_d$ so that (c) and (d) are satisfied.

Furthermore the following statements are valid:
\vsp
\NI (e)  If the normal state $\phi_0(x)=<\Omega,x \Omega>$ on the von-Neumann algebra
$\clm=\{ v_i, v^*_i \}''$ is invariant for the Markov map $\tau(x)=\sum_{1 \le k
\le d} v_ixv^*_i,\;x \in \clm$ then $\psi$ is $\lambda$ invariant and $\phi_0$ is faithful 
on $\clm$.

\vsp
\NI (f) If $P \in \pi(\clo)''$ then following are equivalent:

\NI (i) $\psi$ is an ergodic state for $(\clo_d,\lambda)$; 

\NI (ii) $(\clm,\tau,\phi_0)$ is ergodic.

In such a case $\clm$ is a factor. 
\end{pro} 

\vsp
\NI {\bf PROOF:} We fix a state $\psi$ and consider the GNS space  $(\clh,\pi,\Omega)$
associated with $(\clo_d,\psi)$ and set $S_i=\pi(s_i)$. It is obvious 
that $S^*_kP \subseteq P$ for all $1 \le k \le d $, thus $P$ is the minimal subspace 
containing $\Omega$ and invariant by all $\{ S^*_k;\;1 \le k \le d \}$ i.e.
\be
PS^*_kP=S^*_kP
\ee
Thus $v^*_k=PS^*_kP=S_k^*P$ and so $\sum_k v_kv^*_k = \sum_k PS_kS^*_kP= P$ which is 
identity operator in $\clk$. This completes the proof of (a) (b) and (c).

\vsp
For (d) we note that
$$\psi(s_Is^*_J)=<\Omega, S_IS^*_J\Omega>$$
$$=<\Omega,PS_IS^*_JP\Omega>=<\Omega,v_Iv^*_J\Omega>.$$
Since $\clh$ is spanned by the vectors $\{S_IS^*_J\Omega:|I|,|J| < \infty \}$ and
$\clk$ is spanned by the vectors $\{ S^*_J\Omega=v^*_J\Omega:|I| < \infty \}$, 
$\clk$ is cyclic for $S_I$ i.e. the vectors $\{ S_I\clk: |I| < \infty \}$ spans 
$\clh$. Uniqueness up to isomorphism follows as usual by total property of vectors 
$v^*_I\Omega$ in $\clk$. 

\vsp
Conversely for a Popescu's elements $(\clk,v_i,\Omega)$ satisfying (a) and (b), we consider
the family $(\clh,S_k,\;1 \le k \le d,P)$ of Cuntz's elements defined as in
Theorem 2.1. We claim that $\Omega$ is a cyclic vector for the representation $\pi(s_i) \raro
S_i$. Note that by our construction vectors $\{S_If,f \in \clk: |I| < \infty \}$ are total
in $\clh$ and $v^*_J\Omega=S^*_J\Omega$ for all $|J| < \infty$. Thus by our hypothesis that
vectors $\{v^*_J\Omega:|I| < \infty \}$ are total in $\clk$, we verify that
vectors $\{S_IS^*_J\Omega: |I|,|J| < \infty \}$ are total in $\clh$. Hence $\Omega$ is a cyclic for
the representation $s_i \raro S_i$ of $\clo_d$.

\vsp
We left to prove (e) and (f). It simple to note by (d) that $\psi \lambda=\psi$ i.e.
$$\sum_i <\Omega,S_iS_IS^*_JS^*_i\Omega>= \sum_i <\Omega,v_iv_Iv^*_Jv^*_i\Omega> $$
$$=<\Omega,v_Iv_J^*\Omega> = <\Omega,S_IS^*_J\Omega>$$ for all $|I|,|J| < \infty$ where 
in the second equality we have used our hypothesis that the vector state $\phi_0$ on $\clm$ 
is $\tau$-invariant. In such case we aim now to show that $\phi_0$ is faithful on $\clm$. To that end 
let $p'$ be the support projection in $\clm$ for $\tau$ invariant state $\phi_0$. Thus $\phi_0(1-p') = 0$ 
i.e. $p'\Omega=\Omega$ and by invariance we also have $\phi_0(p'\tau(1-p')p')=\phi_0(1-p')=0$. Since $p'\tau(1-p')p' \ge 0 $ and an element in 
$\clm$, by minimality of support projection, we 
conclude that $p'\tau(1-p')p'=0$. Hence $p'\Omega=\Omega$ and $p'v^*_kp'=v^*_kp'$ for all $1 \le k \le d$. 
Thus $p'v^*_I\Omega=v^*_I\Omega$ for all 
$|I| < \infty$. As $\clk$ is the closed linear span of the vectors $\{ v_I^*\Omega: |I| < \infty \}$, we conclude 
that $p'=p$. In other words $\phi_0$ is faithful on $\clm$. This completes the proof for (e).

\vsp
We are left to show (f) where we assume that $P \in \pi(\clo_d)''$. $\Omega$ being a cyclic vector for $\pi(\clo_d)''$, the weak$^*$ limit of the increasing projection $\Lambda^k(P)$ is $I$. Thus by Theorem 3.6 in [Mo1] we have $(\pi(\clo_d)'',\Lambda,\psi_{\Omega})$ is ergodic if and only if the reduced dynamics $(\clm,\tau,\phi_0)$ is ergodic. Last part of the statement is an easy consequence of a theorem in [Fr] (see also [BJKW],[Mo1]). \qed

\vsp
Before we move to next result we comment here that in general for a $\lambda$ invariant state on $\clo_d$ 
the normal state $\phi_0$ on $\clm=\{v_k,v_k^*:1 \le k \le d \}''$ need not be invariant for $\tau$. To that 
end we consider ( [BR] vol-II page 110 ) the unique KMS state $\psi=\psi_{\beta}$ for the automorphism $\alpha_t(s_i)=e^{it}s_i$ on $\clo_d$. $\psi$ is $\lambda$ invariant and $\psi_{|}\mbox{UHF}_d$ is the unique faithful trace. $\psi$ being a KMS state for an automorphism, the normal state induced by the cyclic vector on $\pi_{\psi}(\clo_d)''$ is also separating for $\pi(\clo_d)''$. As $\psi \beta_z = \psi$ for all $z \in S^1$ we have 
$<\Omega,\pi(s_I)\Omega>=<\Omega,\beta_z(s_I)\Omega> = z^{|I|}<\Omega,\pi(s_I)\Omega>$ for all $z \in S^1$ and 
so $<\Omega,\pi(s_I)\Omega>=0$ for all $|I| \ge 1$. In particular $<\Omega,v_I^*\Omega>=0$ where $(v_i)$ are defined 
as Proposition 2.3 and thus $<v_i\Omega,v^*_I\Omega>=<\Omega,v^*_iv_I\Omega>=0$ for all $1 \le i \le d$. Hence $v_i\Omega=0$. By Proposition 2.3 (e), $\Omega$ is separating for $\clm$ and so we get $v_i=0$ for all $1 \le i \le d$ and this contradicts that $\sum_i v_iv_i^*=1$. Thus we conclude by Proposition 2.3 (e) that $\phi_0$ is not $\tau$ invariant on $\clm$. 
This example also indicates that the support projection of a $\lambda$ invariant state $\psi$ in $\pi(\clo_d)''$ need 
not be equal to the minimal sub-harmonic projection $P$ i.e. the closed span of vectors $\{ S_I^*\Omega: |I| < \infty \}$ containing $\Omega$ and $\{v_Iv_J^*: |I|,|J| < \infty \}$ need not be even an algebra.        

\vsp
Now we aim to deal with another class of Popescu elements associated with an $\lambda$-invariant state on $\clo_d$. 
In fact this class of Popescu elements will play a significant role for the rest of the text and we will repeatedly use
this proposition! 

\vsp
\begin{pro} 

Let $(\clh,\pi,\Omega)$ be the GNS representation of a $\lambda$ invariant state 
$\psi$ on $\clo_d$ and $P$ be the support projection of the normal state $\psi_{\Omega}(X)=<\Omega,X\Omega>$ in the 
von-Neumann algebra $\pi(\clo_d)''$. Then the following holds:

\vsp
\NI (a) $P$ is a sub-harmonic projection for the endomorphism $\Lambda(X)=\sum_k S_kXS^*_k$ on $\pi(\clo_d)''$
i.e. $\Lambda(P) \ge P$ satisfying the following:

\NI (i) $\Lambda_n(P) \uparrow I$ as $n \uparrow \infty$;

\NI (ii) $PS^*_kP=S^*_kP,\;\;1 \le k \le d$;

\NI (iii) $\sum_{1 \le k \le d} v_kv_k^*=I$ 

where $S_k=\pi(s_k)$ and $v_k=PS_kP$ for $1 \le k \le d$;

\vsp
\NI (b) For any $I=(i_1,i_2,...,i_k),J=(j_1,j_2,...,j_l)$ with $|I|,|J| < \infty$ we have $\psi(s_Is^*_J) =
<\Omega,v_Iv^*_J\Omega>$ and the vectors $\{ S_If: f \in \clk,\;|I| < \infty \}$ are total in $\clh$;

\vsp
\NI (c) The von-Neumann algebra $\clm=P\pi(\clo_d)''P$, acting on the Hilbert space
$\clk$ i.e. range of $P$, is generated by $\{v_k,v^*_k:1 \le k \le d \}''$ and the normal state
$\phi_0(x)=<\Omega,x \Omega>$ is faithful on the von-Neumann algebra $\clm$.

\NI (d) The self-adjoint part of the commutant of $\pi(\clo_d)'$ is norm and order isomorphic to the space
of self-adjoint fixed points of the completely positive map $\tau$. The isomorphism takes $X' \in \pi(\clo_d)'$
onto $PX'P \in \clb_{\tau}(\clk)$, where $\clb_{\tau}(\clk) = \{ x \in \clb(\clk): \sum_kv_kxv^*_k=x \}$. Furthermore
$\clm' = \clb_{\tau}(\clk)$.

\vsp
Conversely let $\clm$ be a von-Neumann algebra generated by a family $\{v_k: 1 \le k \le d \}$ of bounded operators
on a Hilbert space $\clk$ so that $\sum_kv_kv_k^*=1$ and the commutant $\clm'=\{x \in \clb(\clk): \sum_k v_kxv_k^*=x
\}$. Then the Popescu dilation $(\clh,P,S_k,\;1 \le k \le d)$ described in Theorem 2.1 satisfies the following:

\NI (i) $P \in \{S_k,S^*_k,\;1 \le k \le d \}''$;

\NI (ii) For any faithful normal invariant state $\phi_0$ on $\clm$ there exists a state $\psi$ on $\clo_d$ defined by
$$\psi(s_Is^*_J)=\phi_0(v_Iv^*_J),\;|I|,|J| < \infty $$
so that the GNS space associated with $(\clm,\phi_0)$ is the support projection
for $\psi$ in $\pi(\clo_d)''$ satisfying (a)-(d). 

\vsp
Further for a given $\lambda$-invariant state $\psi$, the family $(\clk,\clm,v_k\;1 \le k \le d,\phi_0)$ satisfying (a)-(d) 
is determined uniquely up to unitary conjugation. 

\vsp
\NI (e) $\phi_0$ is a faithful normal $\tau$-invariant state on $\clm$. Furthermore the following statements are equivalent:

\NI (i) $(\clo_d,\lambda,\psi)$ is ergodic;

\NI (ii) $(\clm,\tau,\phi_0)$ is ergodic;

\NI (iii) $\clm$ is a factor.

\end{pro} 

\vsp
\NI {\bf PROOF: } $\Lambda(P)$ is also a projection in $\pi_{\psi}(\clo_d)''$ so that
$\psi_{\Omega}(\Lambda(P))=1$ by invariance property. Thus we have $\Lambda(P) \ge P$ i.e.
$P\Lambda(I-P)P=0$. Hence we have
\be
PS^*_kP=S^*_kP
\ee
Moreover by $\lambda$ invariance property we also note that the faithful normal state
$\phi_0(x)=<\Omega,x\Omega>$ on the von-Neumann algebra $\clm=P\pi_{\psi}(\clo_d)''P$
is invariant for the reduce Markov map [Mo1] on $\clm$ given by
\be
\tau(x)=P\Lambda(PxP)P
\ee
\vsp
We claim that $\mbox{lim}_{n \uparrow \infty}\Lambda^n(P)= I$. That $\{ \Lambda^n(P): n \ge 1 \}$
is a sequence of increasing projections follows from sub-harmonic property of $P$ and endomorphism
property of $\Lambda$. Let the limiting projection be $Y$. Then $\Lambda(Y)=Y$ and so $Y \in
\{S_k,S^*_k \}'$. Since by our construction GNS Hilbert space $\clh_{\pi_{\hat{\omega}}}$ is
generated by $S_IS^*_J\Omega$, $Y$ is a scaler, being a non-zero projection, it is
the identity operator in $\clh_{\pi_{\psi}}$.

\vsp
Now it is routine to verify (a) (b) and (c). For the first part of (d) we appeal to Theorem 2.2. For the last part
note that for any invariant element $D$ in $\clb(\clk)$ there exists an element $X'$ in $\pi(\clo_d)'$ so that $PX'P=D$.
Since $P \in \pi(\clo_d)''$ we note that $(1-P)X'P=0$.
Now since $X' \in \{ S_k,S^*_k \}'$, we verify that $Dv^*_k= PXPS^*_kP=PXS^*_kP=PS^*_kXP=PS^*_kPXP=v^*_kD$. Since $D^* \in
\clb_{\tau}(\clk)$ we also have $D^*v^*_k=v^*_kD^*$. Thus $D \in \{v_k,v^*_k: 1 \le k \le d \}'=\clm'$. Since
$P\pi_{\hat{\omega}}(\clo_d)'P=
\clb(\clk)_{\tau}$, we conclude that $\clb(\clk)_{\tau} \subseteq \clm'$. The reverse inclusion is trivial.
This completes the proof for (d).

\vsp
For the converse part of (i), since by our assumption and commutant lifting property self-adjoint elements of the commutant
$\{S_k,S^*_k,1 \le k \le d \}'$ is order isometric with the algebra $\clm'$ via the map $X' \raro PX'P$,
$P \in \{S_k,S^*_k,1 \le k \le d \}''$ by Proposition 4.2 in [BJKW]. For (ii) without loss of generality assume
that $\phi_0(x)=<\Omega,x\Omega>$ for all $x \in \clm$ and $\Omega$ is a cyclic and separating vector for $\clm$.
( otherwise we set state $\psi(s_Is^*_J)=\phi_0(v_Iv_J^*)$ and consider it's GNS representation )
We are left to show that $\Omega$ is a cyclic vector for the representation $\pi(s_i) \raro S_i$. To that end 
let $Y \in \pi(\clo_d)'$ be the projection on the subspace generated by the vectors $\{S_IS^*_J\Omega:  
|I|,|J| < \infty \}$. Note that $P$ being an element in $\pi(\clo_d)''$, $Y$ also commutes with all the element 
$P\pi(\clo_d)''P=P \clm P$. Hence $Yx\Omega=x\Omega$ for all $x \in \clm$. Thus $Y \ge P$.  
Since $\Lambda_n(P) \uparrow I$ as $n \uparrow \infty$ by our construction, we conclude that $Y=\Lambda_n(Y) 
\ge \Lambda_n(P) \uparrow I$ as $n \uparrow \infty$. Hence $Y=I$. In other words
$\Omega$ is cyclic for the representation $s_i \raro S_i$. This completes the proof for (ii). 

\vsp
Uniqueness up to unitary isomorphism follows as GNS representation is determined uniquely unto unitary conjugation 
and so its support projection.  

\vsp
The first part of (e) we note that $PS_IS^*_JP=v_Iv_J^*$ for all $|I|,|J| < \infty$ and thus $\clm=P\pi(\clo_d)''P$ is the von-Neumann 
algebra generated by $\{v_k,v_k^*:1 \le k \le d \}$ and thus $\tau(x)=P\Lambda(PxP)P$ for all $x \in \clm$. That $\phi_0$ is $\tau(x)=\sum_k v_kxv_k^*$ 
invariant follows as $\psi$ is $\lambda$-invariant. We are left to prove equivalence of statements (i)-(iii).  

\vsp
By Theorem 3.6 in [Mo1] Markov semi-group $(\clm,\tau,\phi_0)$ is ergodic if and
only if $(\pi(\clo_d)'',\Lambda,\psi_{\Omega})$ is ergodic ( here we need to recall by (a) that
$\Lambda_n(P) \uparrow I$ as $n \uparrow \infty$ ). By a standard result [Fr, also BJKW] $(\clm,\tau,\phi_0)$
is ergodic if and only if there is no non trivial projection $e$ invariant for $\tau$ i.e. $\cli^{\tau} = \{e \in \clm: e^*=e,e^2=e,\tau(e)=e \}=\{0,1 \}$. If $\tau(e)=e$ for some projection $e \in \clm $ then $(1-e)\tau(e)(1-e)=0$ and so $ev^*_k(1-e)=0$. Same is true 
if we replace $e$ by $1-e$ as $\tau(1)=1$ and $\tau(1-e)=1-\tau(e)=1-e$ and thus $(1-e)v^*_ke=0$. Thus $e$ commutes with $v_k,v_k^*$ for all $1 \le k \le d$. Hence $\cli^{\tau} \subseteq \clm \bigcap \clm'$. Inequality in the reverse direction is trivial and thus
$\cli^{\tau}$ is trivial if and only if $\clm$ is a factor. Thus equivalence of (ii) and (iii) follows by a standard result 
[Fr] in non-commutative ergodic theory. This completes the proof. \qed

\vsp
The following two propositions are essentially easy adaptations of results appeared in [BJKW, Section 6 and Section 7],
crucial in our present framework. 

\begin{pro} 
Let $\psi$ be a $\lambda$ invariant factor state on $\clo_d$ and $(\clh,\pi,\Omega)$
be it's GNS representation. Then the following holds:

\NI (a) The closed subgroup $H=\{z \in S^1: \psi \beta_z =\psi \}$ is equal to 

$$\{z \in S^1: \beta_z \mbox{extends to an automorphism of } \pi(\clo_d)'' \} $$ 

\NI (b) Let $\clo_d^{H}$ be the fixed point sub-algebra in $\clo_d$ under the gauge group $\{ \beta_z: z \in H \}$. Then  
$\pi(\clo_d^{H})'' = \pi(\mbox{UHF}_d)''$.

\NI (c) If $H$ is a finite cyclic group of $k$ many elements and $\pi(\mbox{UHF}_d)''$ is a factor, 
then $\pi(\clo_d)'' \bigcap \pi(\mbox{UHF}_d)' \equiv \!C^m$ where $1 \le m \le k$.  
\end{pro}

\vsp
\NI {\bf PROOF: } It is simple that $H$ is a closed subgroup. For any fix $z \in H $ we define unitary operator $U_z$ 
extending the map $\pi(x)\Omega \raro \pi(\beta_z(x))\Omega$ and check that the map $X \raro U_zXU^*_z$ extends $\beta_z$ 
to an automorphism of $\pi(\clo_d)''$. For the converse we will use the hypothesis that $\psi$ is a $\lambda$-invariant 
factor state and $\beta_z \lambda= \lambda \beta_z$ to guarantee that $\psi \beta_z (X) = {1 \over n} 
\sum_{1 \le k \le n} \psi \lambda^k \beta_z(X) = {1 \over n} \sum_{1 \le k \le n} \psi \beta_z \lambda^k(X) 
\raro \psi(X)$ as $n \raro \infty$ for any $X \in \pi(\clo_d)''$, where we have used the same symbol $\beta_z$ 
for the extension. Hence $z \in H$. 

\vsp
For any $z_1,z_2 \in S^1$ we extend both $\psi \beta_{z_1}$ and $\psi \beta_{z_2}$ to its inductive limit state on $\clo_d^*$ 
using the canonical endomorphism $\clo_d \raro^{\lambda} \clo_d$. Inductive limit state being an affine map, their inductive limit 
states are also factors. The inductive limit of the canonical endomorphism became an automorphism. $(\clb,\theta)$ is asymptotically abelian 
i.e. $||x\theta^n(y)-\theta^n(y)x|| \raro 0$ as $n \raro \infty$ for all $x,y \in \clb$ (see also page 240 in [BR, vol2]). Thus in 
particular $(\clb,\IZ,\omega)$ is $\IZ-$central for any translation invariant ergodic state $\omega$ (see page 380 in [BR vol-2]). 
Thus we may appeal to a general result in $C^*$-non-commutative ergodic theory to conclude that their inductive limit, being translation 
invariant factor states, are either same or orthogonal (Theorem 4.3.19 in [BR vol2]).  

\vsp
In the following instead of working with $\clo_d$ we should be working with the inductive limit $C^*$ algebra and their inductive limit states. For simplicity of notation we still use $UHF_d,\clo_d$ for its inductive limit of $\clo_d \raro^{\lambda} \clo_d$ 
and $UHF_d \raro^{\lambda} UHF_d$ respectively and so for its inductive limit states.   

\vsp
Now we aim to prove (b). $H$ being a closed subgroup of $S^1$, it is either entire $S^1$ or a finite subgroup 
$\{exp({ 2i \pi l \over k})|l=0,1,...,k-1 \}$ where the integer $k \ge 1$. If $H=S^1$ we have nothing to prove 
for (b). When $H$ is a finite closed subgroup, we identify $[0,1)$ with $S^1$ by the usual map and note that 
if $\beta_t$ is restricted to $t \in [0,{1 \over k})$, then by scaling we check that $\beta_t$ defines a 
representation of $S^1$ in automorphisms of $\clo^H_d$. Now we consider the direct integral representation 
$\pi'$ defined by
$$\pi'= \int^{\oplus}_{[0,{1 \over k})} dt \pi_{|_{\clo^H_d}} \beta_t $$
of $\clo_d^{H}$ on $\clh_{|_{\clo_d^H}}\otimes L^2([0,{1 \over k})\;)$, where $\clh_{|_{\clo_d^H}}$ is the cyclic
space of $\pi(\clo_d^H)$ generated by $\Omega$. That it is indeed direct integral follows as states
$\psi \beta_{t_1}$ and $\psi \beta_{t_2}$ are either same or orthogonal for a factor state $\psi$ (see the above paragraph).
Interesting point here to note that the new representation $\pi'$ 
is $(\beta_t)$ co-variant i.e. $\pi' \beta_t=\beta_t \pi' $, hence by simplicity of the $C^*$ algebra $\clo_d$ we 
conclude that $$\pi'(\mbox{UHF}_d)'' = \pi'(\clo_d^{H})''^{\beta_t}$$
 
\vsp
By exploring the hypothesis that $\psi$ is a factor state, 
we also have as in Lemma 6.11 in [BJKW] $I \otimes L^{\infty}([0,{1 \over k})\;) \subset \pi'(\clo_d^H)''$. 
Hence we also have
$$\pi'(\clo_d^H)''= \pi(\clo_d^H)'' \otimes L^{\infty}([0,{1 \over k})\;).$$ 
Since $\beta_t$ is acting as translation on $I \otimes L^{\infty}([0,{1 \over k})\;)$ which being 
an ergodic action, we have 
$$\pi'(\mbox{UHF}_d)'' = \pi(\clo_d^{H})'' \otimes 1$$
Since $\pi'(\mbox{UHF}_d)'' = \pi(\mbox{UHF}_d)'' \otimes 1$, we conclude 
that $\pi(\mbox{UHF}_d)'' = \pi(\clo_d^{H})''$. 

\vsp
A proof for the statement (c) follows from Lemma 7.12 in [BJKW]. The original idea of the proof can be 
traced back to Arveson's work on spectrum of an automorphism of a commutative compact group [Ar1]. \qed

\vsp
Let $\omega'$ be an $\lambda$-invariant state on the $\mbox{UHF}_d$ sub-algebra of $\clo_d$. Following [BJKW, section 7], 
we consider the set 
$$K_{\omega'}= \{ \psi: \psi \mbox{ is a state on } \clo_d \mbox{ such that } \psi \lambda =
\psi \mbox{ and } \psi_{|\mbox{UHF}_d} = \omega' \}$$

By taking invariant mean on an extension of $\omega'$ to $\clo_d$, we verify that $K_{\omega'}$ is non empty and 
$K_{\omega'}$ is clearly convex and compact in the weak topology. In case $\omega'$ is an ergodic state ( extremal state )
$K_{\omega'}$ is a face in the $\lambda$ invariant states. Before we proceed to the next section here we recall Lemma 7.4 
of [BJKW] in the following proposition.

\begin{pro} Let $\omega'$ be ergodic. Then $\psi \in K_{\omega'}$ is an extremal point in
$K_{\omega'}$ if and only if $\psi$ is a factor state and moreover any other extremal point in $K_{\omega'}$
have the form $\psi \beta_z$ for some $z \in S^1$.
\end{pro}

\vsp
\NI {\bf PROOF:} Though Proposition 7.4 in [BJKW] appeared in a different set up, same proof goes through for
the present case. We omit the details and refer to the original work for a proof.  \qed

\section{ Dual Popescu system and pure translation invariant states: } 

\vsp
In this section we review the amalgamated Hilbert space developed in [BJKW] and prove a powerful criteria 
for a translation invariant factor state to be pure. 

\vsp
To that end let $\clm$ be a von-Neumann algebra acting on a Hilbert space $\clk$ and $\{v_k,\;1 \le k \le d\}$ 
be a family of bounded operators on $\clk$ so that $\clm=\{v_k,v_k^*,\;1 \le k \le d \}''$ and 
$\sum_kv_kv^*_k=1$. Furthermore let $\Omega$ be a cyclic and separating vector for $\clm$ so that the normal state 
$\phi_0(x)=<\Omega,x\Omega>$ on $\clm$ is invariant for the Markov map $\tau$ on $\clm$ defined by $\tau(x)=\sum_kv_kxv^*_k$ 
for $x \in \clm$. Let $\omega$ be the translation invariant state on UHF$_d =\otimes_{\IZ}M_d$ defined by
$$\omega(e^{i_1}_{j_1}(l) \otimes e^{i_2}_{j_2}(l+1) \otimes ....\otimes e^{i_n}_{j_n}(l+n-1)) = \phi_0(v_Iv^*_J)$$
where $e^i_j(l)$ is the elementary matrix at lattice sight $l \in \IZ$. 

\vsp
We set $\tilde{v}_k = \overline{ \clj \sigma_{i \over 2}(v^*_k) \clj } \in \clm'$ ( see [BJKW] for details ) 
where $\clj$ and $\sigma=(\sigma_t,\;t \in \IR)$ are Tomita's conjugation operator and modular automorphisms 
associated with $\phi_0$.

By KMS or modular relation [BR vol-1] we verify that 
$$\sum_k \tilde{v}_k \tilde{v}_k^*=1$$ 
and 
\be
\phi_0(v_Iv^*_J)= \phi_0(\tilde{v}_{\tilde{I}}\tilde{v}^*_{\tilde{J}})
\ee
where $\tilde{I}=(i_n,..,i_2,i_1)$ if $I=(i_1,i_2,...,i_n)$. Moreover $\tilde{v}^*_I\Omega = 
\clj \sigma_{i \over 2}(v_{\tilde{I}})^*\clj\Omega= \clj \Delta^{1 \over 2}v_{\tilde{I}}\Omega
=v^*_{\tilde{I}}\Omega$. We also set $\tilde{\clm}$ to be the von-Neumann algebra generated by 
$\{\tilde{v}_k: 1 \le k \le d \}$. Thus $\tilde{\clm} \subseteq \clm'$. A major problem that we 
will have to address when equality holds.  

\vsp
Let $(\clh,P,S_k,\;1 \le k \le d )$ and $(\tilde{\clh},P,\tilde{S}_k,\;1 \le k \le d)$ be the Popescu dilation described 
as in Theorem 2.1 associated with $(\clk,v_k,\;1 \le k \le d)$ and $\clk,\tilde{v}_k,\;1 \le k \le d)$ respectively. 
Following [BJKW] we consider the amalgamated tensor product $\clh \otimes_{\clk} \tilde{\clh}$ of $\clh$ with 
$\tilde{\clh}$ over the joint subspace $\clk$. It is the completion of the quotient of the set $$\!C \bar{I} \otimes 
\!C I \otimes \clk,$$ where $\bar{I},I$ both consist of all finite sequences with elements in $\{1,2, ..,d \}$, 
by the equivalence relation defined by a semi-inner product defined on the set by requiring
$$<\bar{I} \otimes I \otimes f,\bar{I}\bar{J} \otimes IJ \otimes g>=<f,\tilde{v}_{\bar{J}}v_Jg>,$$
$$<\bar{I}\bar{J} \otimes I \otimes f, \bar{I} \otimes IJ \otimes g> =<\tilde{v}_{\bar{J}}f,v_Jg>$$
and all inner product that are not of these form are zero. We also define
two commuting representations $(S_i)$ and $(\tilde{S}_i)$ of $\clo_d$ on
$\clh \otimes_{\clk} \tilde{\clh}$ by the following prescription:
$$S_I\lambda(\bar{J} \otimes J \otimes f)=\lambda(\bar{J} \otimes IJ \otimes f),$$
$$\tilde{S}_{\bar{I}}\lambda(\bar{J} \otimes J \otimes f)=\lambda(\bar{J}\bar{I} \otimes J \otimes f),$$
where $\lambda$ is the quotient map from the index set to the Hilbert space. Note that the subspace generated by
$\lambda(\emptyset \otimes I \otimes \clk)$ can be identified with $\clh$ and earlier $S_I$ can be identified
with the restriction of $S_I$ defined here. Same is valid for $\tilde{S}_{\bar{I}}$. The subspace $\clk$ is
identified here with $\lambda(\emptyset \otimes \emptyset \otimes \clk)$. 
Thus $\clk$ is a cyclic subspace for the representation $$\tilde{s}_j \otimes s_i \raro \tilde{S}_j S_i$$ 
of $\tilde{\clo}_d \otimes \clo_d$ in the amalgamated Hilbert space. Let $P$ be the projection on $\clk$. Then we
have 
$$S_i^*P=PS_i^*P=v_i^*$$
$$\tilde{S}_i^*P=P\tilde{S}_i^*P=\tilde{v}^*_i$$
for all $1 \le i \le d$. 

\vsp
We start with a simple proposition.

\begin{pro} 
The following holds: 

\NI (a) For any $1 \le i,j \le d$ and $|I|,|J|< \infty$ and $|\bar{I}|,|\bar{J}| < \infty$
$$<\Omega,\tilde{S}_{\bar{I}}\tilde{S}^*_{\bar{J}} S_iS_IS^*_JS^*_j \Omega>=<\Omega, 
\tilde{S}_i \tilde{S}_{\bar{I}}\tilde{S}^*_{\bar{J}}\tilde{S}^*_jS_IS^*_J \Omega>;$$

\NI (b) The vector state $\psi_{\Omega}$ on $$\tilde{\mbox{UHF}}_d \otimes \mbox{UHF}_d \equiv
\otimes_{-\infty}^0 M_d \otimes_1^{\infty}M_d \equiv \otimes_{\IZ} M_d$$ is equal to $\omega$;

\NI (c) $\pi(\tilde{\clo}_d \otimes \clo_d)''= \clb(\tilde{\clh} \otimes_{\clk} \clh)$ if and only if
$\{x \in \clb(\clk): \tau(x)=x,\; \tilde{\tau}(x)=x \}= \{zI: z \in  \!C \}$.

\end{pro} 

\vsp
\NI {\bf PROOF: } By our construction $\tilde{S}^*_i\Omega =\tilde{v}^*_i\Omega=v_i^*\Omega=S_i^*\Omega$. Now (a) and (b) 
follows by repeated application of $\tilde{S}^*_i\Omega=S^*_i\Omega$ and
commuting property of the two representation $\pi(\clo_d \otimes I)$ and $\pi(I \otimes \tilde{\clo}_d)$. The last
statement (c) follows from a more general fact proved below
that the commutant of $\pi(\clo_d \otimes \tilde{\clo}_d)''$ is
order isomorphic with the set $\{x \in \clb(\clk): \tau(x)=x,\; \tilde{\tau}(x)=x \}= \{zI: z \in  \!C \}$ via the 
map $X \raro PXP$ where $X$ is the weak$^*$ limit of $\{\Lambda^m\tilde{\Lambda}^n(x)$ as $(m,n) \raro 
(\infty,\infty)$. For details let $Y$ be the strong limit of increasing sequence of projections 
$(\Lambda \tilde{\Lambda})^n(P)$ as $n \raro \infty$. Then $Y$ commutes with $S_i\tilde{S}_j, S^*_i\tilde{S}^*_j$ 
for all $1 \le i,j \le d$. As $\Lambda(P)) \ge P$, we also have $\Lambda(Y) \ge Y$. Hence $(1-Y)S_i^*Y=0$. As $Y$ commutes 
with $S_i\tilde{S}_j$ we get $(1-Y)S_i^*S_i\tilde{S}_jY=0$ i.e. $(1-Y)\tilde{S}_jY=0$ for all $1 \le j \le d$. By symmetry 
of the argument we also get $(1-Y)S_iY=0$ for all $1 \le i \le d$. Hence $Y$ commutes with $\pi(\clo_d)''$ and by symmetry of 
the argument $Y$ commutes as well with $\pi(\tilde{\clo}_d)''$. As $Yf=f$ for all $f \in \clk$ and $\clk$ is cyclic for the 
representation $\pi(\tilde{\clo}_d \otimes \clo_d)$ we conclude that $Y=I$ in $\tilde{\clh} \otimes_{\clk} \clh$.

\vsp
Let $x \in \clb(\clk)$ so that $\tau(x)=x$ and $\tilde{\tau}(x) = x $ then as in the proof of Theorem 2.1 
we also check that $(\Lambda \tilde{\Lambda})^k(P) \Lambda^m \tilde{\Lambda}^n(x) (\Lambda \tilde{\Lambda})^k(P)$ is 
independent of $m,n$ as long as $m,n \ge k$. Hence weak$^*$ limit $\Lambda^m \tilde{\Lambda}^n(x) \raro X$ exists 
as $m,n \raro \infty$. Furthermore limiting element $X \in \pi(\clo_d \otimes \tilde{\clo}_d)'$ and $PXP=x$. That the 
map $X \raro PXP$ is an order-isomorphic on the set of self adjoint elements follows as in Theorem 2.1. This completes the proof. \qed

\vsp
Proposition 3.1 in brief says that $(\tilde{\clh} \otimes_{\clk} \clh,S_i\tilde{S}_j\;1 \le i,j \le d,P)$ 
is the Popescu dilation associated with Popescu elements $(\clk,v_i\tilde{v_j},1 \le i,j \le d \}$. Now we 
will be more specific in our starting Popescu's elements in order to explore the representation $\pi$ of 
$\tilde{\clo}_d \otimes \clo_d$ in the amalgamated Hilbert space $\tilde{\clh} \otimes_{\clk} \clh$.  

\vsp 
Let $\omega$ be a translation invariant factor state on $\clb$ and $\psi$ be an extremal point in $K_{\omega'}$. We consider the Popescu's elements $(\clk,\clm,v_k,1 \le k \le d,\Omega)$ described as in Proposition 2.4 associated 
with support projection of the state $\psi$ in $\pi_{\psi}(\clo_d)''$ and also consider associated dual Popescu's elements $(\clk,\tilde{\clm},\tilde{v_k},\;1 \le k \le d)$ where $\tilde{\clm}$ is the von-Neumann algebra generated by $\{ \tilde{v}_k:\; 1 \le k \le d \}$. Thus in general $\tilde{\clm} \subseteq \clm'$ and an interesting question: when do we have $\clm'=\tilde{\clm}$? Going back to our starting 
example of unique KMS state for the automorphisms $\beta_t(s_i)= ts_i,\;t \in S^1$, we check that $v^*_k=S^*_k$,  
$\clj \tilde{v}^*_k \clj = {1 \over d} S_k$ and thus equality holds i.e. $\tilde{\clm}=\clm'$. But 
the corner vector space $\tilde{\clm}_c=P\pi(\tilde{\clo}_d)''P$ generated by the elements 
$\{\tilde{v}_I\tilde{v}^*_J: |I|,|J| < \infty \}$ fails to be an algebra. Thus two questions sounds reasonable 
here.  

\vsp 
\NI (a) Does equality $\clm'=\tilde{\clm}$ holds in general for an extremal element $\psi \in K_{\omega'}$ and 
a factor state $\omega$?  

\vsp 
\NI (b) When can we expect $\tilde{\clm}_c$ to be a $*$-algebra and so equal to $\tilde{\clm}$?  

\vsp
The dual condition on support projection and equality $\tilde{\clm}=\clm'$ are rather deep and will lead us to a far reaching consequence on the state $\omega$. In the paper [BJKW] these two conditions are implicitly assumed to give a criteria for a translation invariant factor state to be pure. Apart from this refined interest, we will address the converse problem that terns out to be crucial for our main results. In the following we prove a crucial step towards that goal fixing the basic structure which will be repeatedly used in the computation using Cuntz relations. 

\vsp
\begin{pro} 
Let $\omega$ be a translation invariant factor state on $\clb$ and $\psi$ be an extremal point in $K_{\omega'}$. 
We consider the amalgamated representation $\pi$ of $\tilde{\clo}_d \otimes
\clo_d$ in $\tilde{\clh} \otimes_{\clk} \clh$ where the Popescu's elements $(\clk,\clm,v_k,\;1 \le k \le d)$ are 
taken as in Proposition 2.4. Then the following statements hold:

\vsp
\NI (a) $\pi(\tilde{\clo}_d \otimes \clo_d)''= \clb(\tilde{\clh} \otimes_{\clk} \clh)$. Furthermore $\pi(\clo_d)''$ 
and $\pi(\tilde{\clo}_d)''$ are factors and the following sets are equal:

\NI (i) $H=\{ z \in S^1: \psi \beta_z = \psi \}$;

\NI (ii) $H_{\pi}  = \{ z : \beta_z \mbox{ extends to an automorphisms of } \pi(\clo_d)'' \}$;

\NI (iii) $\tilde{H}_{\pi}= \{ z : \beta_z \mbox{ extends to an automorphisms of } \pi(\tilde{\clo}_d)'' \}$.
Moreover $\pi(\tilde{\mbox{UHF}}_d \otimes I)''$ and $\pi(I \otimes \mbox{UHF}_d)''$ are factors.

\vsp 
\NI (b) $z \raro U_z$ is the unitary representation of $H$ in the Hilbert space $\tilde{\clh} \otimes_{\clk} \clh$ 
defined by
$U_z(\pi(\tilde{s}_j \otimes s_i )\Omega=\pi(z\tilde{s}_j \otimes zs_i)\Omega$
 
\vsp  
\NI (c) The commutant of $\pi(\tilde{\mbox{UHF}}_d \otimes \mbox{UHF}_d )''$ is invariant by the canonical endomorphisms
$\Lambda(X)=\sum_i S_iXS_i^*$ and $\tilde{\Lambda}(X) = \sum_i \tilde{S}_iX \tilde{S}^*_i$. Same is true 
for each $i$ that the surjective map $X \raro S^*_iXS_i$ keeps the commutant of $\pi(\tilde{\mbox{UHF}}_d \otimes 
\mbox{UHF}_d)''$ invariant. Same holdss for the map $X \raro \tilde{S}^*_iX\tilde{S}_i$.

\vsp 
\NI (d) The centre of $\pi(\tilde{\mbox{UHF}}_d \otimes \mbox{UHF}_d )''$ is invariant by the canonical endomorphisms
$\Lambda(X)=\sum_i S_iXS_i^*$ and $\tilde{\Lambda}(X) = \sum_i \tilde{S}_iX \tilde{S}^*_i$. Moreover for each $i$ the
surjective map $X \raro S^*_iXS_i$ keeps the centre of $\pi(\tilde{\mbox{UHF}}_d \otimes 
\mbox{UHF}_d)''$ invariant. Same holdss for the map $X \raro \tilde{S}^*_iX\tilde{S}_i$.
\end{pro} 

\vsp
\NI {\bf PROOF: } $P$ being the support projection by Proposition 2.4 we have $\{x \in \clb(\clk): 
\sum_k v_k x v_k^*= x\} = \clm'$. That $(\clm',\tilde{\tau},\phi_0)$ is ergodic follows from a general 
result [Mo1] ( see also [BJKW] for a different proof ) as $(\clm,\tau,\phi_0)$ is ergodic for a factor state
$\psi$ being extremal in $K_{\omega'}$ (Proposition 2.6). Hence $\{x \in \clb(\clk): \tau(x)=\tilde{\tau}(x)=x \} = \!C$. Hence by Proposition 3.1 we conclude that $\pi(\tilde{\clo}_d \otimes \clo_d)''= \clb(\tilde{\clh} \otimes_{\clk} \clh)$. That both $\pi(\clo_d)''$ and $\pi(\tilde{\clo}_d)''$ are factors follows trivially as
$\pi(\tilde{\clo}_d \otimes \clo_d)''=\clb(\tilde{\clh} \otimes_{\clk} \clh)$ and $\pi(\clo_d)'' \subseteq \pi(\tilde{\clo}_d)'$. 

\vsp
By our discussion above we first recall that $\Omega$ is a cyclic vector for the representation of 
$\pi(\tilde{\clo}_d \otimes \clo_d)$. Let $G = \{ z=(z_1,z_2) \in S^1 \times S^1: \beta_z \mbox{
extends to an automorphism }$ $\mbox{on}\; \pi(\tilde{\clo}_d \otimes \clo_d)'' \}$ be the closed 
subgroup where
$$\beta_{(z_1,z_2)}(\tilde{s}_j \otimes s_i)=z_1 \tilde{s}_j \otimes z_2s_i .$$
By repeated application of the fact that $\pi(\clo_d)''$ commutes with $\pi(\tilde{\clo}_d)''$ and 
$S_i^*\Omega=\tilde{S}^*_i\Omega$ as in Proposition 3.1 (a) we verify that $\psi \beta_{(z,z)}=\psi$ 
on $\clo_d \otimes \tilde{\clo}_d$ if $z \in H$. For $z \in H $ we set unitary operator $U_z \pi(x \otimes y) \Omega 
= \pi(\beta_z(x) \otimes \beta_z(y))\Omega$ for all $x \in \tilde{\clo}_d$ and $y \in \clo_d$. 
Thus we have $U_z\pi(s_i)U_z^*=z\pi(s_i)$ and also $U_z\pi(\tilde{s}_i)U_z^*=z\tilde{s}_i.$  By taking 
it's restriction to $\pi(\clo_d)''$ and $\pi(\tilde{\clo}_d)''$ respectively we check that
$H \subseteq \tilde{H}_{\pi} $ and $H  \subseteq H_{\pi}$.

\vsp
For the converse let $z \in H_{\pi}$ and we use the same symbol $\beta_z$ for the extension to an automorphism of
$\pi(\clo_d)''$. By taking the inverse map we check easily that $\bar{z} \in H_{\pi}$ and in fact $H_{\pi}$ is
a subgroup of $S^1$. Since $\lambda$ commutes with $\beta_z$ on $\clo_d$, the canonical endomorphism $\Lambda$ defined
by $\Lambda(X) = \sum_k S_kXS_k^*$ also commutes with extension of $\beta_z$ on $\pi(\clo_d)''$. Note that the map
$\pi(x)_{{|}_{\clh}} \raro \pi(\beta_z(x))_{{|}_{\clh}} $ for $x \in \clo_d$ is a well defined linear $*$-homomorphism.
Since same is true for $\bar{z}$ and $\beta_z\beta_{\bar{z}}=I$, the map is an isomorphism. Hence $\beta_z$ extends
uniquely to an automorphism of $\pi(\clo_d)''_{{|}_{\clh}}$ commuting with the restriction of the canonical endomorphism
on $\pi(\clo_d)''_{|\clh}$. Since $\pi(\clo_d)''_{{|}_{\clh}}$ is a factor, we conclude as in Proposition 2.5 (a) that 
$z \in H$. Thus $H_{\pi} \subseteq H$. As $\pi(\tilde{\clo}_d)''$ is also a factor, we also have $\tilde{H}_{\pi} \subseteq H$. Hence we have $H=H_{\pi}=\tilde{H}_{\pi}$ and $\{(z,z): z \in H \} \subseteq G \subseteq H \times H$. 

\vsp
For the second part of (a) we will adopt the argument used for Proposition 2.5. To that end we first note that $\Omega$
being a cyclic vector for the representation $\tilde{\clo}_d \otimes \clo_d$ in the Hilbert space $\tilde{\clh} \otimes_{\clk} \clh $, by Lemma 7.11 in [BJKW] (note that the proof only needs the cyclic property ) the representation of UHF$_d$ on $\tilde{\clh} \otimes_{\clk} \clh $ is quasi-equivalent to it's sub-representation on the cyclic space generated by $\Omega$. On the other hand by our hypothesis that $\omega$ is a factor state, Power's theorem [Po1] ensures that the state $\omega'$
(i.e. the restriction of $\omega$ to $\clb_R$ which is identified here with $\mbox{UHF}_d$ ) is also a factor state
on $\mbox{UHF}_d$. Hence quasi-equivalence ensures that $\pi(I \otimes \mbox{UHF}_d)''$ is a factor. We also note that
the argument used in Lemma 7.11 in [BJKW] is symmetric i.e. same argument is also valid for $\tilde{\mbox{UHF}}_d$.
Thus $\pi(\tilde{\mbox{UHF}}_d \otimes I)''$ is also a factor. This completes the proof of (a). We have proved (b) 
while giving proof of (a).

\vsp
For $X \in \clb(\tilde{\clh} \otimes_{\clk} \clh)$, as $\Lambda(X)$ commutes with $\pi(\lambda(\tilde{\mbox{UHF}}_d \otimes \mbox{UHF}_d))''$ and
$\{ S_iS^*_j: 1 \le i,j \le d \}$ we verify by Cuntz's relation that $\Lambda(X)$ is also an 
element in the commutant of $\pi(\lambda(\tilde{\mbox{UHF}}_d \otimes \mbox{UHF}_d))''$ once 
$X$ is so. It is also obvious that $\Lambda(X)$ is an element in $\pi(\tilde{\mbox{UHF}}_d \otimes \mbox{UHF}_d)''$ if 
$X$ is so. Thus $\Lambda(X)$ is an element in the commutant/centre of $\pi(\lambda(\tilde{\mbox{UHF}}_d \otimes \mbox{UHF}_d)''$ once $X$ is so. For the last statement consider the map $X \raro S_i^*XS_i$
on $\pi(\tilde{\mbox{UHF}}_d \otimes \mbox{UHF}_d)''$ which is clearly onto by Cuntz relation. Hence 
we need to show that $S^*_iXS_i$ is an element in the commutant whenever $X$ is so. To that end note that $S_i^*XS_i S_i^*YS_i=S_i^*S_iS_i^*XYS_i=S_i^*YXS_i=S^*_iYS_iS_i^*XS_i$ since $X$ commutes with $S_iS_i^*$. 
Thus onto property of the map ensures that $S_i^*XS_i$ is an element in the commutant/centre of $\pi(\tilde{\mbox{UHF}}_d \otimes \mbox{UHF}_d)''$ once $X$ is so. This completes the proof of (c) and (d). \qed

\vsp
One interesting problem here how to describe the von-Neumann algebra $\cli$ consists of invariant elements of the gauge action $\{ \beta_z: z \in H \}$ in $\clb(\tilde{\clh} \otimes_{\clk} \clh)$. A general result 
due to E. Stormer [So] says that the algebra of invariant elements are von-Neumann algebra of type-I with 
centre completely atomic. Here the situation is much simple because we know explicitly 
that $\cli=\{U_z: z \in H \}'$ and we may write 
spectral decomposition as $U_z=\sum_{ k \in \hat{H} } z^k F_k$ for $z \in H$, $\hat{H}$ is the dual group of $H$ either $\hat{H}=\{z:z^n=1 \}$ or $\IZ$. Thus the centre of $\cli$ is equal to $\{F_k: k \in \hat{H} \}$. 

\vsp
As a first step we describe the center $\clz$ of $\pi(\tilde{\mbox{UHF}}_d \otimes \mbox{UHF}_d)''$ 
by exploring Cuntz relation that it is also non-atomic even for a factor state $\omega$. In fact we
will show that the centre $\clz$ is a sub-algebra of the centre of $\cli$. In the following proposition 
we give an explicit description.

\vsp
\begin{pro}  
Let $\omega,\psi$ be as in Proposition 3.2 with Popescu system $(\clk,\clm,$ $v_k,\Omega)$ be 
taken as in Proposition 2.4 i.e. on support projection. Then the centre of 
$\pi(\tilde{\mbox{UHF}}_d \otimes \mbox{UHF}_d)''$ is completely atomic and the element 
$E_0=[\pi(\tilde{\mbox{UHF}}_d \otimes \mbox{UHF}_d)' \vee \pi(\tilde{\mbox{UHF}}_d \otimes 
\mbox{UHF}_d)''\Omega]$ is a minimal projection in the centre of $\pi(\tilde{\mbox{UHF}}_d \otimes \mbox{UHF}_d)''$ and centre is invariant for both $\Lambda$ and $\tilde{\Lambda}$. Furthermore the 
following holds:

\vsp 
\NI (a) The centre of $\pi(\tilde{\mbox{UHF}}_d \otimes \mbox{UHF}_d)''$ has the following two disjoint
possibilities:

\vsp 
\NI (i) There exists a positive integer $m \ge 1$ such that the centre is generated by the family of minimal
orthogonal projections $\{ \Lambda_k(E_0): 0 \le k \le m-1 \}$ where $m \ge 1$ is the least positive integer so that
$\Lambda^m(E_0)=E_0$. In such a case $\{z: z^m = 1 \} \subseteq H$;

\vsp 
\NI (ii) The family of minimal nonzero orthogonal projections $\{ E_k: k \in \IZ \}$ where $E_k= \Lambda^k(E_0)$
for $k \ge 0$ and $E_k= S^*_IE_0S_I$ for $k < 0$ where $|I|=-k$ and independent of multi-index $I$ generates
the centre and $H=S^1$; 

\vsp
\NI (b) $\Lambda(E)=\tilde{\Lambda}(E) $ for any $E$ in the centre of $\pi(\tilde{\mbox{UHF}}_d \otimes \mbox{UHF}_d)''$

\vsp 
\NI (c) If $\Lambda(E_0)=E_0$ then $E_0=1$.

\end{pro} 

\vsp
\NI {\bf PROOF: } Let $E' \in \pi(\tilde{\mbox{UHF}}_d \otimes \mbox{UHF}_d)'$ be the projection on the subspace generated
by the vectors $\{ S_IS^*_J\tilde{S}_{I'}S^*_{J'}\Omega,\; |I|=|J|,|I'|=|J'| < \infty \}$ and $\pi_{\Omega}$ be the restriction
of the representation $\pi$ of $\tilde{\mbox{UHF}}_d \otimes \mbox{UHF}_d$ to the cyclic subspace $\clh_{\Omega}$ generated by
$\Omega$. Identifying $\clb$ with $\tilde{\mbox{UHF}}_d \otimes \mbox{UHF}_d$ we check that $\pi_{\omega}$ is unitary
equivalent with $\pi_{\Omega}$. Thus $\pi_{\Omega}$ is a factor representation. 

\vsp
For any projection $E$ in the centre of $\pi(\tilde{\mbox{UHF}}_d \otimes \mbox{UHF}_d)''$, via the unitary equivalence we note 
that $EE'=E'EE'$ is an element in the centre of $\pi_{\Omega}(\tilde{\mbox{UHF}}_d \otimes \mbox{UHF}_d)''$. $\omega$ being a 
factor state we conclude that $EE'$ is a scaler multiple of $E'$ and so we have
\be
EE'=\omega(E)E'
\ee
Thus we also have $EYE'=\omega(E)YE'$ for all $Y \in \pi(\tilde{\mbox{UHF}}_d \otimes \mbox{UHF}_d)'$ and so 
\be
EE_0=\omega(E)E_0
\ee
 
\vsp
Since $EE'$ is a projection and $E' \ne 0$, we have $\omega(E)=\omega(E)^2$. Thus $\omega(E)=1$ or $0$. So for such an element $E$
the following is true:

\NI (i) If $E \le E_0$ then either $E=0$ or $E=E_0$ i.e. $E_0$ is a minimal projection in the centre of $\pi(\tilde{\mbox{UHF}}_d \otimes \mbox{UHF}_d)''$

\NI (ii) $\omega(E)=1$ if and only if $E \ge E_0$ 

\NI (iii) $\omega(E)=0$ if and only if $EE_0=0$.

\vsp
As $\Lambda(E_0)$ is a projection in the centre of $\pi(\tilde{\mbox{UHF}}_d \otimes \mbox{UHF}_d)''$ by our last proposition i.e. Proposition 
3.2 (c), we have either $\omega(\Lambda(E_0))=1$ or $0$. Since $\Lambda(E_0) \neq 0$ by injective property of the endomorphism, 
we have either $\Lambda(E_0) \ge E_0$ or $\Lambda(E_0)E_0=0$. In case $\Lambda(E_0) \ge E_0$ we have $S^*_iE_0S_i \le S^*_i\Lambda(E_0)S_i=E_0$ 
for all $1 \le i \le d$. If so $S^*_iE_0S_i$ being a non-zero projection in the centre of $\pi(\mbox{UHF}_d \otimes \tilde{\mbox{UHF}}_d)''$ (Proposition 3.2 (c) ), 
by (i) we have $E_0=\Lambda(E_0)$. Thus we have either $\Lambda(E_0)=E_0$ or $\Lambda(E_0)E_0=0$. 

\vsp
If $\Lambda(E_0)E_0=0$, we have $\Lambda(E_0) \le I-E_0$ and by Cuntz's relation we check that
$E_0 \le I - S^*_iE_0S_i$ and $S^*_jS^*_iE_0S_iS_j \le I- S^*_jE_0S_j$ for all $1 \le i,j \le d$.
So we also have $E_0S^*_jS^*_iE_0S_iS_jE_0 \le E_0- E_0S^*_jE_0S_jE_0=E_0$. Thus we have either
$E_0S^*_jS^*_iE_0S_iS_jE_0= 0$ or $E_0S^*_jS^*_iE_0S_iS_jE_0=E_0$ as $S^*_jS^*_iE_0S_iS_j$ is an element
in the centre by Proposition 3.2 (c). So either we have $\Lambda^2(E_0)E_0=0$ or $\Lambda^2(E_0) \le E_0$. 
$\Lambda$ being an injective map we either have $\Lambda^2(E_0)E_0=0$ or $\Lambda^2(E_0)= E_0.$

\vsp
More generally we check that if $\Lambda(E_0)E_0=0, \Lambda^2(E_0)E_0=0,..\Lambda^k(E_0)E_0=0$ for some $k \ge 1$
then either $\Lambda^{k+1}(E_0)E_0=0$ or $\Lambda^{k+1}(E_0)=E_0$. To verify that first we check that in such a
case $E_0 \le I- S^*_IE_0S_I$ for all $|I|=n$ and then following the same steps as before to check that $S^*_iS^*_IE_0S_IS_i
\le I-S^*E_0S_i$ for all $i$. Thus we have $E_0S^*_iS^*_IE_0S_IS_iE_0 \le E_0$ and arguing as before we complete the proof
of the claim that either $\Lambda^{k+1}(E_0)E_0=0$ or $\Lambda^{k+1}(E_0)=E_0.$ 

\vsp
We summarize now by saying that $E_0,\Lambda(E_0),..,\Lambda^{m-1}(E_0)$ are mutually orthogonal projections with $m \ge 1$ possibly 
be infinite if not then $\Lambda^m(E_0)=E_0$. 

\vsp
Let $\pi_k,\;k \ge 0$ be the representation $\pi$ of $\tilde{\mbox{UHF}}_d \otimes \mbox{UHF}_d$ restricted to the
subspace $\Lambda^k(E_0)$. The representation $\pi_0$ of $\tilde{\mbox{UHF}}_d \otimes \mbox{UHF}_d$ is 
isomorphic to the representation $\pi$ of $\tilde{\mbox{UHF}}_d \otimes \mbox{UHF}_d$ restricted to $E'$ and thus
quasi-equivalent. For a general discussion on quasi-equivalence we refer to section 2.4.4 in [BR vol-1]. $\omega$ being a factor state, $\pi_0$ is a factor representation. We claim now that each $\pi_k$ is a factor representation. We fix any $k \ge 1$ and let $X$ be an element in the centre of $\pi_k(\mbox{UHF}_d \otimes \tilde{\mbox{UHF}}_d)$. Then for any $|I|=k$, $S^*_IE_kS_I=E_0$ and so $S^*_IXS_I$ is an element in the centre of $\pi_0(\mbox{UHF}_d \otimes \tilde{\mbox{UHF}}_d)$ by Proposition 3.2 (d). Further $S_I^*XS_I= S^*_IXS_IS^*_JS_J = S^*_JXS_J$ for all $|J|=|I|=k$. $\pi_0$ being a factor representation, we have $S^*_IXS_I= c E_0$ for some scaler $c$ independent of the multi-index we choose $|I|=k$. Hence $c \Lambda_k(E_0) = \sum _{|J|=k } S_JS^*_IXS_IS_J^*=\sum _{|J|=k } S_JS^*_IS_IS_J^*X= X$ as $X$ is an element in the centre of $\pi(\tilde{\mbox{UHF}}_d \otimes \mbox{UHF}_d)$. Thus for each $k \ge 1$, $\pi_k$ is a factor representation as $\pi_0$ is so.  

\vsp
We also note that $\Lambda(E_0)\tilde{\Lambda}(E_0) \neq 0$. Otherwise we have $<S_i\Omega,\tilde{S}_j\Omega>=0$
for all $i,j$ and so $<\Omega,\tilde{S}_jS^*_i\Omega>=0$ for all $i,j$ as $\pi(\clo_d)''$ commutes with 
$\pi(\tilde{\clo}_d)''$. However $\tilde{S}^*_i\Omega=S^*_i\Omega$ and $\sum_i \tilde{S}_i\tilde{S}^*_i=1$ 
which leads a contradiction. Hence $\Lambda(E_0)\tilde{\Lambda}(E_0) \neq 0$. As $\pi$ restricted to 
$\Lambda(E_0)$ is a factor state and both $\Lambda(E_0)$ and $\tilde{\Lambda}(E_0)$ are elements in the 
centre of $\pi(\tilde{\mbox{UHF}}_d \otimes \mbox{UHF})''$ by Proposition 3.2 (d), 
we conclude that $\Lambda(E_0)=\tilde{\Lambda}(E_0).$ Using commuting property of the endomorphisms $\Lambda$ 
and $\tilde{\Lambda}$, we verify by a simple induction method that $\Lambda^k(E_0) = \tilde{\Lambda}^k(E_0)$ 
for all $k \ge 1$. Thus the sequence of orthogonal projections $E_0, \tilde{\Lambda}(E_0),...$ are also periodic with same 
period or aperiodic according as the sequence of orthogonal projections $E_0,\Lambda(E_0),...$. 

\vsp
If $\Lambda^m(E_0)=E_0$ for some $m \ge 1$ then we check that $\sum_{0 \le k \le m-1} \Lambda^k(E_0)$ is a 
$\Lambda$ and as well $\tilde{\Lambda}$ invariant projection and thus equal to $1$ by cyclic property of 
$\Omega$ for $\pi(\clo_d \otimes \tilde{\clo}_d)''$. In such a case we set $V_z=\sum_{0 \le k \le m-1} z^kE_k$ for $z \in S^1$ 
for which $z^m=1$ and check that $\Lambda(V_z)= \sum_{0 \le k \le m-1}z^k\Lambda(E_k)=\sum_{0 \le k \le m-1} z^k E_{k+1}= \bar{z} V_z$ 
where $E_m=E_0$ and so by Cuntz relation we have $V^*_zS_iV_z= \bar{z}S_i$ for all $1 \le i \le d$. Following the same 
steps we also have $\tilde{\Lambda}(V_z)=\bar{z}V_z$ and so $V^*_z\tilde{S}_iV^*_z=\bar{z}\tilde{S}_i$ for $ 1 \le i \le d$. 
Thus $V_z=U_z$ for all $z \in H_0 = \{z: z^m=1 \} \subseteq H$. 

\vsp
Now we consider the case where $E_0,\Lambda(E_0),..\Lambda^k(E_0),..$ is a sequence of aperiodic orthogonal
projections. We extend family of projections $\{E_k:\;k \in \!Z \}$ to all integers by 
 
$$E_k=\Lambda^k(E_0)\;\; \mbox{for all}\;\; k \ge 1$$
and
$$E_k = S^*_IE_0S_I\;\; \mbox{for all}\;\;k \le 1,\;\;\mbox{where}\;\; |I|=-k$$
We claim that the definition of $\{ E_k;\;k \le -1 \}$ does depends only on length of the multi-index 
$I$ that we choose. We may choose any other $J$ so that $|J|=|I|$ and check the following identity: 
$$S^*_IE_0S_I = S^*_IE_0S_IS^*_JS_J=S^*_IS_IS^*_JE_0S_J=S^*_JE_0S_J$$ 
where $E_0$, being an element in the centre of $\pi(\tilde{\mbox{UHF}}_d \otimes \mbox{UHF}_d)''$, 
commutes with $S_IS_J^*$ as $|I|=|J|$. Further $\Lambda^k(E_0)=\tilde{\Lambda}^k(E_0)$ ensures that 
$S_I\tilde{S}_J^*$ commutes with $E_0$ for all $|I|=|J|=k$ and $k \ge 1$. Hence we also have 
$$E_{-k}=S^*_IE_0S_I\tilde{S}^*_J\tilde{S}_J=\tilde{S}^*_JE_0\tilde{S}_J$$ 
for all $|J|=|I|=k$ and $k \ge 1$. Now we claim that 
$$\Lambda(E_k) = \tilde{\Lambda}(E_k) = E_{k+1}$$ 
for all $k \in \!Z$. For $k \ge 0$ we have nothing to prove. For $k \le -1$ we check that the following steps
$$\Lambda(S^*_IE_0S_I)=\sum_k S_kS^*_iS^*_{I'}E_0S_{I'}S_iS_k^*$$
$$=\sum_kS^*_{I'}E_0S_{I'} S_kS^*_iS_iS_k^*=S^*_{I'}E_0S_{I'}$$
where we wrote $I=(I',i)$ and used elements $S_kS^*_i$ commutes with $\{ E_k: k \in \!Z \}$, elements in the centre of 
$\pi(\mbox{UHF}_d \otimes \tilde{\mbox{UHF}}_d)''$. For a proof that 
$\tilde{\Lambda}(E_k)=E_{k+1}$ we may follow the same steps as 
$E_k=\tilde{S}^*_IX\tilde{S}_I$ where $|I|=-k$ and $k \le -1$. 

\vsp
We also claim that $\{ E_k: k \in \!Z \}$ is an orthogonal family of non-zero projections. To that end we choose any two elements
say $E_k,E_m,\;k \ne m$ and use endomorphism $\Lambda^n$ for $n$ large enough so that both $n+k \ge 0 , n+m  \ge 0$
to conclude that $\Lambda^n(E_kE_m) = E_{k+n}E_{k+m}=0$ as $k+n \ne k+m$. $\Lambda$ being an injective map we get
the required orthogonal property. Thus $\sum_{ k \in \!Z} E_k$ being an invariant projection for both $\Lambda$ and
$\tilde{\Lambda}$ we get by cyclicity of $\Omega$ that $\sum_{ k \in \!Z} E_k=I$. Let $\pi_k,\; k \le -1$ be the representation $\pi$ 
of $\tilde{\mbox{UHF}}_d \otimes \mbox{UHF}_d$ restricted to the subspace $E_k$. Going along the same line as above, we verify that for 
each $k \le -1$, $\pi_k$ is a factor representation of $\tilde{\mbox{UHF}}_d \otimes \mbox{UHF}_d$. We also set $V_z= \sum_{-\infty < k 
< \infty} z^k E_k$ for all $z \in S^1$ and check that $\Lambda(V_z)=\bar{z}V_z$ and also $\tilde{\Lambda}(V_z)=\bar{z}V_z$. Hence $S^1=H$ 
and as $H$ is a closed subset of $S^1$. This completes the proof of (a). Proof for (b) and (c) are now simple consequence of the proof of (a). \qed

\vsp
It is clear that $\cli$ contains $\cli_0:=^{\mbox{def}} \pi(\tilde{\mbox{UHF}}_d \otimes \mbox{UHF}_d)'' \vee \{U_z: z \in H \}''$. 
By the last proposition the centre of $\cli$, which is equal to $\{U_z: z \in H \}''$, contains the centre of $\pi(\tilde{\mbox{UHF}}_d \otimes \mbox{UHF}_d)''$ 
and thus by taking commutant we also have $\cli \subseteq \pi(\tilde{\mbox{UHF}}_d \otimes \mbox{UHF}_d)''\vee \pi(\tilde{\mbox{UHF}}_d \otimes 
\mbox{UHF}_d)'$. In the last proposition we have described explicitly the factor decomposition of the representation 
$\pi$ of $\pi(\mbox{UHF}_d \otimes \tilde{\mbox{UHF}}_d)''$. One central issue when such an factor decomposition 
is also an extremal decomposition. A clear answer at this stage seems to be bit hard. However the following 
proposition makes an attempt for our purpose. To that end we set few more notations and elementary properties. 

\vsp
For each $k \in \hat{H}$, let $\pi'_k$ be the representation $\pi$ of $\tilde{\mbox{UHF}}_d \otimes \mbox{UHF}_d$ restricted to $F_k$. We claim that each $\pi'_k$ 
is pure once $\pi'_0$ is pure. Fix any $k \in \hat{H}$ and let $X$ be an element in the commutant of $\pi'_k(\tilde{\mbox{UHF}}_d \otimes \mbox{UHF}_d)''$ then 
$S^*_IF_kS_I = S^*_I\Lambda^k(F_0)S_I = F_0$ as $S^*_JS_I$ commutes with $F_0$ for $|I|=|J|=k$ and further for any $|I|=|J|$, $S_I^*XS_IS_J^*S_J=S_I^*S_IS_J^*XS_J=S^*_JXS_J$ 
as $X$ commutes with $S_IS_J^*$ with $|I|=|J|$. Thus by Proposition 3.2 (c) $S^*_IXS_I$ is an element in commutant of $\pi'_0(\tilde{\mbox{UHF}}_d \otimes \mbox{UHF}_d)''$ 
for any $|I|=k$ and thus $S_I^*XS_I= c F_0$ for some scaler $c$ independent of $|I|=k$ as $\pi'_0$ is pure. We use 
commuting property of $X$ with $\pi(\tilde{\mbox{UHF}}_d \otimes \mbox{UHF}_d)''$ to conclude that 
$X = c \Lambda^k(E_0)$ for some scaler $c$. If $k \le -1$ we employ the same method but with endomorphism 
$\Lambda^{-k}$ so that $\Lambda^{-k}(X)$ is an element in the commutant of 
$\pi'_0(\tilde{\mbox{UHF}}_d \otimes \mbox{UHF}_d)''$. Thus $\sum_{I:|I|=-k }S_IXS^*_I = c I$ and by 
injective property of the endomorphism we get $X$ is a scaler. Thus we conclude that each $\pi'_k$ is 
a pure once $\pi'_0$ is pure. 

\vsp 
Next we claim that for each fix $k \in \hat{H_0}$, representation $\pi_k$ of $\tilde{\mbox{UHF}}_d \otimes \mbox{UHF}_d$ defined in Proposition 3.3 is quasi-equivalent to representation $\pi'_k$ ( here we recall 
$\hat{H_0} \subseteq \hat{H}$ as $H_0 \subseteq H$ ). That $\pi'_0$ is quasi-equivalent to 
$\pi_0$ follows as they are isomorphic. More generally for any $k \in \hat{H}$, we abuse the notation and 
extend $\pi_k$ for the restriction of $\pi$ to the minimal central projections $E_k$ on the subspace span by 
$\{\pi(\tilde{\mbox{UHF}}_d \otimes \mbox{UHF}_d)'f: \forall \;\; f\;\; \in \clh \otimes_{\clk} \tilde{\clh}, \;\;F_kf=f \}$. Each $E_k$ is a minimal central element containing $F_k$ and however two such elements i.e. 
$E_k$ and $E_j$ are either equal or mutually orthogonal. Thus $\{E_k: k \in \hat{H} \} = \{E_k: k \in \hat{H}_0 \}$
and quasi-equivalence follows as $\pi_k$ is isomorphic with $\pi'_k$ for all $k \in \hat{H_0}$.     

\vsp 
At this stage we also set for the time being 
$$F'_0=[\pi(\tilde{\mbox{UHF}}_d \otimes \mbox{UHF}_d)''\Omega]$$ 
It is obvious that $F'_0 \le F_0$. We prove in following text that equality holds if $\omega$ is pure.
 
\vsp 
First we consider the case when $H=\{z: z^n=1\}$. Projections $\Lambda(F'_0)$ and $\tilde{\Lambda}(F'_0)$ are elements in $\pi(\tilde{\mbox{UHF}}_d \otimes \mbox{UHF}_d)'$ by Proposition 3.3. The representation $\pi(\tilde{\mbox{UHF}}_d \otimes \mbox{UHF}_d)''$ restricted to both the projections $\Lambda(F'_0),\tilde{\Lambda}(F'_0)$ are as well pure. A pull back by the map $X \raro S^*_iXS_i$ with any $1 \le i \le d$ will do the job for the projection $\Lambda(F'_0)$. Thus 
$\Lambda(F'_0)\tilde{\Lambda}(F'_0)\Lambda(F'_0) = c \Lambda(F'_0)$ for 
some scaler. By pulling back with the action $X \raro S^*_i X S_i$ we get
$F'_0 S^*_i\tilde{\Lambda}(F'_0)S_i F'_0 = cF'_0$ and so 
$$c=<\Omega,S^*_i\tilde{\Lambda}(F'_0)S_i\Omega>$$
$$=\sum_k <\Omega,S_k S^*_iF'_0S_iS^*_k\Omega>$$
as $\tilde{S}^*_k\Omega=S^*_k\Omega$ and further $F_0$ commutes 
with $\pi(\mbox{UHF}_d$ and thus $c=\sum_k <\Omega, S_k S_k^*\Omega>=1$.
This shows that $\tilde{\Lambda}(F'_0) \ge \Lambda(F'_0)$. Interchanging the role
of $\Lambda$ and $\tilde{\Lambda}$ we conclude that $\Lambda(F'_0)=\tilde{\Lambda}(F'_0)$. Proof essentially follows along the same line for 
$\Lambda^k(F'_0)=\tilde{\Lambda}^k(F'_0)$ for all $k \ge 1$. 
By Proposition 2.5 we also note that $\Lambda^n(F'_0)=F'_0=
\tilde{\Lambda}^n(F_0)$ as $H=\{z : z^n =1 \}.$ 

Thus $F'=\sum_{0 \le k \le n-1}\Lambda(F'_0)$ 
is a $\Lambda$ and as well $\tilde{\Lambda}$ invariant projection. Since $F'\Omega=\Omega$ we conclude by the 
cyclic property of $\Omega$ for $\pi(\clo_d \otimes \tilde{\clo}_d)''$ 
that $F'=1$. Since $\Lambda^k(F'_0) \le F_k$ and $\sum_k F_k=1 $ we conclude that 
$\Lambda^k(F_0)=F_k$. In such a case we may check that  
$$F_k = [\pi(\tilde{\mbox{UHF}}_d \otimes \mbox{UHF}_d)''S^*_I\Omega: |I|=n-k]$$ 
for $1 \le k \le n-1$. 

\vsp
Similarly in case $H=S^1$ and $\omega$ is pure we also have $F_0=F'_0$ and 
for $k \ge 1$ 
$$F_k=[\pi(\tilde{\mbox{UHF}}_d \otimes \mbox{UHF}_d)''S_I\Omega: |I|= k]$$
$$F_{-k}=[\pi(\tilde{\mbox{UHF}}_d \otimes \mbox{UHF}_d)''S^*_I\Omega: |I|=k]$$
Thus we have got an explicit description of the complete atomic centre of $\cli$ 
when $\omega$ is a pure state. 

\begin{pro}
Let $\omega,\psi$ and Popescu system $(\clk,\clm,v_k,\Omega)$ be as in Proposition 3.3. Then 

\NI (a) $\{\beta_z: z \in H \}$ invariant elements in $\pi(\tilde{\mbox{UHF}}_d \otimes \clo_d)''$ ( as well as in 
$\pi(\tilde{\clo}_d \otimes \mbox{UHF}_d)''$ ) are equal to $\pi(\mbox{UHF}_d \otimes \tilde{\mbox{UHF}}_d)''$. 

\NI (b) $\cli=\cli_0$ if and only if $\omega$ is pure. 

Further the following statements are equivalent: 

\NI (c) $\cli = \pi(\tilde{\mbox{UHF}}_d \otimes \mbox{UHF}_d)''$; 

\NI (d) $\pi(\tilde{\mbox{UHF}}_d \otimes \clo_d)''= \clb(\tilde{\clh} \otimes_{\clk} \clh)$;

\NI (e) $\pi(\tilde{\clo}_d \otimes \mbox{UHF}_d )''= \clb(\tilde{\clh} \otimes_{\clk} \clh)$;

\vsp
In such a case ( if any of (c),(d) and (e) is true ) the following statements are also true: 

\NI (f) $\pi(\tilde{\mbox{UHF}}_d \otimes \mbox{UHF}_d)''$ is a type-I von-Neumann algebra with centre equal 
to $\{ U_z: z \in H \}''$ where $U_z$ is defined in Proposition 3.2. 

\NI (g) $\omega$ is a pure state on $\clb$. 

\vsp
Conversely if $\omega$ is a pure state then $\pi(\tilde{\mbox{UHF}}_d \otimes \mbox{UHF}_d)''$ is a type-I von-Neumann 
algebra with centre equal to $\{U_z: z \in H_0 \}''$ where $H_0$ is a subgroup of $H$.     
\end{pro} 

\vsp
\NI {\bf PROOF: } Along the same line of the proof of Proposition 2.5 (b) we get $\{\beta_z: z \in H \}$ invariant elements 
in $\pi(\clo_d \otimes \tilde{\mbox{UHF}}_d)''$ is $\pi(\tilde{\mbox{UHF}}_d \otimes \mbox{UHF}_d)''$ where factor property 
of $\pi(\clo_d)''$ is crucial as in proof of Proposition 2.5 (b). Same holds for $\pi(\mbox{UHF}_d \otimes \tilde{\clo}_d)''$ as 
$\pi(\tilde{\clo}_d)''$ is a factor. Here we comment that factor property of $\pi(\tilde{\clo}_d)''$ can be
ensured whenever $\psi$ is an extremal element in $K_{\omega'}$ (See Proposition 3.2 (a) ). 

\vsp
For (b) we will first prove $\cli_0=\cli$ if $\omega$ is pure. As by definition $\cli_0 \subseteq \cli$, it is enough if we show 
$\cli'_0 \subseteq \cli'$. Let $X \in \cli'_0$ i.e. $X$ commutes with $\{ U_z:z \in H\}''$ and $\pi(\mbox{UHF}_d \otimes 
\tilde{\mbox{UHF}}_d)''$. For each $k \in \hat{H}$, $F_kXF_k$ is an element in the commutant of $F_k\pi(\mbox{UHF}_d \otimes 
\tilde{\mbox{UHF}}_d)''F_k$. $\omega$ being pure each representation $\pi$ restricted to $F_k$ is irreducible and thus
$F_kXF_k = c_k F_k$ for some scalers $c_k$. Hence $X=\sum_k c_k F_k \in \cli'=\{U_z: z \in H \}''$. 

\vsp
For the converse we need to show that the restriction of $\pi(\mbox{UHF}_d \otimes 
\tilde{\mbox{UHF}}_d)''$ to $F'_0$ is pure. Let $X$ be an element on the subspace $F'_0$ and in the commutant 
of $F'_0\pi(\mbox{UHF}_d \otimes \tilde{\mbox{UHF}}_d)''F'_0$, ( which in our earlier notation $E'$ in Proposition 3.3 ). 
Then $X$ commutes with each $F_k$ for $k \in \hat{H}$ and $\pi(\mbox{UHF}_d \otimes 
\tilde{\mbox{UHF}}_d)''F_k$ as $F'_0 \le F_0$ and $F_k$ are orthogonal to 
$F'_0$ for $k \neq 0$. So $X$ commutes with $\{U_z: z \in H \}''$ and 
$\pi(\mbox{UHF}_d \otimes \tilde{\mbox{UHF}}_d)''$ i.e. $X \in \cli'_0$. By our assumption $\cli_0=\cli$, we have now
$X \in \cli'$ which is equal to $\{F_k: k \in \hat{H} \}''$ and so $X=cF_0$ 
for some scaler $c_0$. This shows that $F'_0=F_0$ and $\omega$ is pure.    

\vsp
(c) implies (d): $\{U_z: z \in H \}$ is a commuting family of unitaries such that $\beta_z(X)=U_zXU_z^*$ and 
thus by (c) $\{U_z: z \in H \}'' \subseteq \pi(\tilde{\mbox{UHF}}_d \otimes \mbox{UHF}_d)''$. Let $X$ be an 
element in the commutant of $\pi(\tilde{\mbox{UHF}}_d \otimes \clo_d)''$. Then 
$X$ commutes also with $\{U_z: z \in H \}''$ and thus $X \in \pi(\tilde{\mbox{UHF}}_d \otimes \mbox{UHF}_d)''$ by (c).
Hence $X$ is an element in the centre of $\pi(\tilde{\mbox{UHF}}_d \otimes \mbox{UHF}_d)''$ and so 
$X=\sum_k c_kE_k$ where $E_k$ are the minimal projections in the centre of $\pi(\tilde{\mbox{UHF}}_d \otimes \mbox{UHF}_d)''$ 
given in Proposition 3.3. However $X$ also 
commutes with $\pi(\clo_d)''$ by our assumption (c) and $\Lambda(E_k)=E_{k+1}$ for $k \in \hat{H}$.  
So $c_k=c_{k+1}$ and $X$ is a scaler multiple of unit operator. Hence (d) follows from (c). Along the same line we prove 
(c) implies (e). For a proof for (d) implies (c) and (e) implies (c), we simply apply (a).   

\vsp
Now we will prove (f) and (g). That $\pi(\tilde{\mbox{UHF}}_d \otimes \mbox{UHF}_d)''$ is 
a type-I von-Neumann algebra ( with completely atomic centre ) follows by a theorem of [So] once we use (c). In the proof 
of Proposition 3.3 we have proved that the centre of $\pi(\tilde{\mbox{UHF}}_d \otimes \mbox{UHF}_d)''$ is 
$\{U_z: z \in H_0 \}''$ where $H_0 \subseteq H$. For equality in the present situation we simply use (c), as $\beta_w(U_z)=U_z$
for all $w,z \in H$, to conclude that $U_z$ is in the centre of $\pi(\tilde{\mbox{UHF}}_d \otimes \mbox{UHF}_d)''$. 

\vsp
If (c) holds then $\cli_0=\cli$ and thus (g) follows by (b). Here we will give another proof using the same idea to prove (f). 
Let $X$ be an element in the commutant of $\pi_0(\tilde{\mbox{UHF}}_d \otimes \mbox{UHF}_d)''$, where $\pi_0$ is the factor 
representation on the minimal central projection $E_0$ defined in Proposition 3.3. Then $X$ commutes with 
$\{U_z: z \in H \}''$ and so by (c) $X$ in an element in $\pi(\tilde{\mbox{UHF}}_d \otimes \mbox{UHF}_d)''$. 
So $X$ is in the centre of $\pi_0(\tilde{\mbox{UHF}}_d \otimes \mbox{UHF}_d)''$. $\pi_0$ being a factor representation
$X$ is a scaler multiple of $E_0$. Thus $\pi_0$ is an irreducible representation and so $\omega$ is pure. 

\vsp
By Proposition 3.1 we recall that $\pi'_0$ is unitarily equivalent to GNS representation of $(\clb,\omega)$. 
Thus $\pi'_0$ is irreducible if and only if $\omega$ is pure. So for a pure state $\omega$, 
for each $k \in \hat{H}_0$, $\pi_k$ being quasi-equivalent to $\pi'_k$, $\pi_k$ is a type-I factor 
representation of $\pi(\mbox{UHF}_d \otimes \tilde{\mbox{UHF}}_d)''$. This completes the proof. \qed
  
\vsp
The following theorem is the central step that will be used repeatedly.

\vsp
\begin{pro} 
Let $\omega$ be an extremal translation invariant state on $\clb$ and $\psi$ be an extremal 
element $\psi$ in $K_{\omega}$. We consider the Popescu elements $(\clk,v_k:1 \le k \le d, \clm,\Omega)$
as in Proposition 2.4 for the dual Popescu elements and associated amalgamated representation $\pi$ of
$\clo_d \otimes \tilde{\clo}_d$ as described in Proposition 3.1. Then the following holds:

\vsp 
\NI (a) $\pi(\tilde{\clo}_d \otimes \clo_d)''= \clb( \tilde{\clh} \otimes_{\clk} \clh)$;

\vsp 
\NI (b) $\pi(\tilde{\clo}_d)''=\pi(\clo_d)'$ if and only if $\pi(\tilde{\clo}_d)''E=\pi(\clo_d)'E$. 

\vsp 
\NI (c) $Q=\cle \tilde{\cle}$ is the support projection of the state $\psi$ in $\pi(\clo_d)''\tilde{\cle}$
and also in $\pi(\tilde{\clo}_d)''E$ where $\cle$ and $\tilde{\cle}$ are the support projections of the state $\psi$
in $\pi(\clo_d)''$ and $\pi(\tilde{\clo}_d)''$ respectively; 

\vsp 
\NI (d) If $\cle \clf = \tilde{\cle}\tilde{\clf}$ then $\cle=\tilde{\clf},\;\tilde{\cle}=\clf,\;P=Q$. 

\vsp 
\NI (e) If $P=Q$ then the following statements are true:

\vsp 
\NI (i) $\clm'= \tilde{\clm}$ where $\tilde{\clm} = \{ P\tilde{S}_iP: 1 \le i \le d
\}''$;

\vsp 
\NI (ii) $\pi(\clo_d)' = \pi(\tilde{\clo}_d)''$.

\vsp 
\NI (f) If $P=[\tilde{\clm}\Omega]$ then $\clm'=\tilde{\clm}$. 

\vsp 
\NI (g) $\omega$ is pure on $\clb$ if and only if there exists a sequence of elements $x_n \in \clm$ such that for each $m \ge 0$ 
$x_{n+m}\tau_n(x) \raro \phi_0(x)1$ as $n \raro \infty$ in strong operator topology, equivalently $\phi_0(\tau_n(x)x_{n+m}^*x_{n+m}\tau_n(y)) \raro \phi_0(x)\phi_0(y)$ as $n \raro \infty$ for 
all $x,y \in \clm$ where $\clm = \{ v_i=PS_iP: 1 \le i \le d \}''$ and $\tau(x)=\sum_{1 \le k \le d}v_kxv_k^*,\;\;x \in \clm$; Same holds true if we replace $\clm_0$ for 
$\clm$ where $\clm_0=\{x \in \clm: \beta_z(x)=x; z \in H \}$.  

\end{pro} 

\vsp
\NI {\bf PROOF: } (a) is a restatement of Proposition 3.2 (a). $\cle$ ( $\tilde{\cle}$ ) being the support projection of
the state $\psi$ in $\pi(\clo_d)''$ ( $\pi(\tilde{\clo}_d)''$ ) and $\psi = \psi \Lambda$ we have $\Lambda(\cle) \ge \cle$ 
and further we have $\cle=[\pi(\clo_d)'\Omega] \ge [\pi(\tilde{\clo}_d)''\Omega]$ and hence $\Lambda^n(\cle) \uparrow I$ 
as $n \raro \infty$ because $\Omega$ is cyclic for $\pi(\clo_d \otimes \tilde{\clo}_d)''$ in $\clh \otimes_{\clk} \tilde{\clh}$. 

\vsp 
We set von-Neumann algebras $\cln_1=\pi(\clo_d)'\cle$ and $\cln_2=\pi(\tilde{\clo}_d)''\cle$. By our construction in general $\pi(\tilde{\clo}_d)'' \subseteq \pi(\clo_d)' $ and so $\cln_2 \subseteq \cln_1$. Since $\Lambda^n(\cle) \uparrow I$ as $n \raro \infty$ in strong operator topology, two operators in $\pi(\clo_d)'$ are same if their actions are same on $\cle$.  
So (b) is true.

For (c) we note that $Q=\cle\tilde{\cle} \in \cln_2 \subseteq \cln_1$ and claim that $Q$ is the support projection of the state $\psi$ in $\cln_2$. To that end let $x\cle \ge 0$ for some $x \in \pi(\tilde{\clo}_d)''$ so that $\psi(QxQ)=0$. 
As $\Lambda^k(x\cle) \ge 0$ for all $k \ge 1$ and $\Lambda^k(\cle) \raro I$ we conclude that $x \ge 0$. As $\cle \Omega=\Omega$ and thus $\psi(\tilde{\cle}x\tilde{\cle})=\psi(QxQ)=0$, we conclude
$\tilde{\cle}x\tilde{\cle}=0$, $\tilde{\cle}$ being the support projection for $\pi(\tilde{\clo}_d)''$. Hence $QxQ=0$.
As $\psi(Q)=1$, we complete the proof of the claim that $Q$ is the support of $\psi$ in $\cln_2$. Similarly $Q$ is also
the support projection of the state $\psi$ in $\pi(\clo_d)''\tilde{\cle}$. This completes the proof of (c).

\vsp 
Thus if $\cle \clf=\tilde{\cle}\tilde{\clf}$, we get $\Lambda^n(\cle)\clf=\tilde{\cle} \Lambda^n(\tilde{\clf})$ and 
$\cle \tilde{\Lambda}(\clf) = \tilde{\Lambda}(\tilde{\cle})\tilde{\clf}$ and thus taking limit we get 
$\clf=\tilde{\cle}$ and $\cle=\tilde{\clf}$. It is obvious now that $P=\cle \clf = \cle \tilde{\cle}=Q$. This completes 
the proof of (d).

\vsp
As $\cle \in \pi(\clo_d)''$ and $\tilde{E} \in \pi(\tilde{\clo}_d)''$ we check that von-Neumann algebras 
$\clm^1 = Q \pi(\clo_d)'' Q$ and $\tilde{\clm}^1 = Q \pi(\tilde{\clo}_d) Q$ acting on $Q$ satisfies 
$\tilde{\clm}^1 \subseteq \clm^{1'}$. Now we explore that $\pi(\tilde{\clo}_d \otimes \clo_d)''= 
\clb(\clh \otimes_{\clk} \tilde{\clh})$ and note that in such a case $Q\pi( \tilde{\clo}_d \otimes \clo_d )''Q$ 
is the set of all bounded operators on the Hilbert subspace $Q$. As $\cle \in \pi(\clo_d)''$ and 
$\tilde{\cle} \in \pi(\tilde{\clo}_d)''$ we check that together $\clm^1 = Q \pi(\clo_d)'' Q$ and 
$\tilde{\clm}^1 = Q \pi(\tilde{\clo}_d) Q$ generate all bounded operators on $Q$. Thus both $\clm^1$ 
and $\tilde{\clm}^1$ are factors. The canonical states $\psi$ on $\clm^1$ and $\tilde{\clm}^1$ are 
faithful and normal. We set $l_k=QS_kQ$ and $\tilde{l}_k=Q\tilde{S}_kQ,\; 1 \le k \le d$ and recall 
that $v_k=PS_kP$ and $\tilde{v}_k=P\tilde{S}_kP,\; 1 \le k \le d$. We note that $Pl_kP=v_k$ and 
$P\tilde{l}_kP = \tilde{v}_k$ where we recall by our construction $P$ is the support projection of 
the state $\psi$ in $\pi(\clo_d)''_{|}[\pi(\clo_d)\Omega]$. $Q$ being the support projection of 
$\pi(\clo_d)\tilde{\cle}$, by Theorem 2.4 applied to Cuntz elements $\{S_i\tilde{E}: 1 \le i \le d \}$, 
$\tilde{\cle} \pi(\clo_d)'\tilde{\cle}$ is order isomorphic to $\clm^{1'}$ via the map $X \raro QXQ$. As the projection 
$F=[\pi(\clo_d)''\Omega] \in \pi(\clo_d)'$, we check that the element $Q\clf\tilde{\cle}Q \in \clm^{1'}$. However 
$Q\clf \tilde{\cle}Q=\cle\tilde{\cle}\clf\tilde{\cle}\cle=QPQ=P$ and thus $P \in \clm^{1'}$. We also check that 
$\clm^1\Omega=\clm^1P\Omega=P\clm^1\Omega=\clm \Omega$ and thus $P=[\clm^1\Omega]$. We set 
$\tilde{\clm}$ for the von-Neumann algebra generated by $\{\tilde{v}_k:\; 1 \le k \le d \}$. 

\vsp 
In such a case $\clm^1=\clm$ and $\tilde{\clm}^1=\tilde{\clm}$. By order isomorphic property 
we get (i) is equivalent to $\tilde{\cle} \pi(\clo_d)' \tilde{\cle} = \tilde{\cle} \pi(\tilde{\clo}_d)''\tilde{\cle}$ 
and taking commutant again we get $\pi(\clo_d)''\tilde{\cle} = \pi(\tilde{\clo}_d)'\tilde{\cle}$. Now we invoke the first 
part of the argument changing the role or using the endomorphism $\tilde{\Lambda}$ we conclude that 
$\pi(\clo_d)''=\pi(\tilde{\clo}_d)'$. This completes the proof of (e) provided we find a proof for (i) 
which is not so evident. 

\vsp 
Now we explore the representation $\pi$ of $\tilde{\clo}_d \otimes \clo_d$ which is pure to prove (i). To that end we note since $P=Q$ by
our assumption, $\Omega$ is a common cyclic and separating vector for $\tilde{\clm}$ and $\clm'$.   
Thus we can get an endomorphism $\alpha: \clm' \raro \tilde{\clm}$ defined by $\alpha(y)=\clj\tilde{\clj}y\tilde{\clj}\clj$ 
where $\tilde{\clj}$ is the Tomita's conjugate operator associated with cyclic and separating vector $\Omega$ for $\tilde{\clm}$. 
We note that the general theory does not guarantee [AcC] that the endomorphism be Takesaki's 
canonical conditional expectation associated with $\phi_0$. If so then the modular 
automorphism group $(\sigma_t)$ of $\clm'$ also preserves $\tilde{\clm}$. Thus $\sigma_z(x) \in \tilde{\clm}$ for $-1 \le Im(z) \le 1$ 
if $x$ is an analytic element in $\tilde{\clm}$. Thus we would have got $\clj v_k(\delta) \clj = 
\sigma_{i \over 2}(\tilde{v}_k(\delta)) \in \tilde{\clm}$ where $x(\delta)$ is average of 
$\sigma_t(x)$ with respect to Gaussian measure with variance $\delta > 0$. That $\tilde{v}_k(\delta)$ 
is an analytic element follows from the general Tomita-Takesaki theory [BR1]. Since $v_k(\delta) \raro v_k$ in strong 
operator topology as $\delta \raro 0$ and $\clm = \{v_k: 1 \le k \le d \}''$ together with 
$\clj \clm \clj = \clm'$ we arrived at $\tilde{\clm}=\clm'$. In the following we avoid this tempted 
route and aim to explore the general representation theory of $C^*$-algebras [BR1,chapter 2].   

\vsp 
We claim that $\clm'=\tilde{\clm}$. Suppose not. Then $\alpha(\tilde{\clm})$ is a proper von-Neumann subalgebras 
of $\alpha(\clm') \subseteq \tilde{\clm} $ being an into map and hence $\alpha(\tilde{\clm})$ is a proper von-Neumann 
subalgebra of $\tilde{\clm}$. Now consider the Popescu elements $(\clk,\alpha(\tilde{v}_i),\Omega)$ and its dilation 
as in Theorem 2.1. Then by the commutant lifting theorem applied to pairs 
$(\tilde{v}_i),\alpha(\tilde{v}_i)$ we find an unitaryoperator $U$ on $\tilde{\clh}$ so that 
$U\pi(\tilde{\clo}_d)''U^*$ is strictly contained in $\pi(\tilde{\clo}_d)''$ ( Without loss of generality we can 
take the dilated Hilbert space for $(\clk,\alpha(\tilde{v}_i),\Omega)$ to be same as $\tilde{\clh}$ as there exists an 
isomorphism preserving $\clk$, see the remark that follows after Theorem 2.1 ). We extend $U$ to an unitary operator 
on $\tilde{\clh} \otimes_{\clk} \clh$ and denote $\pi_u(x)=U \pi(x) U^*$ for $x \in \tilde{\clo}_d \otimes \clo_d$ which 
is unitary equivalent to the pure representation $\pi$ and $\pi_u(\tilde{\clo}_d)''$ is strictly contained in 
$\pi(\tilde{\clo}_d)''$. Now $\pi_u$ is also an amalgamated representation over the subspace $\clk$ with 
$P_u=Q_u$. Thus we can repeat now same with $\pi_u$ and so on. Note that the process won't terminate in finite 
time. Our aim is to bring a contradiction from this using formal set theory.  

\vsp 
To that end we reset for $\pi$ as $\pi_0$ as temporary notation as $\pi$ will be used as notation for a generic representation. 
Let $\clp$ be the collection of representation $(\pi,H_{\pi},\Omega)$ quasi-equivalent to $\pi_0:\tilde{\clo}_d \otimes \clo_d  
\raro \clb(\tilde{\clh} \otimes_{\clk} \clh )$ with a shift invariant vector state $\omega(x)=<\Omega,\pi(x)\Omega>$ i.e. $\omega(\pi(\theta(x))=\omega(\pi(x))$. So there exists cardinal numbers $n_{\pi},n_0(\pi)$ so that $n_{\pi}H_{\pi}$ is unitary equivalent to $n_0(\pi)\pi_0$. Thus given an element $(\pi,H_{\pi},\Omega)$ we can associate two cardinal numbers $n_{\pi}$ and $n_0(\pi)$ and without loss of generality we assume that $H_{\pi} \subseteq n_0(\pi)H_0$ and $n_{\pi}H_{\pi}=n_0(\pi)H_0.$  $\pi_0$ being a pure representation, any element $\pi \in \clp$ is a type-I factor representation 
of $\tilde{\clo}_d \otimes \clo_d $. The interesting point here that $\dsp{\oplus_{\pi \in \clp }}\pi$ 
is also an element in $\clp$ with associated cardinal numbers $\sum_{\pi} n_{\pi}$ and $\sum_{\pi} n_0(\pi)$. We say $(\pi_1,H_1,\Omega^1) << (\pi_2,H_2,\Omega^2)$ if there exists an isometry $U:n_{\pi_1} H_1 \raro n_{\pi_2} H_2$ so that 

\vsp 
\NI (C1) For each $1 \le \alpha \le n_{\pi_1} $ we have 
$U \Omega^1_{\alpha} = \Omega^2_{\alpha'}$ for some 
$1 \le \alpha' \le n_{\pi_2}$;

\NI (C2) $n_{\pi_2}\pi_2(x)E'_2 = U n_{\pi_1}\pi_1(x) U^* $ where $\cle'_2 \in n_{\pi_2}\pi_2(\tilde{\clo}_d)'$; 
 
\NI (C3) $U\oplus_{1 \le \alpha \le n_{\pi_1} }\pi^{\alpha}_1(\tilde{\clo}_d)''U^* \subset \oplus_{1 \le 
\alpha \le n_{\pi_2}}\pi^{\alpha}_1(\tilde{\clo}_d)''E'_2$. 

\vsp 
That the partial order is non-reflexive follows as $(\pi,H,\Omega) << (\pi,H,\Omega)$ contradicts (C3) as $I=E_2'$. By our starting assumption that $\clm' \neq \tilde{\clm}$ we check that $\pi_0 << \pi_u$. Thus going via the isomorphism we also check that for a given element $\pi \in \clp$ there exists an element $\pi' \in \clp$ so that $\pi << \pi'$. Thus $\clp_0$ is a non empty set and has at least one infinite chain. Partial order property follows easily. If $\pi_1 << \pi_2$ and $\pi_2 << \pi_3$ then $\pi_1 << \pi_3$. If $U_{12}$ and $U_{23}$ are isometric operators that satisfies (C1)-(C3) respectively, then $U_{13}=U_{23}U_{12}$ will do the job for $\pi_1$ and $\pi_3$. 
 
\vsp 
However by Hausdorff maximality theorem there exists a non-empty maximal totally ordered subset $\clp_0$ of $\clp$. We claim that 
$\pi_{max}=\oplus_{\pi \in \clp_0} \pi$ on $H_{\pi_{max}}=\oplus_{\pi \in \clp_0} H_{\pi} $ is an upper bound 
in $\clp_0$. That $\pi_{max} \in \clp$ is obvious. Further given an element $(H_1,\pi_1,\Omega_1) \in \clp_0$ there exists 
an element $(H_2,\pi_2,\Omega_2) \in \clp_0$ so that $\pi_1 << \pi_2$ by our starting remark as $\pi_0 << \pi_u$. 
By extending isometry $U_{12}$ to an isometry from $H_1 \raro n_{\pi_{max}} H_{\pi_{max}}$ trivially we get the 
required isometry that satisfies (C1),(C2) and (C3) where cardinal numbers $n_{\pi_{max}} = \sum_{\pi \in \clp_0} 
n_{\pi} \in \aleph_0$. Thus by maximal property of $\clp_0$ we have $\pi_{max} \in \clp_0$. This brings a contradiction 
as by our construction $(\pi_{max},H_{\pi_{max}},\Omega) << (\pi_{max},H_{\pi_{max}},\Omega)$ 
as $\pi_{max} \in \clp_0$ but partial order is strict. This contradicts our starting hypothesis that $\tilde{\clm}$ is 
a proper subset of $\clm'$. This completes the proof for (i) of (e) $\clm'=\tilde{\clm}$ when $P=Q$.  

\vsp 
In the proof of $\clm'=\tilde{\clm}$ in (e), we have used equality $P=Q$ just to ensure that $\Omega$ is also a cyclic for $\tilde{\clm}$ and $P=Q$ is used to prove $\pi(\clo_d)' = \pi(\tilde{\clo}_d)''$. So (f) follows by the proof of (e). 

\vsp 
A proof for (g) is given in [Mo3] with $\clm_0$. Here we will also give an alternative proof relating the criteria obtained in Proposition 3.4. 
To that end we claim that $$\bigcap_{n \ge 1}\tilde{\Lambda}^n(\pi(\mbox{UHF}_d)')=\pi(\tilde{\mbox{UHF}}_d)'
\bigcap \pi(\mbox{UHF}_d)'.$$ That $\tilde{\Lambda}^n(\pi(\mbox{UHF}_d)') \subseteq \{ \tilde{S}_I\tilde{S}^*_J: |I|=|J| < \infty \}'$ follows by 
Cuntz relation and thus $\bigcap_{n \ge 1}\tilde{\Lambda}^n(\pi(\mbox{UHF}_d)') \subseteq \pi(\tilde{\mbox{UHF}}_d)'\bigcap 
\pi(\mbox{UHF}_d)'$. For the reverse inclusion let $X \in E\pi(\tilde{\mbox{UHF}}_d)'\bigcap \pi(\mbox{UHF}_d)'E$. For $n \ge 1$, 
we choose $|I|=n$ and set $Y_n=\tilde{S}^*_IX\tilde{S}_I$.  We check that it is independent of the index that we have chosen as 
$Y_n=\tilde{S}^*_IX\tilde{S}_I\tilde{S}^*_J\tilde{S}_J=\tilde{S}^*_I\tilde{S}_I\tilde{S}^*_JX \tilde{S}_J=\tilde{S}^*_JX\tilde{S}_J$ 
where in second equality we have used $X \in \pi(\tilde{\mbox{UHF}}_d)'$ and also 
$\Lambda^n(Y_n)=\sum_{|J|=n } \tilde{S}_J\tilde{S}^*_IX\tilde{S}_I\tilde{S}^*_J=X$. This proves the equality 
in the claim. Going along the same line we also get $$\bigcap_{n \ge 1}\tilde{\Lambda}^n(\pi(\clo_d)') = 
\pi(\tilde{\mbox{UHF}}_d)'\bigcap \pi(\clo_d)'= \pi(\tilde{\mbox{UHF}}_d \otimes \clo_d)'.$$ 

\vsp 
By Proposition 3.4 $\omega$ is pure if and only if the set above is trivial. Thus once more by Proposition 1.1 in [Ar2] and 
Theorem 2.4 in [Mo2], purity is equivalent to asymptotic relation $||\Psi \tilde{\tau}^n -\phi_0|| \raro 0$ as $n \raro \infty$ 
for any normal state on $\clm'$ ( Here we recall by Proposition 2.4 $P\pi(\clo_d)'P=\clm'$ as $P$ is also the support projection 
in $\pi(\clo_d)''F$ ), where commutant is taken in $\clb(\clk)$. By duality argument [Mo3] we conclude that $\omega$ is pure if 
and only if there exists a sequence of elements $x_n \in \clm$ so that for each $m \ge 0$, $x_{m+n}\tau_n(x) \raro \phi_0(x)1$ as $n \raro \infty$ 
for all $x \in \clm \subseteq \clb(\clk)$. This completes the proof of (d) with $\clm$. For the proof with $\clm_0$ we need to 
show if part as only if part follows $\clm_0$ being a subset of $\clm$ and $\tau$ 
takes elements of $\clm_0$ to itself. For if part we refer to Theorem 3.2 in [Mo3]. \qed   

\vsp
We set 
$$(\clm')_0 = \{x \in \clm' : \beta_z(x)=x,\;\; z \in H \}.$$
Similarly we also set $\tilde{\clm}_0$ and $(\tilde{\clm}')_0$ as $(\beta_z:\;z \in H)$ invariant elements
of $\tilde{\clm}$ and $(\tilde{\clm}')$ respectively. We note that as a set $(\tilde{\clm}_0)'$ could be different from 
$(\tilde{\clm}')_0$. We note also that $P\tilde{\clm}^1P \subseteq \tilde{\clm}$ and unless $P$ is an element 
in $\tilde{\clm}^1$, equality is not guaranteed for a factor state $\omega$. The major problem is to show that 
$P$ is indeed an element in $\tilde{\clm}^1$ when $\omega$ is a pure state.   

\vsp
We warn here an attentive reader that in general for a factor state $\omega$, the set $\clf \pi(\tilde{\clo}_d)'' \clf$, which is a 
subset of $\clf \pi(\clo_d)'\clf$, need not be an algebra. However by commutant lifting theorem applied to dilation $v_i \raro S_i\clf$, 
$\pi(\clo_d)'\clf$ is order isomorphic to $\clm'$ as $P=\clf \cle$ is the support projection. Thus the von-Neumann sub-algebra generated 
by the elements $\clf\pi(\tilde{\clo}_d)''\clf$ is order isomorphic to $\tilde{\clm}$. However 
$\tilde{\clm}_0$ may properly include $\tilde{\clm}_{00}=\{ P\pi(\tilde{\mbox{UHF}}_d)P \}''$ 
( as an example take $\psi$ to be the unique KMS state on $\clo_d$ and $\omega$ be the unique trace on $\clb$ 
for which we get $\tilde{\clm}_{00} = \!C$ and $P\pi(\tilde{\clo}_d)''P$ is the linear span of 
$\{\tilde{v}_J^*,\;I,\;\tilde{v}_J:|J| < \infty \}$. 

\vsp 
Existence of a $\phi_0$ preserving norm one projection $\int_{z \in H} \beta_z dz $ 
ensures that modular operator of $\phi_0$ preserves $\clm_0$ [Ta] and so does on 
$(\clm')_0$. However there is no reason to take it granted for $\tilde{\clm}_0$ to be 
invariant by the modular group of $((\clm')_0,\phi_0)$. By Takesaki's theorem such a property 
is true if and only if there exists a $\phi_0$-invariant norm one projection from $(\clm')_0$ 
onto $\tilde{\clm}_0$. In the following we avoid this tempted route.  

\vsp 
At this stage it is not clear how we can ensure existence of a norm one projection from $\clm'$ to 
$\tilde{\clm}$ directly and so the equality $\clm'=\tilde{\clm}$ when $\omega$ is a pure state. Further interesting point here that the equality $\tilde{\clm}=\clm'$ holds when $\omega$ is the unique trace on 
$\clb$ as $v^*_k=S^*_k$ and $\clj \tilde{v}^*_k \clj = {1 \over d} S_k$ for all $1 \le k \le d$ where $P \neq Q$ and $\pi(\clo_d)' \supset \pi(\tilde{\clo}_d)''$. In the last proposition we have also proved if $[\tilde{\clm}\Omega]=P$ then 
$\clm'=\tilde{\clm}$. Thus a natural question that arises here: how the equality $P=Q$ is related to 
purity of $\omega$?  We are now in a position to state the main mathematical result of this section.

\vsp 
\begin{thm} Let $\omega$ be as in Theorem 3.5. Then the following holds:

\vsp 
\NI (a) $P$ is also the support projection of $\psi$ in $\pi(\tilde{\clo}_d))''_{|}\tilde{\clh}$ if and only if $\omega$ is pure. 

\vsp 
\NI (b) If $\omega$ is pure then the following holds:   

\NI (i) $\clm'= \tilde{\clm}$ where $\tilde{\clm} = \{ P\tilde{S}_iP: 1 \le i \le d
\}''$;

\NI (ii) $\pi(\clo_d)' = \pi(\tilde{\clo}_d)''$.

\NI (iii) $\pi_{\omega}(\clb_R)' = \pi_{\omega}(\clb_L)''$;

\end{thm}

\vsp
\NI {\bf PROOF: } First we will prove that $\omega$ is pure if $P$ is also the support projection of the state $\psi$ in $\pi(\tilde{\clo}_d)''\tilde{\clf}$, where $\tilde{\clf}=[\pi(\tilde{\clo}_d)''\Omega]$. 
The support projection of $\psi$ in $\pi(\tilde{\clo}_d)''\tilde{\clf}$ is $\tilde{\cle}\tilde{\clf}$ and thus we also 
have $P = \tilde{\cle}\tilde{\clf}$ by our hypothesis. Since $\Lambda^n(P)=\Lambda^n(\cle)\clf \uparrow \clf$ and now
$\Lambda^n(P)=\tilde{\cle}\Lambda^n(\tilde{\clf}) \uparrow \tilde{\cle}$ as $n \uparrow \infty$, we also have 
$\clf = \tilde{\cle}$. Similarly we also have for each $n$, $\cle\tilde{\Lambda}^n(\clf)=\tilde{\Lambda}^n(\tilde{\cle})\tilde{\clf}$ 
and thus taking limit we also get $\cle=\tilde{\clf}$.  

\vsp 
So we have $P=\cle \clf=\cle \tilde{\cle}=Q$. $\tilde{\clm}=P\pi(\tilde{\clo}_d)''P$ is cyclic in $\clk$ i.e. 
$[\tilde{\clm}\Omega]=[P\pi(\tilde{\clo}_d)''P\Omega]= P\tilde{\clf}=P\cle=P$ as $\tilde{\clf}=\cle$.  

\vsp 
However $\bigcap_{n \raro \infty} \tilde{\Lambda}^n(\pi(\tilde{\mbox{UHF}}_d)) = \pi(\tilde{\mbox{UHF}}_d)'' \bigcap \pi(\tilde{\mbox{UHF}}_d)'$ ( for a proof which is a simple application of Cuntz relation, we refer to section 5 of [Mo2]). 
Further $\psi$ being a factor state in $K_{\omega'}$, by Proposition 3.2 $\pi(\tilde{\mbox{UHF}}_d)''$ is a factor. 
In particular we have $\bigcap_{n \raro \infty} \tilde{\Lambda}^n(\pi(\tilde{\mbox{UHF}}_d))\tilde{\clf} = \!C \tilde{\clf}$. Thus by Proposition 1.1 in [Ar2] we conclude that $||\Psi \circ \tilde{\tau}_n -\phi_0|| \raro 0$ 
as $n \raro \infty$ for all normal state $\Psi$ on $\tilde{\clm}_0$ where $\tilde{\clm}_0=P \pi(\tilde{\mbox{UHF}}_d)''P$ as $\tilde{\clf}=\cle$ and support projection of $\psi$ in $\pi(\tilde{\mbox{UHF}}_d)''$ 
is $\tilde{\cle}$ and $P=\cle\tilde{\cle}\cle$. 

\vsp 
Note that $\tilde{\clm}_0 \subseteq \clm_0'$ where $\clm_0=P\pi(\mbox{UHF}_d)''P$. Further by Proposition 2.5 
$\clm_0 = \{x \in \clm:\beta_z(x)=x;z \in H \}\;\;\;$ and  $\tilde{\clm}_0=\{x \in \tilde{\clm}:
\beta_z(x)=x;z \in H \}$. Once we set $P_0=[\clm_0\Omega]$ then we also have $P_0=[\tilde{\clm}_0\Omega]$ 
as $[\tilde{\clm}\Omega]=P=[\clm\Omega]$ by expending $u_z=\sum_{k \in \hat{H}} z^kP_k$ where $z \raro u_z=
PU_zP$ is an unitary representation of group $H$.      

\vsp
For $x \in \tilde{\clm}, y \in \tilde{\clm}'$ we have  
$$\phi_0(\tilde{\tau}(x)y)=\sum_k <\tilde{v}^*_k\Omega,x\tilde{v}^*_ky\Omega> = \sum_k <v^*_k\Omega,x y \tilde{v}^*_k\Omega>$$  
(as $v_k^*\Omega=\tilde{v}^*_k\Omega$ )
$$=\sum_k <\Omega, x v_kyv_k^* \Omega>=\phi_0(x \tau(y))$$
The dual group of $(\tilde{\clm},\tilde{\tau},\phi_0)$ is given on the commutant by 
$(\tilde{\clm}',\tau,\phi_0)$ where $\tau(x)=\sum_k v_k x v_k^*$ for $x \in \tilde{\clm}'$. 
where commutant is taken in $\clb(\clk)$. 
Now moving to $\{\beta_z: z \in H \}$ invariant elements in the duality relation above, we verify that 
adjoint Markov map of $(\tilde{\clm}_0,\tilde{\tau},\phi_0)$ is given by 
$(\tilde{\clm}'_0,\tau,\phi_0)$ where $\tilde{\clm}'_0$, the commutant of $\tilde{\clm}_0$ 
is taken in $\clb(\clk_0)$ and $\clk_0$ is the Hilbert subspace $P_0$ with $\Omega$ as
cyclic and separating vector for $\tilde{\clm}_0$ in $\clk_0$. Thus by Theorem 2.4 in [Mo3], 
there exists a sequence of elements $y_n \in \tilde{\clm}'_0$ such that $y_n\tau_n(y) \raro 
\phi_0(y)1$ as $n \raro \infty$ for all $y \in \tilde{\clm}'_0 \subseteq \clb(\clk_0)$. 
Thus $\omega$ is pure by Proposition 3.5 (e) once we recall $\clm'=\tilde{\clm}$ as $P=Q$ 
( and so $\clm'_0 = \tilde{\clm}_0$ ) by Proposition 3.5 (d). 
 
\vsp  
Now we aim to prove $\tilde{\clf}=\cle$ and $\clf=\tilde{\cle}$ if $\omega$ is pure. We set unitary operator $V=\sum_k S_k\tilde{S}^*_k$. That $V$ is an unitaryoperator follows by Cuntz's relations and commuting property of $(S_i)$ and $(\tilde{S}_i)$. Further a simple computation shows that $V\pi(x)V^*=\pi(\theta(x))$ for all $x \in \clb = \clb_L \otimes \clb_R$, 
identified with $\tilde{\mbox{UHF}}_d \otimes \mbox{UHF}_d$ and $\theta$ is the right shift. We also have    
$$V \cle V^*= \sum_{k,k'}S_k\tilde{S}^*_k \cle S^*_{k'}\tilde{S}_{k'}$$
$$=\Lambda(\cle) \ge \cle$$ So $V(I-\cle)V^* \le I-\cle$, i.e. $(I-\cle)V^*\cle=0$.  
Also for any $X \in \pi(\clo_d)'$ we have $V^*\tilde{\clf}X\Omega = 
\tilde{\clf} \sum_k \tilde{S}_k X \tilde{S}_k^*\Omega$ as $S_k^*\Omega 
=\tilde{S}_k^*\Omega$. Thus $(I-\tilde{\clf})V^*\tilde{\clf}=0$ i.e. 
$V\tilde{\clf}V^*=\Lambda(\tilde{\clf}) \ge \tilde{\clf}$. Similarly we also 
have $V^*\tilde{\cle}V \ge \tilde{\cle}$ and $V^* \clf V \ge \clf$. We set two family 
of increasing projections for all natural numbers $n \in \!Z$ as follows 
$$\cle_n=V^n \cle (V^n)^*,\;\;\tilde{\clf}_n=V^n\tilde{\clf}(V^n)^*$$

\vsp 
Since $\beta_z(V)=V$ for all $z \in H$, $V \in \pi(\tilde{\mbox{UHF}}_d \otimes \mbox{UHF}_d)''$ by Proposition 3.4 as $\omega$ is pure. $\omega$ being also a factor state, we have $<f,V^ng> \raro <f,\Omega><\Omega,g>$ as $n \raro + \mbox{or} - \infty$ for any $f,g \in \pi(\clb_{loc}))\Omega$ by Power's criteria [Po1]. Since such vectors are dense in the Hilbert space topology and the family $\{V^n:n \ge 1\}$ is uniformly bounded, we get $V^n \raro |\Omega><\Omega|$ in weak operator topology as $n \raro + \mbox{or} - \infty$.   

\vsp 
For the time being we assume that $H$ is trivial. Otherwise the argument that follows here we can use for the representation $\pi_0$ of $\tilde{\mbox{UHF}}_d \otimes \mbox{UHF}_d$ i.e. $\pi$ restricted to $[\pi(\tilde{\mbox{UHF}}_d \otimes \mbox{UHF}_d)\Omega]$. 

\vsp 
We have following distinct cases: 

\vsp 
\NI {\bf Case 1.} $\cle \neq I( \tilde{\cle} \neq I)$. Let $\cle_n \raro \cle_{-\infty}$ as $n \raro -\infty$ and thus $V\cle_{-\infty}V=\cle_{-\infty}$. We claim that either $\cle_{-\infty} = |\Omega><\Omega|$ or $\cle_{-\infty}$ is a proper infinite dimensional projection i.e. if $\cle_{- \infty}$ is a finite projection then $\cle_{- \infty}=|\Omega><\Omega|$. Suppose not then the finite subspace is shift invariant. In particular there exists an unit  vector $f$ orthogonal to $\Omega$ such that $Vf=zf$ for some $z \in S^1$ and this contradicts weak mixing property i.e. $V^n \raro |\Omega><\Omega|$ in weak operator topology proved above as point spectrum of $V$ has only $1$ with spectral multiplicity $1$. 

\vsp 
If $\cle_{-\infty}$ is infinite dimensional we can get an unitaryoperator $U_0$ from 
$F_0=\tilde{\clh} \otimes_{\clk} \clh$ onto $\cle_{- \infty}$ and 
via the unitary map we can get a sequence of increasing projections $U_0\cle_nU^*_0$ in $\cle_{- \infty}$ and note that $U_0\cle_nU_0^*=V^nU_0 \cle U_0^*(V^n)^*$. Note that if $\cle_{- \infty}$ is infinite dimension the process will not stop in finite step. Thus we have $F_0 \ominus \Omega = \oplus_{1 \le k \le n_{\cle}} F(k)$ where the index set is either singleton or infinity and each $F(k)$ will give a system of imprimitivity with respect to $V$, where $F(1)
=F_0-\cle_{-\infty}$. Further $\tilde{\mbox{UHF}}_d$ being a simple $C^*$-algebra, each such imprimitivity sysyem 
is of Mackey index $\aleph_0$ [Mo3,section 4]: We fix a nonzero $f \in \cle-\theta^{-1}(\cle) \neq 0$ otherwise $\cle=I$ as $\theta^n(\cle) \uparrow I$ as $n \raro \infty$. $\pi_f:x \raro \theta^{-1}(x)f$ gives a representation of $\tilde{\mbox{UHF}}_d=\clb_L$ and we check that $[\pi_f(\theta^{-1}(\clb_L)''f] \le \cle-\theta^{-1}(\cle)$ 
as $f \perp [\theta^{-1}(\pi(\clb_R)')\Omega] \ge [\theta^{-1}(\pi(\clb_L))\Omega]$. Thus simplicity 
ensures that $\cle - \theta^{-1}(\cle)$ is a projection of dimension $\aleph_0$. Further $n_\cle$ is either $1$
or $\aleph_0$ since $F_0$ is separable.   

\vsp 
Since $\tilde{\clf}$ is also a proper projection, same argument is valid for $\tilde{\clf}$ with $\tilde{\clf}_{-\infty}=
\mbox{lim}_{n \raro -\infty} \theta^n(\tilde{\clf})$ i.e. we can write $F_0 \ominus \Omega = \oplus_{1 \le k \le n_{\tilde{\clf}}} \tilde{F}(k)$, where each $\tilde{H}(k)$ give rises to a system of imprimitivity with respect to $V$ where each system of imprimitivity is of Mackey index $\aleph_0$ where $\tilde{F}(1) = F_0-\tilde{\clf}_{-\infty}$ 
and $n_{\tilde{\clf}}$ is either $1$ or $\aleph_0$.        

\vsp 
In the following we use temporary notation $H$ for Hilbert subspace $F_0$. For a cardinal number $n$, we 
amplify a representation $\pi:\clb \raro \clb(H)$ of the $C^*$ algebra $\clb$ to $n$ fold direct sum 
$n\pi=\oplus_{1 \le k \le n } \pi_k $ acting on $n \clh=\oplus_{1 \le k \le n } H_k$  defining by 
$$n\pi(x) (\oplus \zeta_k) = \oplus (\pi(x)\zeta_k)$$ 
where $\pi_k = \pi$ is the representation of $\clb=\mbox{UHF}_d \otimes \tilde{\mbox{UHF}}_d$ 
on $H_k=H$ where $H =[\pi(\tilde{\mbox{UHF}}_d \otimes \mbox{UHF}_d)\Omega]$. We also 
extend $\bar{\tilde{F}} = \oplus \tilde{F}_{\alpha}$, $\bar{E} = \oplus E_{\alpha}$ and $\bar{V}=\oplus_{1 \le k \le n} V_k$ respectively. We also set notation $\Omega_k=\oplus_{1 \le k \le n}\delta^k_j\Omega$. 

\vsp 
Thus by Mackey's theorem, there exists a cardinal number $n \in \aleph_0$ and an unitary operator $U: nH \raro nH$ so that 
$\bar{V}=U\bar{V}U^*$ and $\bar{\cle}=U\bar{\tilde{\clf}}U^*$. We set a representation 
$\pi^U:\clb \raro \clb(nH)$ by $\pi^U(x)= Un\pi(x)U^*$
and rewrite the above identity as 
$$\oplus_{1 \le k \le n}[\pi_k(\mbox{UHF}_d)'\Omega_k] = \oplus_{1 \le k \le n}
[\pi^U_k(\tilde{\mbox{UHF}}_d)''\Omega_k]$$
where $\pi^U_k(x)=U\pi_k(x)U^*$. Note that by our construction we 
can ensure $U\Omega_k=\Omega_k$ for all $1 \le k \le n$ as the operator 
intertwining between two imprimitivity systems are acting on the orthogonal 
subspace of the projection generated by vectors $\{\Omega_k:1 \le k \le n\}$.

\vsp 
We claim $\cle=\tilde{\clf}$. Suppose not i.e. $\tilde{\clf} < \cle$. In such a case we 
have
$$\oplus_{1 \le k \le n}[\pi_k(\mbox{UHF}_d)'\Omega_k] < \oplus_{1 \le k \le n}
[\pi^u_k(\mbox{UHF}_d)'\Omega_k]$$
Alternatively 
$$\oplus_{1 \le k \le n}[\pi_k(\tilde{\mbox{UHF}}_d)''\Omega_k] < \oplus_{1 \le k \le n}
[\pi^u_k(\tilde{\mbox{UHF}}_d)''\Omega_k]$$

\vsp 
Thus in principle we can repeat our construction now with $\pi^U$ and so we 
get a strict partial ordered set of quasi-equivalent representation of $\clb$. 
In the following we now aim to employ formal set theory to bring a contradiction 
on our starting assumption that $\tilde{\clf} < \cle$. 

\vsp 
To that end we need to deal with more then one representation of $\clb$. For the rest of the proof we reset notation $\pi_0$ for $\pi$ used for the pure representation of 
$\clb$ in $H_0 =[\pi_0(\clb)\Omega_0]$ where $\Omega_0$ is the cyclic vector, the reset notation for $\Omega$. Let $\clp$ be the collection of representation $(\pi,H_{\pi},\Omega)$ 
quasi-equivalent to $\pi_0:\clb \raro \clb(H_0)$ with a shift invariant vector state $\omega(x)=<\Omega,\pi(x)\Omega>$ i.e. $\omega(\pi(\theta(x))=\omega(\pi(x))$. So there exists 
minimal cardinal numbers $n_{\pi},n_0(\pi)$ so that $n_{\pi}H_{\pi}$ is unitary equivalent to $n_0(\pi)\pi_0$. Thus for such an element $(\pi,H_{\pi},\Omega_{\pi})$ we can associate two cardinal numbers $n_{\pi}$ and $n_0(\pi)$ and without loss of generality we assume that $H_{\pi} \subseteq n_0(\pi)H_0$ and $n_{\pi}H_{\pi}=n_0(\pi)H_0.$  $\pi_0$ being a pure 
representation, any element $\pi \in \clp$ is a type-I factor representation of $\clb$. The interesting point here that $\dsp{\oplus_{\pi \in \clp }}\pi$ is also an element in $\clp$ with associated 
cardinal numbers $\sum_{\pi} n_{\pi}$ and $\sum_{\pi} n_0(\pi)$. We say $(\pi_1,H_1,\Omega^1) << (\pi_2,H_2,\Omega^2)$ if there exists an 
isometry $U:n_{\pi_1} H_1 \raro n_{\pi_2} H_2$ so that 

\vsp 
\NI (C1) For each $1 \le \alpha \le n_{\pi_1} $ we have 
$U \Omega^1_{\alpha} = \Omega^2_{\alpha'}$ for some 
$1 \le \alpha' \le n_{\pi_2}$;        

\NI (C2) $n_{\pi_2}\pi_2(x)E'_2 = U n_{\pi_1}\pi_1(x) U^* $ where $\cle'_2 \in n_{\pi_2}\pi_2(\clb)'$; 
 
\NI (C3) $\oplus_{1 \le \alpha \le n_{\pi_1} }[\pi^{\alpha}_1(\mbox{UHF}_d)'\Omega^1_{\alpha}] < \oplus_{1 \le 
\alpha \le n_{\pi_2}}[\pi^{\alpha}_2(\mbox{UHF}_d)'\Omega^2_{\alpha}]E'_2$.

\vsp 
In the inequality we explicitly used that both Hilbert spaces are subspaces of $nH_0$ for some possibly larger cardinal 
number $n$. That the partial order is non-reflexive follows as $(\pi,H,\Omega) << (\pi,H,\Omega)$ contradicts 
(C3) as $I=E_2'$. Partial order property follows easily. If $\pi_1 << \pi_2$ and $\pi_2 << \pi_3$ then $\pi_1 << \pi_3$. If $U_{12}$ and $U_{23}$ are isometric operators that satisfies (C1)-(C3) respectively, then $U_{13}=U_{23}U_{12}$ will do the job for $\pi_1$ and $\pi_3$. Thus $\pi^U \in \clp$ and by our starting assumption that $\tilde{\clf} \neq \cle$ we also 
check that $\pi_0 << \pi^U$. Thus going via the isomorphism we also check that for a given element $\pi \in \clp$ there exists an element $\pi' \in \clp$ so that $\pi << \pi'$. Thus $\clp_0$ is a non empty set and has at least one infinite chain containing $\pi_0$.
 
\vsp 
However by Hausdorff maximality theorem there exists a non-empty maximal totally ordered subset $\clp_0$ of $\clp$ 
containing $\pi_0$. We claim that $\pi_{max} = \oplus_{\pi \in \clp_0} \pi$ on $H_{\pi_{max}}=\oplus_{\pi \in \clp_0} H_{\pi} $ is an upper bound in $\clp_0$. That $\pi_{max} \in \clp$ is obvious. Further given an element $(H_1,\pi_1,\Omega_1) \in \clp_0$ there exists an element $(H_2,\pi_2,\Omega_2) \in \clp_0$ so that $\pi_1 << \pi_2$ by our starting remark as $\pi_0 << \pi^U$. By extending isometry $U_{12}$ to an isometry from $H_1 \raro n_{\pi_{max}} H_{\pi_{max}}$ trivially we get the required isometry that satisfies (C1),(C2) and (C3) where cardinal numbers $n_{\pi_{max}} = \sum_{\pi \in \clp_0} n_{\pi} \in \aleph_0$. Thus by maximal property of $\clp_0$ we have $\pi_{max} \in \clp_0$. This brings a contradiction as by our construction $(\pi_{max},H_{\pi_{max}},\Omega) << (\pi_{max},H_{\pi_{max}},\Omega)$ as $\pi_{max} \in \clp_0$ but partial order is strict. This contradicts our starting hypothesis that $\tilde{\clf} < \cle$. This completes the proof that $\tilde{\clf}=\cle$ when $\cle \neq I$. By symmetry of the argument we also get $\clf=\tilde{\cle}$ when $\tilde{\cle} < 1$.  

\vsp 
\NI {\bf Case 2:} $\cle=I ( \tilde{\cle}=I )$. We need to show $\tilde{\clf}=I (\clf=I)$ respectively. Suppose 
not and assume that both $\tilde{\clf}$ is a proper non-zero projection. 

\vsp 
We set projection $G$ on the closed linear span of elements in the subspaces $[\theta^{-n}(\clf)\pi(\tilde{\mbox{UHF}}_d)''\Omega]$ for all $n \ge 0$. We recall that $\theta(X)=VXV^*$ 
where $V=\sum_k S_k\tilde{S}_k^*$ and $\theta^{-1}(X)=\tilde{\Lambda}(X)$ 
for $X \in \pi(\tilde{\mbox{UHF}}_d)''$. Thus we have 
$$V^*\theta^{-n}(\clf)\pi(\tilde{\mbox{UHF}}_d)''\Omega$$
$$=\theta^{-n-1}(\clf)V^*\pi(\tilde{\mbox{UHF}}_d)''V\Omega$$
$$=\theta^{-n-1}(\clf) \tilde{\Lambda}(\pi(\tilde{\mbox{UHF}}_d)''\Omega.$$
Thus $(1-G)V^*G=0$ i.e. $\theta(G) \ge G$. It is also clear that $\tilde{\clf} \le G$ 
as the defining sequence of subspaces of $G$ goes to precisely $\tilde{\clf}$ as $n \raro \infty$ ( recall 
that $\theta^{-n}(\clf)=\tilde{\Lambda}_n(\clf) \uparrow I$ strongly as $n \uparrow \infty$ ). Once more we have 
$\theta^n(G) \ge \theta^n(\tilde{\clf}) = \Lambda^n(\tilde{\clf}) \uparrow I$ as $n \uparrow 
\infty$.  

\vsp 
If $G$ is a proper projection we can follow the steps as in the case 1 to find an unitary operator 
$U:nH \raro nH$ with $U\Omega_k=\Omega_k$ and $U\bar{V}U^*=\bar{V}$ so that $U\bar{G}U^*=\bar{\tilde{\clf}}$. 
We consider the subset $\clp_G$ of elements in $\clp$ for which $\cle_{\pi}=1$ and 
$\{\theta^{-n}(\clf_{\pi}):n \ge 0 \}$ commutes with $\tilde{\clf}_{\pi}$ 
and modify the strict partial ordering by modifying (C3) as  

\vsp 
\NI (C3') $\oplus_{1 \le \alpha \le n_{\pi_1} }[\pi^{\alpha}_1(\tilde{\mbox{UHF}}_d)''\Omega^1_{\alpha}] <  \oplus_{1 \le 
\alpha \le n_{\pi_2}}[\pi^{\alpha}_2(\tilde{\mbox{UHF}}_d)''\Omega^2_{\alpha}]E'_2$

\NI So we also get $\pi^U \in \clp_G$ and $\pi_0 << \pi^U$ and going along the same line we conclude that 
$G=\tilde{\clf}$. Thus we conclude that $G$ is either equal to $1$ or $G=\tilde{\clf}$. 

\vsp 
\NI {\bf Sub-case 1 of case 2:} If $G=I$ then $\clf G = \clf$ and so $[\clf\pi(\tilde{\mbox{UHF}}_d)''\Omega]=\clf$ as $\theta^{-n}(\clf) \ge \clf$. Thus $\tilde{\clf} \ge \clf$. So $\tilde{\clf} \ge \tilde{\Lambda}^n(\clf)$ for all $n \ge 1$ and taking limit we get $\tilde{\clf} \ge I$ i.e. $\tilde{\clf}=I$. This contradicts our starting assumption that $\tilde{\clf}$ is a proper projection.  

\vsp 
\NI {\bf Sub-case 2 of case 2:} Now we consider the case $G=\tilde{\clf} < I$. In such a case we have 
$(1- \tilde{\clf})\theta^{-n}(\clf)\tilde{\clf}=0$ and so $\theta^{-n}(\clf)$ commutes with $\tilde{\clf}$ for 
all $n \ge 0$. 

\vsp 
Now we set projection $\clf'$ defined by $\clf'=\clf-\clf\tilde{\clf}+|\Omega><\Omega|$ and check 
by commuting property of $\tilde{\clf}$ with $\clf$ that 
$$\clf'\theta^{-1}(\clf')\clf' = \clf(I-\tilde{\clf})\theta^{-1}(\clf)(I-\theta^{-1}(\tilde{\clf}))\clf(I-\tilde{\clf}) 
+ |\Omega><\Omega| $$
$$= \clf(I-\tilde{\clf})+|\Omega><\Omega|$$ 
as $\theta^{-1}(\clf) \ge \clf$ and $\theta^{-1}(\tilde{\clf}) \le \tilde{\clf}$. 
If $\clf'-|\Omega><\Omega|=\clf(I-\tilde{\clf}) \neq 0$ and so nor equal to $I$. 
If $\theta^{-n}(\clf') \uparrow I+ |\Omega><\Omega| - \mbox{lim}_{n \raro \infty} 
\theta^{-n}(\tilde{\clf}) \neq I$, we get the orthogonal 
projection i.e. $\mbox{lim}_{n \raro \infty} \theta^{-n}(\tilde{\clf})-|\Omega><\Omega|$ is
$\theta$ invariant and as in case-1 it would be an infinite dimensional separable 
Hilbert subspace. Thus we can follow the steps of case-1 with elements $\clf',\clf$ 
replacing the role of $\tilde{\clf},\cle$ to get an unitary operator $U:nH \raro nH$ 
so that $U\bar{V}=\bar{V}U$ and $U\bar{\clf}U^*=\bar{\clf'}$ for a cardinal number $n$.    

\vsp 
Now we consider a further subset $\clp_{G'}$ of $\clp_G$ consist of quasi-equivalent representations $\pi$ to $\pi_0$ 
of $\clb$ where $\pi$ admits the additional property: $\cle_{\pi}=I$ and $\{ \theta^{-n}(\tilde{\clf_{\pi}}): n \ge 0 \}$ 
commutes with $\clf_{\pi}$ satisfying $\clf_{\pi}\tilde{\clf}_{\pi} \neq \clf_{\pi}$ with the strict partial ordering 
$\pi_1 << \pi_2$ given by modifying condition (C3') as

\vsp 
\NI (C3'') $\oplus_{1 \le \alpha \le n_{\pi_1} }[\pi^{\alpha}_1(\tilde{\mbox{UHF}}_d)''\Omega^1_{\alpha}] >  \oplus_{1 \le 
\alpha \le n_{\pi_2}}[\pi^{\alpha}_2(\tilde{\mbox{UHF}}_d)''\Omega^2_{\alpha}]E'_2$

\vsp 
Since $\pi^U$ also satisfies the conditions that of $\pi_0 \in \clp_{G'}$ by covariance relation 
of $U$ with respect to shifts once more we get $\pi^U \in \clp_{G'}$ and $\pi_0 << \pi^U$. Thus we can 
repeat the process and so $\clp_{G'}$ has at least one infinite chain of totally ordered containing $\pi_0$.   

\vsp 
Once more by Hausdorff maximality principle we bring a contradiction to our starting assumption 
that $\clf' < \clf$. Thus we have shown either $\tilde{\clf}\clf=\clf$ or $\tilde{\clf}\clf=|\Omega><\Omega|$. 

\vsp 
However $\Lambda(\tilde{\clf}\clf)=\Lambda(\clf\tilde{\clf}\clf)= \clf\Lambda(\tilde{\clf})\clf \ge \clf\tilde{\clf}\clf=\tilde{\clf}\clf$ 
and similarly $\tilde{\Lambda}(\tilde{\clf}\clf) \ge \tilde{\clf}\clf$ 
i.e. subharmonic projections for both the endomorphisms. If $\clf\tilde{\clf}=|\Omega><\Omega|$ 
we get by sub-harmonic property that $\omega(s_Is^*_J) = l_Il_J^*$ where $(l_i)$ is a scalers i.e. $\omega$ 
is a pure product state on $\cla_R$ [BJP]. So $\cle\clf=|\Omega><\Omega|$ and 
$\clf=|\Omega><\Omega| \le \tilde{\clf}$ as $\cle=1$ and $\tilde{\clf}\Omega=\Omega$. 

On the other hand $\clf$ commutes with $\tilde{\clf}$ and $\clf\tilde{\clf}=\clf$ 
also says that $\clf \le \tilde{\clf}$. Thus in any case we have $\tilde{\clf} \ge \clf$ 
and so $\tilde{\Lambda}^n(\clf) \le \tilde{\clf}$ for all $n \ge 1$. Taking limit 
$n \raro \infty$, we get $I \le \tilde{\clf}$ i.e. $\tilde{\clf}=I$. This brings a 
contradiction to our starting hypothesis that $\tilde{\clf}$ is a proper projection 
i.e. $\tilde{\clf} < \cle=I$. This completes the proof of $\tilde{\clf}=\cle$ for 
the case when $H$ is trivial closed subgroup of $S^1$.   

\vsp 
Now we remove the assumption that $H$ is trivial. By our construction
$Q_0=[\pi_0(\mbox{UHF})'\Omega] [\pi_0(\tilde{\mbox{UHF}}_d)'\Omega]$ as 
$\cle=[\pi(\mbox{UHF})'\Omega]$, $\tilde{\cle} = [\pi(\tilde{\mbox{UHF}}_d)'\Omega]$ 
and $Q=\cle\tilde{\cle}$. As $\clf=[\pi(\clo_d)''\Omega]$ and $P=\cle\clf$ we get $P_0=[\pi_0(\mbox{UHF})'\Omega] [\pi_0(\mbox{UHF}_d)''\Omega]$.  

\vsp 
For the situation case-1 (i.e. $\cle \neq I$ ) one can simply replace $\pi$ by $\pi_0$ to get a proof from above that we have $P_0=Q_0$. For case-2 we need to modify the argument with $\clf$ and $\tilde{\clf}$ replaced by 
$\clf F_0$ and $\tilde{\clf} F_0$. Since $\clf$ and $\tilde{\clf}$ are $(\beta_z:z \in H)$ invariant, by purity we have 
$\clf$ and $\tilde{\clf}$ are elements in $\pi(\tilde{\mbox{UHF}}_d \otimes \mbox{UHF}_d)''$ and thus commutes with $F_0$. As $V$ commutes with $F_0$, the argument will go through for the situation to conclude that either $\tilde{\clf}F_0 \ge \clf F_0$ or $\theta^{-n}(\clf)F_0$ commutes with $\tilde{\clf}F_0$. 

\vsp 
In case $\tilde{\clf}F_0 \ge \clf F_0$ we get $\tilde{\clf} F_k \ge \tilde{\Lambda}^k(\clf)F_k \ge \clf F_k$ for all $k \in \hat{H}$ but $k \ge 0$. In case $\hat{H}=\IZ$, we also deduce $\tilde{S}^*_I\tilde{\clf}\tilde{S}_I\tilde{S}^*_IF_0\tilde{S}_I \ge \clf F_{-k}$ where $|I|=k$. However $\tilde{\clf} \ge \tilde{S}^*_I\tilde{\clf} \tilde{S}_I$ and thus we also have $\tilde{\clf} F_{-k} \ge \clf F_{-k}$ for all $k \ge 1$. Thus summing up we get $\tilde{\clf} \ge \clf$ and as before we 
conclude that $\tilde{\clf}=I$. This bring a contradiction. 

\vsp 
Thus we are left to consider the situation case -2 where $\cle=I$ and $\{\theta^{-n}(\clf)F_0: n \ge \}$ 
commutes with $\tilde{\clf}F_0$ ( $H \neq \{ 1 \}$ ). In such a case we will have as before 
$\clf\tilde{\clf}F_{0}=|\Omega><\Omega|$ and so $\omega$ is a pure product state and $P=\cle\clf=|\Omega><\Omega|$.
Thus $\clf \le \tilde{\clf}$ as $\cle=I$ and $\tilde{\clf}\Omega=\Omega$. So $\tilde{\Lambda}^n(\clf) \le \tilde{\clf}$ 
for all $n \ge 1$. Taking limit we conclude $\tilde{\clf}=I$ which contradicts our starting assumption 
that $\tilde{\clf} < \cle=I$.   

\vsp 
Now we write the equality $P_0=Q_0$ as $\cle \clf F_0=\cle\tilde{\cle}F_0$ and apply $\Lambda$ on both side to 
conclude that $\Lambda(\cle)\clf F_1=\Lambda(\cle)\tilde{\cle} F_1$ and multiplying by $\cle$ from left we get 
$\cle\clf F_1=\cle\tilde{\cle}F_1$ as $\Lambda(\cle)\cle=\cle$ and thus we get $PF_1=QF_1$. By repeated 
application of $\Lambda$, we get $PF_m=QF_m$. Thus we get $P=\sum_k PF_k=\sum_k QF_k=Q$.
This completes the proof for $P=Q$. Similarly $\tilde{\clf} = \sum_k \tilde{\clf}F_k = \sum_k \cle F_k = \cle$ and also $\clf=\tilde{\cle}$.    

\vsp
Now we are left to prove those three statements given in (b). $\omega$ being pure we have $P=Q$ and thus 
by Proposition 3.5 we have $\tilde{\clm}=\clm'$ and $\pi(\clo_d)'=\pi(\tilde{\clo}_d)''$. We are left to show
$\pi_{\omega}(\clb_R)'=\pi_{\omega}(\clb_L)''$. For that we recall $F_0$ and check few obvious relation 
$F_0\pi(\clo_d)'F_0=F_0\pi(\tilde{\clo}_d)''F_0$ and $\pi_{\omega}(\clb_R)' \subseteq F_0\pi(\clo_d)'F_0$.
Since $F_0\pi(\tilde{\mbox{UHF}}_d)''F_0$ is equal to $(\beta_z:z \in H)$ invariant elements in 
$F_0\pi(\tilde{\clo}_d)''F_0$ and elements in $\pi_{\omega}(\clb_R)'$ are $(\beta_z: z \in H)$ invariant 
we conclude that $\pi_{\omega}(\clb_R)' \subseteq \pi_{\omega}(\clb_L)''$. Inclusion in other direction 
is obvious and thus Haag duality property (iii) holds. \qed

\section{Symmetry of a translation invariant pure state on $\clb$ }

\vsp 
In this section we investigate $\omega$ on $\clb$ with some additional natural discrete symmetry. Let $\psi$ be a $\lambda$-invariant state on $\clo_d$ and $\tilde{\psi}$ be the state on $\clo_d$ defined by 
$$\tilde{\psi}(s_Is_J^*)=\psi(s_{\tilde{J}}s_{\tilde{I}}^*)$$ 
for all $|I|,|J| < \infty$ and $(\clh_{\tilde{\psi}},\pi_{\tilde{\psi}},\Omega_{\tilde{\psi}})$ be the GNS 
space associated with $(\clo_d,\tilde{\psi})$. That $\tilde{\psi}$ is well defined follows once we check by (14) that $$\psi(s_{\tilde{J}}s_{\tilde{I}}^*)=\phi_0(v_{\tilde{J}}v_{\tilde{I}}^*) = 
\phi_0(\tilde{v}_I\tilde{v}_J^*)$$ and appeal to Proposition 2.3 by observing that cyclicity condition i.e. 
the closed linear span $P_0$ of the set of vectors $\{\tilde{v}_I^*\Omega: |I| < \infty \}$ is $\clk$, can be 
ensured if not true already by taking a new set of Popescu elements $\{P_0\tilde{v}_kP_0:1 \le k \le d \}$. 
Otherwise one may also recall Proposition 2.3 that the map $s_Is_J^* \raro \tilde{v}_I\tilde{v}^*_{J}$ being unital and completely positive [Po] ( in particular positive ), $\tilde{\psi}$ is a well defined state on $\clo_d$.       

\vsp
Similarly for any translation invariant state $\omega$ on $\clb$ we set translation invariant state $\tilde{\omega}$ by reflecting 
around the point ${1 \over 2}$ on $\clb$ by $$\tilde{\omega}(Q_{-l}^{(-l)} \otimes Q_{-l+1}^{(-l+1)} \otimes ... \otimes Q_{-1}^{(-1)} \otimes Q_0^{(0)} \otimes Q_1^{(1)} ... \otimes Q_n^{(n)})$$
\be 
= \omega(Q_n^{(-n+1)}... \otimes Q_1^{(0)} \otimes Q_0^{(1)} \otimes Q_{-1}^{(2)} 
\otimes ... Q_{-l+1}^{(l)} \otimes Q_{-l}^{(l+1)})
\ee
for all $n,l \ge 1$ and $Q_{-l},..Q_{-1},Q_0,Q_1,..,Q_n \in 
M_n(\IC)$ where $Q^{(k)}$ is the matrix $Q$ at lattice point $k$. We define $\tilde{\omega}$ on $\clb$ by extending linearly to any $Q \in \clb_{loc}$.

\vsp 
Note first that the map $\psi \raro \tilde{\psi}$ is a one to one and onto affine map in the convex set of $\lambda$ invariant
state on $\clo_d$. In particular the map $\psi \raro \tilde{\psi}$ 
takes an element from $K_{\omega}$ to $K_{\tilde{\omega}}$
and the map is once more one to one and onto. Hence for any extremal point $\psi \in K_{\omega}$, $\tilde{\psi}$ is also an 
extremal point in $K_{\tilde{\omega}}$. Using Power's criteria we also verify here that $\omega$ is an extremal state if 
and only if $\tilde{\omega}$ is an extremal state. However such a conclusion for a pure state $\omega$ is not so obvious. We 
have the following useful proposition. 

\vsp
\begin{pro}

Let $\omega$ be an extremal translation invariant state on $\clb$ and $\psi \raro \tilde{\psi}$ 
be the map defined for $\lambda$ invariant states on $\clo_d$. Then the following holds:

\NI (a) $\psi \in K_{\omega}$ is a factor state if and only if $\tilde{\psi} \in K_{\tilde{\omega}}$ is a factor state. 

\NI (b) $\omega$ is pure if and only if $\tilde{\omega}$ is pure. 

\NI (c) A Popescu systems $(\clk,\clm,v_k,\Omega)$ of $\psi$ satisfies Proposition 2.4 with 
$(\pi_{\psi}(s_k),\;1 \le k \le d,\;P,\;\Omega)$ i.e. the projection $P$ on the subspace $\clk$ 
is the support projection of the state $\psi$ in $\pi(\clo_d)''$ and $v_i=P\pi_{\psi}(s_i)P$ 
for all $1 \le i \le d$, then the dual Popescu systems $(\clk,\clm',\tilde{v}_k,\Omega)$ satisfies 
Proposition 2.4 with $(\pi_{\tilde{\psi}}(s_k),\;1 \le k \le d,\;P,\;\Omega)$ i.e. the projection
$P$ on the subspace $\clk$ is the support projection of the state $\tilde{\psi}$ in 
$\pi_{\tilde{\psi}}(\clo_d)''$ and $\tilde{v}_i=P\pi_{\tilde{\psi}}(s_i)P$ for all 
$1 \le i \le d$, if and only if $\{x \in \clb(\clk): \sum_k \tilde{v}_kx \tilde{v}_k^* = x  \} = \clm$.     

\end{pro}

\vsp
\NI {\bf PROOF: } Since $\omega$ is an extremal translation invariant state, by Power's criteria $\tilde{\omega}$ is also an 
extremal state. As an extremal point of $K_{\omega}$ is map to an extremal point in $K_{\tilde{\omega}}$ by one to one property of the map $\psi \raro \tilde{\psi}$, we conclude by Proposition 2.6 that $\psi$ is a factor state if and only if $\tilde{\psi}$
is a factor state. For (b) note that $\tilde{xy}=\tilde{x} \tilde{y}$ and $\tilde{x^*}=\tilde{x}^*$ by our definition. Thus 
$\tilde{\omega}(x^*y)= \omega(\tilde{x^*y})= \omega((\tilde{x})^*\tilde{y})$. Thus one can easily construct an unitary operator between the two GNS spaces associated with $(\clb,\omega)$ and $(\clb,\tilde{\omega})$ intertwining two representation modulo a reflection i.e. $U \pi_{\omega}(x) U^* = \pi_{\tilde{\omega}}(\tilde{x})$ and $U\Omega_{\omega} =\Omega_{\tilde{\omega}}$. Thus (b) is now obvious. (c) follows by the converse part of the Proposition 2.4 once 
applied to the dual Popescu systems $(\clk,\clm',\tilde{v}_k,\Omega)$. \qed 

\vsp 
Thus the state $\tilde{\omega}$ is translation invariant, ergodic, factor state, pure if and only if $\omega$ is translation invariant, ergodic, factor state, pure respectively. We say $\omega$ is {\it lattice symmetric } if $\tilde{\omega}=\omega$.

\vsp
For a $\lambda$ invariant state $\psi$ on $\clo_d$ we define as before a $\lambda$ invariant state
$\tilde{\psi}$ by
\be
\tilde{\psi}(s_Is^*_J)=\psi(s_{\tilde{I}}s^*_{\tilde{J}})
\ee
for all $|I|,|J| < \infty.$ It is obvious that $\psi \in K_{\omega'}$ if and only if $\tilde{\psi}
\in K_{\tilde{\omega}'}$ and the map $\psi \raro \tilde{\psi}$ is an affine map. In particular an
extremal point in $K_{\omega'}$ is also mapped to an extremal point of $K_{\tilde{\omega}'}$. It is also
clear that $\tilde{\psi} \in K_{\omega'}$ if and only if $\omega$ is lattice symmetric. Hence a lattice symmetric state $\omega$ determines an affine map $\psi \raro \tilde{\psi}$ on the compact convex set $K_{\omega'}$. Furthermore, if $\omega$ is also extremal on $\clb$, then the affine map takes 
extremal elements to extremal elements of $K_{\omega'}$. The set of extremal elements in $K_{\omega'}$ 
can be identified with $S^1 / H \equiv S^1$ or $\{1 \}$ and the restriction of the affine map on the set of extremal element is continuous in weak topology ( by Proposition 2.6 the map $z \raro \psi \beta_z$ is one to one and onto the set of extremal elements of $K_{\omega'}$ for a fixed extremal element $\psi \in K_{\omega'}$ ). 

\vsp
Thus there exists $z_0 \in S^1$ so that $\tilde{\psi}=\psi \beta_{z_0}$ and as $\tilde{\beta \beta_z}=\tilde{\beta} \beta_z$ for all $z \in S^1$, we get the affine map taking $\psi \beta_z \raro \psi \beta_{z_0} \beta_z $ and thus determines a continuous one to one and onto map on $S^1 / H$ and as $\tilde{\tilde{\psi}}=\psi$ its 
inverse is itself. Thus either the affine map has a fixed point or $z_0^2=1$ i.e. it is a rotation map 
by an angle $2\pi$ ( Here we have identified $S^1 / H $ with $S^1$ in case $H \neq S^1$ ). Thus there exists 
an extremal element $\psi \in K_{\omega'}$ so that either $\tilde{\psi}=\psi \beta_{\zeta}$ where 
$\zeta $ is either $1$ or $-1$ where we recall that we have identified $S^1 / H =S^1$ when $H \neq S^1$. Note that if we wish to remove the identification, then for $H=\{ z : z^n =1 \}$ for some $n \ge 1$, $\zeta$ is either $1$ or $exp^{\pi i \over n }$. Note that in case $H=S_1$ then $\tilde{\psi}=\psi$ for $\psi \in K_{\omega'}$ as $K_{\omega}$ is a singleton set by Proposition 2.6.

\vsp
\begin{pro} Let $\omega$ be a translation invariant lattice symmetric state on $\clb$. Then the following holds:

\NI (a) If $\omega$ is also an extremal translation invariant state on $\clb$ then $H=\{z \in S^1: \psi \beta_z =\psi \}$ 
is independent of $\psi \in K_{\omega'}$.  

\NI (b) If $H=\{z: z^n =1 \}$ for some $n \ge 0$ then $\tilde{\psi}=\psi \beta_{\zeta}$ for all $\psi \in K_{\omega'}$ where $\zeta$ is fixed either $1$ or $exp^{\pi i \over n }$. Let $(\clh,S_k,\;1 \le k \le d,\Omega)$ be the GNS space associated 
with $(\clo_d,\psi)$, $P$ be the support projection of the state $\psi$ in $\pi(\clo_d)''$ and $\clk= P \clh$ with Popescu 
systems $(\clk,\clm,v_k, \;1 \le k \le d,\;\Omega)$ as in Proposition 2.4 where $v_k=PS_kP$ and associated normal state $\phi_0$ 
on $\clm=\{ v_k,v_k^*: 1 \le k \le d \}''$ is invariant for $\tau(x)=\sum_k v_k x v_k^*$. Let $(\tilde{\clh},\tilde{S}_k \;1 \le k \le d,\Omega)$ be the Popescu minimal dilation in Theorem 2.1 of the dual Popescu systems 
$(\clk,\tilde{\clm},\tilde{v}_k,\;1 \le k \le d,\;\Omega)$ defined in Proposition 3.2. 
Then there exists an unitaryoperator $U_{\zeta}: \clh \otimes_{\clk} \tilde{\clh}$ so that
\be
U_{\zeta}^*=U_{\bar{\zeta}}, U_{\zeta}\Omega=\Omega,\;\;\;U_{\zeta} S_k U_{\zeta}^*= \beta_{\bar{\zeta}}(\tilde{S}_k)
\ee
for all $1 \le k \le d$. 

\vsp 
Furthermore if $\omega$ is also pure then there exists an unitaryoperator $u_{\zeta}: \clk \raro \clk$ 
so that
\be
u_{\zeta}\Omega=\Omega, \;\;\; u_{\zeta}v_ku_{\zeta}^*= \beta_{\bar{\zeta}}(\tilde{v}_k)
\ee
for all $1 \le k \le d$ and $u_{\zeta} \clj u_{\zeta}^*=\clj,\;\;u_{\zeta} \Delta^{1 \over 2} u_{\zeta}^*=\Delta^{-{1 \over 2}}$, $u_{\zeta}^*=u_{\bar{\zeta}}$ and $u_{\zeta} \clm u_{\zeta}^* = \clm',\;u^*_{\zeta}\clm u_{\zeta}=\tilde{\clm}$. Moreover $\clm' = \tilde{\clm}$. Further if $\zeta=1$ then $u_{\zeta}$ is self-adjoint and otherwise if $\zeta \ne 1$ then $u_{\zeta}^{2n}$ is self adjoint.

\NI (c) If $H=S^1$ then $K_{\omega'}$ is having only one element $\psi$, so $\psi=\tilde{\psi}$ and
(19) and (20) are valid with $\zeta=1$. 

\end{pro} 

\vsp
\NI {\bf PROOF: } (a) follows by Proposition 2.6. Now we aim to prove (b). For existence of an extremal state $\psi \in K_{\omega'}$ 
so that $\tilde{\psi}=\psi \beta_{\zeta}$ we refer to the paragraph preceding the statement of this proposition. As $\tilde{(\psi \beta_z)}
=\tilde{\psi} \beta_z$ for all $z \in S^1$, a simple application of Proposition 2.6 says that $\tilde{\psi} = \psi \beta_{\zeta} $ for 
all extremal points in $K_{\omega'}$ if it holds for one extremal element. Hence existence part in 
(b) is true by Krein-Millmann theorem.

\vsp
$\Omega$ is a cyclic vector for $\pi(\clo_d \otimes \tilde{\clo}_d)$ and thus we define
$U_{\zeta}:\clh \otimes_{\clk} \tilde{\clh} \raro \clh \otimes_{\clk} \tilde{\clh}$ by
$$U_{\zeta}: S_IS^*_J\tilde{S}_{I'}\tilde{S}_{J'}^*\Omega \raro 
\beta_{\bar{\zeta}}(S_{I'}S_{J'}^*\tilde{S}_I\tilde{S}^*_J) \tilde{\Omega}$$
That $U_{\zeta}$ is an unitary operator follows from (14) and the dual relation (18) along with our condition that $\tilde{\psi}=\psi \beta_{\zeta}$. By our construction we also have 
$U_{\zeta}S_k = \beta_{\bar{\zeta}}(\tilde{S}_k)U_{\zeta} $ for all $1 \le k \le d$. 
In particular $U_{\zeta}\pi(\clo_d)''U_{\zeta}^*=\pi(\tilde{\clo}_d)''$. 

\vsp 
$\omega$ being pure we have $P=Q$ by Theorem 3.6 as $\clf=\tilde{\cle}$ and so 
$UPU^*=UQU^*=U\cle \tilde{\cle}U^*=\tilde{\cle}\cle=Q=P$ which ensures an unitary operator $u_{\zeta}=PU_{\zeta}P$ on $\clk$ and a routine calculation shows that
\be
u_{\zeta}v^*_k u_{\zeta}^*= \beta_{\bar{\zeta}}(\tilde{v}^*_k)
\ee
for all $1 \le k \le d$. As $U_{\zeta}^*=U_{\bar{\zeta}}$ we have $u_{\zeta}^*=u_{\bar{\zeta}}$. If $\zeta \ne 1$, then
$\zeta^{2n}=1$ and thus $U^{2n}_{\zeta}$ is inverse of its own. Thus $u^{2n}_{\zeta}$ is self-adjoint. That
$\clm'=\tilde{\clm}$ by Theorem 3.6 (b). 

\vsp
In the following we consider the case $\zeta=1$ for simplicity of notation and otherwise for the case $\zeta \ne 1$ very little modification is needed in the symbols or simply reset temporary notation $\tilde{v}_k$ for 
$\bar{\zeta}\tilde{v}_k$ i.e. include the phase factor. 

\vsp
We denote $u_1=u$ in the following for simplicity. It is simple to verify now the following steps 
$uSv_Iv_J^*\Omega=uv_Jv_I^*\Omega = \tilde{v}_J\tilde{v}^*_I\Omega=F\tilde{v}_I\tilde{v}^*_J\Omega$ 
where $Sx\Omega=x^*\Omega,\;x \in \clm$ and $Fx'\Omega=x'^*\Omega,\;x' \in \clm'$ are the Tomita's 
conjugate operator. Hence $u \clj \Delta^{1 \over 2}=\clj \Delta^{-{1 \over 2}}u$, i.e $u \clj u^* 
u\Delta^{1 \over 2}u^* = \clj \Delta^{-{1 \over 2}}$ and by uniqueness of polar decomposition we 
conclude that $u \clj u^* = \clj$ and $u\Delta^{1 \over 2}u^*= \Delta^{-{1 \over 2}}$. That $u\clm u^*=
\tilde{\clm}$ is obvious. For $u\clm u^*=\tilde{\clm}$ we note that by our construction 
$U \tilde{S}_k U^* = S_k$ and so $U\pi(\tilde{\clo}_d)U^* = \pi(\clo_d)$ 
and hence projecting to its support projection we get the required relation. \qed 

\vsp
Now we introduce another useful symmetry on $\omega$. If $Q= Q^{(l)}_0 \otimes Q^{(l+1)}_1 \otimes ....\otimes Q^{(l+m)}_m$ we set $Q^t={Q^t}^{(l)}_0 \otimes {Q^t}^{(l+1)}_1 \otimes ..\otimes {Q^t}^{(l+m)}_m$
where $Q_0,Q_1,...,Q_m$ are arbitrary elements in $M_d$ and $Q_0^t,Q^t_1,..$ stands for transpose
with respect to an orthonormal basis $(e_i)$ for $\IC^d$ (not complex conjugate) of $Q_0,Q_1,..$ 
respectively. We define $Q^t$ by extending linearly for any
$Q \in \clb_{loc}$. For a state $\omega$ on UHF$_d$ $C^*$ algebra $\otimes_{\IZ} M_d$ we define
a state $\bar{\omega}$ on $\otimes_{\IZ} M_d$ by the following prescription
\be
\bar{\omega}(Q) =
\omega(Q^t)
\ee
Thus the state $\bar{\omega}$ is translation invariant, ergodic, factor state if and only if $\omega$ is
translation invariant, ergodic, factor state respectively. We say $\omega$ is {\it real } if
$\bar{\omega}=\omega$. In this section we study a translation invariant real state.

\vsp
For a $\lambda$ invariant state $\psi$ on $\clo_d$ we define a $\lambda$ invariant state $\bar{\psi}$ by
\be
\bar{\psi}(s_Is^*_J)=\psi(s_Js^*_I)
\ee
for all $|I|,|J| < \infty$ and extend linearly. That it defines a state follows as 
for an element $x= \sum c(I,J)s_Is^*_J$ we have $\bar{\psi}(x^*x)=\psi(y^*y) \ge 0$ 
where $y= \sum \overline{c(I,J)}s_Js_I^*$. It is obvious that $\psi \in K_{\omega'}$ if and only if $\bar{\psi}
\in K_{\bar{\omega}'}$ and the map $\psi \raro \bar{\psi}$ is an affine map. In particular an
extremal point in $K_{\omega'}$ is also mapped to an extremal point in $K_{\bar{\omega}'}$. It is also
clear that $\bar{\psi} \in K_{\omega'}$ if and only if $\omega$ is real. Hence a real
state $\omega$ determines an affine map $\psi \raro \bar{\psi}$ on the compact convex set $K_{\omega'}$. 
Furthermore, if $\omega$ is also extremal on $\clb$, then the affine map, being continuous on the set of extremal 
elements in $K_{\omega'}$, which can be identified with $S^1 / H \equiv S^1$ or $\{ 1 \}$ 
( by Proposition 2.6 ) by fixing an extremal element $\psi_0 \in K_{\omega'}$. In such a case there exists a unique   $z_0 \in S^1$ so that 
$\bar{\psi}_0 = \psi_0 \beta_{z_0}$ 

\vsp 
Now $\bar{\psi_0 \beta_z}= \bar{\psi}_0 \beta_{\bar{z}}$ for all $z \in S^1$, the affine map takes 
$\psi_0 \beta_z \raro \psi_0 \beta_{z_0 \bar{z}}$. If $z_0=1$ we get that the map fixes two point namely 
$\psi_0$ and $\psi_0 \beta_{-1}$. 

\vsp 
Even otherwise we can choose $z \in S^1$ so that $z^2=z_0$ and for such a choice we get an extremal element namely $\psi_0 \beta_z$ gets 
fixed by the map. What is also crucial here that we can as well choose $z \in S^1$ so that $z^2 =-z_0$, if so then 
$\psi_0 \beta_z $ gets mapped into $\psi_0 \beta_{z_0} \beta_{\bar{z}} = \psi_0 \beta_{-z}=\psi_0 \beta_z \beta_{-1}$. 
Thus in any case we also have an extremal element $\psi \in K_{\omega'}$ so that 
$\bar{\psi} = \psi \beta_{\zeta}$ where $\zeta \in \{1,-1\}$.   

\vsp
Thus going back to the original set up, we sum up the above by saying that if  $H=\{z: z^n= 1 \} \subseteq S^1$ and 
$\zeta \in \{1,exp^{i\pi \over n } \}$ then there exists an extremal element $\psi \in K_{\omega'}$ so that $\bar{\psi}
=\psi \beta_{\zeta}$.   

\vsp
\begin{pro} 
Let $\omega$ be a translation invariant real factor state on $\otimes_{\IZ}M_d$. Then the 
following holds: 

\NI (a) if $H=\{z: z^n= 1 \} \subseteq S^1$ and $\zeta \in \{1,exp^{i\pi \over n } \}$ then there exists an extremal 
element $\psi \in K_{\omega'}$ so that $\bar{\psi}=\psi \beta_{\zeta}$. Let $(\clh,\pi_{\psi}(s_k)=S_k,\;1 \le k \le d, 
\Omega)$ be the GNS representation of $(\clo_d,\psi)$, $P$ be the support projection of the state $\psi$ in $\pi(\clo_d)''$ 
and $(\clk,\clm,v_k,\;1 \le k \le d, \Omega)$ be the associated Popescu systems as in Proposition 2.4. Let $\bar{v}_k= \clj v_k \clj$ for all $1 \le k \le d$ and $(\bar{\clh},\bar{S}_k,P,\Omega)$ be the Popescu minimal dilation as described in Theorem 2.1 associated with 
the systems $(\clk,\clm',\bar{v}_k,\;1 \le k \le d,\Omega)$. Then there exists an unitary
operator
$W_{\zeta}: \clh \raro \bar{\clh}$ so that
\be
W_{\zeta}\Omega=\Omega,\;\;\;W_{\zeta}S_kW_{\zeta}^* = \beta_{\bar{\zeta}}(\bar{S}_k) 
\ee
for all $1 \le k \le d$. Furthermore $P$ is the support projection of the state $\bar{\psi}$ in $\bar{\pi}(\clo_d)''$ 
and there exists an unitaryoperator $w_{\zeta}$ on $\clk$ so that
\be
w_{\zeta}\Omega=\Omega, \;\;\; w_{\zeta}v_kw_{\zeta}^* = \beta_{\bar{\zeta}}(\bar{v}_k)= \clj \beta_{\zeta}(v_k) \clj
\ee
for all $1 \le k \le d$ and  $w_{\zeta} \clj w_{\zeta}^*=\clj$ and $w_{\zeta} \Delta^{1 \over 2} w_{\zeta}^* = 
\Delta^{-{1 \over 2}}$. $w_{\zeta}$ is self adjoint if and only if $\zeta=1$; 

\NI (b) If $H=S^1$, $K_{\omega'}$ is a set with unique element $\psi$ so that $\bar{\psi}=\psi$ and relations 
(24) and (25) are valid with $\zeta=1$.   
\end{pro}

\vsp
\NI {\bf PROOF: } For existence part in (a) we refer the paragraph above preceded the statement of the proposition. 
We fix a state $\psi \in K_{\omega'}$ so that $\bar{\psi}=\psi  \beta_{\zeta}$ and define $W:\clh \raro \bar{\clh}$ by
$$W_{\zeta}: S_IS^*_J\Omega \raro \beta_{\bar{\zeta}}(\bar{S}^*_{I}\bar{S}^*_{J})\Omega$$
That $W_{\zeta}$ is an unitaryoperator follows from (3.10) and thus
$W_{\zeta}S_k= \beta_{\bar{\zeta}}(\bar{S}_k)W_{\zeta}$ for all $1 \le k \le d$. For simplicity of notation we take the case $\zeta=1$
as very little modification is needed to include the case when $\zeta \ne 1$ or reset Cuntz elements by 
absorbing the phase factor in the following computation and use notation $W$ for $W_{\zeta}$.

\vsp
$P$ being the support projection we have by Proposition 2.4 that $\clm'= \{ x \in \clb(\clh): \sum_k v_k x v_k^*= x \}$
and thus $\clm = \{ x \in \clb(\clk): \sum_k \clj v_k \clj x \clj v_k^* \clj = x \}$. Hence by the converse part of 
Proposition 2.4 we conclude that $P$ is also the support projection of the state $\bar{\psi}$ in $\bar{\pi}(\clo_d)''$. 
Hence $W_{\zeta}PW_{\zeta}^*=P$. Thus we define an unitary operator $w_{\zeta}: \clk \raro \clk$ by
$w_{\zeta}= PW_{\zeta}P $ and verify that
$$\bar{v}^*_k = P\bar{S}^*_kP$$
$$ = PW_{\zeta} \beta_{\zeta}(S^*_k)W_{\zeta}^*P = PW_{\zeta}P \beta_{\zeta}(S^*_k)PW_{\zeta}^*P$$
$$= PW_{\zeta}P\beta_{\zeta}(v^*_k)PW_{\zeta}^*P=w_{\zeta}\beta_{\zeta}(v^*_k)w_{\zeta}^*.$$

\vsp
We recall that Tomita's conjugate linear operators $S,F$ [BR] are
the closure of the linear operators defined by $S:x\Omega \raro x^*\Omega$ for $x \in \clm$ and $F:y\Omega \raro y^*\Omega$
for $y \in \clm'$. We check the following relations for $\zeta=1$ with simplified notation $w_1=w$,
$$wSv_Iv^*_J\Omega=wv_Jv^*_I\Omega= \bar{v}_J \bar{v}^*_I\Omega$$
$$=F \bar{v}_I\bar{v}^*_J\Omega=Fw v_Iv^*_J\Omega$$ 
for $|I|,|J| < \infty$. Since such vectors are total, we have $wS=Fw$ on the domain of $S$. Thus
$wSw^*=F$ on the domain of $F$. We write $S=\clj \Delta^{1 \over 2}$ as the unique polar decomposition. Then $F=S^* =
\Delta^{1 \over 2} \clj= \clj \Delta^{-{1 \over 2}}$. Hence $w \clj w^*w\Delta^{1 \over 2}w^*=\clj \Delta^{-{1 \over 2}}$.
By the uniqueness of polar decomposition we get $w \clj w^*=\clj$ and $w \Delta^{1 \over 2} w^*= \Delta^{-{1 \over 2}}$.
Same algebra is valid in case $\zeta \ne 1$ if we reset the notations $\tilde{v}_k$ on the right hand side absorbing 
the phase factor. In the following we repeat it for completeness of the proof as the phase factor could be delicate.

$$w_{\zeta}Sv_Iv^*_J\Omega= w_{\zeta} v_Jv^*_I\Omega= w_{\zeta}v_Jv^*_Iw^*_{\zeta}\Omega$$
$$=\zeta^{|I|-|J|} \bar{v}_J \bar{v}^*_I\Omega= \zeta^{|I|-|J|} F \bar{v}_I\bar{v}^*_J\Omega$$
$$=F\zeta^{-|I|+|J|} \bar{v}_I\bar{v}^*_J\Omega\;\;\mbox{for}\;\; |I|,|J| < \infty$$
$$=F w_{\zeta} v_Iv^*_J\Omega$$
for all $|I|,|J| < \infty$. 

\vsp
Now we are going to show that $w_{\zeta}$ is self-adjoint if and only if $\zeta=1$. Note that $\beta_z(w_{\zeta}xw_{\zeta}^*)=w_{\zeta} \beta_{\bar{z}}(x)w_{\zeta}^*$ and thus 
applying $\beta_{\bar{\zeta}}$ on both side of the following identity
\be
w_{\zeta}v_kw^*_{\zeta}= \clj \beta_{\zeta}(v_k) \clj
\ee
for all $1 \le k \le d$, we also get 
$w_{\zeta} \beta_{\zeta}(v_k) w_{\zeta}^* = \clj \beta_{\bar{\zeta}}^2(v_k)\clj$ and thus
$w_{\zeta}^2 v_k (w_{\zeta}^*)^2 = w_{\zeta} \clj \beta_{\zeta}(v_k) \clj w_{\zeta}^*= \beta_{\bar{\zeta}^2}(v_k)$ 
as $\clj$ commutes with $w_{\zeta}$. 

\vsp 
$\zeta^2=1$ if and only if $\zeta=1$ ( as $\zeta=1$ or $exp^{i \pi \over n}$ where $n \ge 2$ ). In such a case we get $w_{\zeta}^2 \in \clm'$ and further 
as $w_{\zeta}$ commutes with $\clj$, $w_{\zeta}^2 \in \clm$. $\omega$ being an extremal element in $K_{\omega'}$ we have $\clm \vee \tilde{\clm}=\clb(\clk)$ 
by Proposition 3.5 and as $\tilde{\clm} \subseteq \clm'$, we get that $\clm$ is a factor. Thus for a factor $\clm$, $w_{\zeta}^2$ is a scaler. Since 
$w_{\zeta}\Omega=\Omega$ we get $w_{\zeta}^2=1$ i.e. $w^*_{\zeta}=w_{\zeta}$. This completes the proof. \qed 

\vsp
A state $\omega$ on $\otimes_{\IZ}M_d$ is said be in detailed balance if $\omega$ is both lattice symmetric and real.
In the following proposition as before we identified once more $S^1 / H \equiv S^1$ in case $H \neq S^1$ and set $\zeta$
be the least value in $S^1 \ H \equiv S^1 $ so that $\zeta^2 \in H$.   

\vsp
\begin{thm} 
Let $\omega$ be a translation invariant factor state on $\clb = \otimes_{\IZ}M_d$. Then the following 
are equivalent:

\NI (a) $\omega$ is real and lattice symmetric;

\NI (b) There exists an extremal element $\psi \in K_{\omega'}$ so that $\tilde{\psi}=\psi \beta_{\zeta}$ and 
$\bar{\psi}=\psi \beta_{\zeta}$, where $\zeta$ is either $1$ or $exp^{i\pi \over n }$.    

\vsp
Furthermore if $\omega$ is a pure state then the following holds:

\NI (c) There exists a Popescu elements $(\clk,v_k,\;1 \le k \le d,\Omega)$ for $\omega$ with relation  
$v_k = \clj_v \tilde{v}_k \clj_v $ for all $1 \le k \le d$, where $\clj_v=v\clj$ and $v$ is a self-adjoint unitary 
operator on $\clk$ commuting with modular operators $\Delta^{1 \over 2}$ and conjugate operator $\clj$ associated 
with cyclic and separating vector $\Omega$ for $\clm$. Further $\beta_z(v)=v$ for all $z \in H$ and  
$H \subseteq \{1, -1 \}$;

\NI (d) The map $\clj_v: \clh \otimes_{\clk} \tilde{\clh} \raro \clh \otimes_{\clk} \tilde{\clh}$ defined by 
$$\pi(s_Is^*_J\tilde{s}_{I'}\tilde{s}^*_{J'})\Omega 
\raro \pi(s_{I'}s^*_{J'}\tilde{s}_{I}\tilde{s}^*_{J})\Omega,$$
$\;|I|,|J|,|I'|,|J'| < \infty $ extends the map $\clj_v: \clk  \raro  \clk$ to an anti-unitary map so that $\clj_v 
\pi(s_i) \clj_v = \bar{\pi}(\tilde{s}_i)$ for all $ 1 \le i \le d$ where
$\bar{\pi}$ is the conjugate linear extension of $\pi$ from the generating 
set $(\tilde{s}_i)$, i.e. $\bar{\pi}(\tilde{s}_I\tilde{s}^*_J) = \pi(\tilde{s}_I\tilde{s}^*_J)$ for $|I|,|J| < \infty$ 
and then extend it anti-linearly for its linear combinations. 
\end{thm} 

\vsp
\NI {\bf PROOF: } Since $\omega$ is lattice symmetric, by Proposition 4.2 $\tilde{\psi}= \psi \beta_{\zeta}$ 
for all $\psi \in K_{\omega'}$ where $\zeta$ is fixed number either $1$ or $exp^{i \pi \over n}$ for some $n \ge 1$. 
Now we use real property of $\omega$ and choose by Proposition 4.3 an extremal element $\psi \in K_{\omega'}$ so 
that $\bar{\psi}=\psi \beta_{\zeta}$. This proves that (a) implies (b). That (b) implies (a) is obvious.

\vsp
Now we aim to prove the last statements which is the main point of the proposition. For simplicity of notation we consider 
the case $\zeta=1$ and leave it to reader to check that a little modification needed to include the case $\zeta \neq 1$ and 
all the algebra stays valid if $\tilde{v}_k$ is replaced by $\beta_{\zeta}(\tilde{v}_k)$. We consider the Popescu system $(\clk,\clm,v_k,
\;1 \le k \le d,\Omega)$ as in Proposition 2.4 associated with $\psi$. Thus by Proposition 4.3 and Proposition 4.4 there exists unitary operators $u_{\zeta},w_{\zeta}$ on $\clk$ so that
$$u_{\zeta}v_ku_{\zeta}^*=\beta_{\bar{\zeta}}(\tilde{v}_k)$$
$$w_{\zeta}v_kw_{\zeta}^* = \beta_{\bar{\zeta}}(\bar{v}_k) = \clj \beta_{\zeta}(v_k) \clj $$
where $\;\;u_{\zeta} \clj u_{\zeta}^*=\clj, \;\; 
\;w_{\zeta} \clj w_{\zeta}^*=\clj$ and
$u_{\zeta} \Delta^{1 \over 2}u_{\zeta}^*=w_{\zeta} \Delta^{1 \over 2}w^*_{\zeta}=\Delta^{-{1 \over 2}}$.
Thus
\be
u_{\zeta} w_{\zeta} v_k w_{\zeta}^* u_{\zeta}^* = u_{\zeta} \clj \beta_{\zeta}(v_k) \clj u_{\zeta}^* 
= \clj \beta_{\zeta}(u_{\zeta}v_ku_{\zeta}^*) \clj = \clj \beta_{\zeta}(\beta_{\bar{\zeta}}(\tilde{v}_k)\clj=
\clj \tilde{v}_k \clj
\ee

We also compute that
\be
w_{\zeta}u_{\zeta}v_ku_{\zeta}^*w_{\zeta}^* = w_{\zeta} \beta_{\bar{\zeta}}(\tilde{v}_k) w_{\zeta}^* = 
\clj \beta_{\zeta} \beta_{\bar{\zeta}}(\tilde{v}_k) \clj = \clj \tilde{v}_k \clj 
\ee

\vsp
By Proposition 3.2, for a factor state $\omega$ we also have $\clm \vee \tilde{\clm} = \clb(\clk)$. As $\tilde{\clm} \subseteq \clm'$, in particular we note that $\clm$ is a factor. So $u_{\zeta}^*w_{\zeta}^* u_{\zeta}w_{\zeta} \in \clm'$ commuting also with $\clj$ and thus a scaler as $\clm$ is a factor. As $u_{\zeta}\Omega=w_{\zeta}\Omega=\Omega$, we conclude that $u_{\zeta}$ commutes with $w_{\zeta}$.  

\vsp
Now we set $v_{\zeta}=u_{\zeta}w_{\zeta}$ which is an unitary operator commuting with both $\clj$ and
$\Delta^{1 \over 2}$. That $v_{\zeta}$ commuting with $\Delta^{1 \over 2}$ follows as $u_{\zeta}w_{\zeta}\Delta^{1 \over 2}=
u_{\zeta}\Delta^{- {1 \over 2}}w_{\zeta}= \Delta^{1 \over 2} u_{\zeta}w_{\zeta}$. 

\vsp
Next claim that we make now that $v_{\zeta}$ is a self-adjoint element. To that end note that the relations
(28) says that $v_{\zeta}\clm v_{\zeta}^* \subseteq \clm$ and so 
$v_{\zeta}\clm'v_{\zeta}^* \subseteq \clm'$. We check the following identity:
$v_{\zeta}\tilde{v}^*_kv_{\zeta}^*\Omega
=v_{\zeta}v_k^*\Omega=v_{\zeta}v_k^*v_{\zeta}^*\Omega
=\clj \tilde{v}_k^* \clj \Omega
=\clj v_k \clj \Omega$
and thus separating property we deduce that 
$$v_{\zeta}\tilde{v}_k^*v_{\zeta}^*=\clj v_k \clj$$
for all $1 \le k \le d$. So we conclude that $v^2_{\zeta} \in \clm'$ and as $v_{\zeta}$ commutes with $\clj$, $v_{\zeta}^2$ is an element in the 
centre of $\clm$. The centre of $\clm$ being trivial as $\omega$ is a factor state ( here we have more namely pure )
and $v_{\zeta}\Omega=\Omega$, we conclude that $v^2_{\zeta}$ is the unit operator. Hence $v_{\zeta}$ is a self-adjoint element. 

\vsp 
For simplicity of notation we set $v$ for $v_{\zeta}$. $\beta_z(v)=v$ for all $z \in H$ is equivalent to the property that 
$v$ keeps the subspaces $PF_k$ invariant. Since $vv_Iv_J^*v^*=\clj \tilde{v}\tilde{v}^*_J\clj$ and $v$ is self-adjoint 
and $v\Omega=\Omega$ we get $vv_Iv_J^*\Omega=\clj \tilde{v}\tilde{v}^*_J\Omega$. $\beta_z$ being an automorphism on $\clm$ preserving 
the state $\phi_0$, modular elements $\clj,\Delta^{1 \over 2}$ commutes with $\clj$ and in particular $\clj$ commutes with 
$F_kP$ for all $k \in \hat{H}$. Since $[v_Iv^*_J\Omega: |I|-|J|=k]=PF_k$ and $[\tilde{v}_I\tilde{v}^*_J\Omega: |I|-|J|=k ]=PF_k$ 
we get $vF_kP=\clj F_k P=F_kP$ as $\clj$ commutes with $F_k$. 

\vsp 
Fix any $z \in H$. By taking action of $\beta_z$ on both side of the relation $vv_kv^*=\clj \tilde{v}_k \clj$, we have 
$vv_kv^* = \bar{z}^2  \clj \tilde{v}_k \clj  = \bar{z}^2 vv_kv^*$. Thus $z^2v_k=v_k$ for all 
$1 \le k \le d$. Since $\sum_k v_kv_k^*=1$, we have $z^2=1$.   

\vsp 
The last statement (d) follows by a routine calculation as shown below for a special vectors.
$$<\Omega,\pi(s_Is^*_J\tilde{s}_{I'}\tilde{s}^*_{J'}\Omega>$$
$$=<\Omega,v_Iv_J^*\tilde{v}_{I'}\tilde{v}_{J'}\Omega>$$
$$=<\Omega,\clj_v \tilde{v}_I\tilde{v}_J^*v_{I'}v^*_{J'} \clj _v\Omega> $$ 
( as $\clj_v v_i \clj_v = \tilde{v}_i$)
$$=< \tilde{v}_I\tilde{v}_J^*v_{I'}v^*_{J'}\Omega,\Omega>$$ ($\clj_v$ being anti-unitary )
$$=<\pi(s_{I'}s^*_{J'}\tilde{s}_I\tilde{s}^*_J)\Omega,\Omega>$$ 

For anti-unitary relation involving more general vectors, we use Cuntz relations and the 
above special cases. The statement is obvious as $\clj_v$ is anti-linear. This completes the proof. \qed

\vsp 
We set an anti-linear $*$-automorphism $\clj_v: \clo_d \otimes \tilde{\clo}_d \raro  \clo_d \otimes \tilde{\clo}_d$ 
defined by 
$\clj_v(s_Is^*_J \otimes \tilde{s}_{I'}\tilde{s}^*_{J'}) = s_{I'}s^*_{J'} \otimes \tilde{s}_I\tilde{s}^*_{J}$
for $|I|,|J|,|I'|,|J'| < \infty$ by extending anti-linearly.

\vsp 
We say a state $\psi$ on $\clo_d \otimes \tilde{\clo}_d$ 
is reflection positive if $\psi(\clj_v(x)x) \ge 0$ for all $x \in \clo_d$ and equality holds if and only if $x=0$. 
Similarly for a state $\omega$ on $\clb$ we define reflection positivity. Note that this notion extended to 
$\tilde{\clo}_d \otimes \clo_d$ is an abstract version of the concept ``reflection positivity" of a state on 
$\clb$ introduced in [FILS] for any involution (linear or conjugate linear ) taking element from future algebra 
to past algebra. However this notion is different from Lieb's spin flip reflection symmetric [FILS]. This hidden 
symmetry $v$ will play an important role to determine properties of $\omega$ (Section 5). Now we aim to make $v$ little more general primarily motivated with reflection symmetry with a twist introduced in [FILS]. 

\vsp 
To that end we fix any $g_0 \in U_d(\!C)$ so that $g_0^2=1$ and $\beta_{g_0}$ is the natural action on $\clo_d$ and $\tilde{\clo}_d$. We say $\omega$ is {\it lattice reflection symmetric with twist } $g_0$ if $\omega(\beta_{g_0}(r(x))=\omega(x)$ for all $x \in \clb$ where $r$ is the refection automorphism around ${-{1 \over 2}}$. So when 
$g_0=1$ we get back to our notion of lattice reflection symmetric. We fix now such a lattice reflection 
$g_0$-twisted factor state $\omega$. Since $\beta_{g_0}\beta_z=\beta_z \beta_{g_0}$ for all $z \in S^1$, 
by going along the same line as in Proposition 4.2, any extremal element in $\psi$ in $K_{\omega}$ will 
admit $\tilde{\psi}^{g_0}=\psi \circ \zeta$ where $\zeta=1$ or $\zeta = exp^{\pi i \over n }$ 
where $H=\{z \in S^1: z^n=1 \}$ and $\tilde{\psi}^{g_0} = \tilde{\psi} \beta_{g_0}$. Thus we can follow the same
steps that of Proposition 4.4 to have a modified statements in the proof of Proposition 4.4 with $v_k$ replaced by 
$\beta_{g_0}(v_k)$ for such a pure real state i.e. there exists unitary operators $u_{\zeta},w_{\zeta}$ on $\clk$ 
so that
$$u_{\zeta}\beta_{g_0}(v_k)u_{\zeta}^*=\beta_{\bar{\zeta}}(\tilde{v}_k)$$
$$w_{\zeta}v_kw_{\zeta}^* = \beta_{\bar{\zeta}}(\bar{v}_k) = \clj \beta_{\zeta}(v_k) \clj $$
where $\;\;u_{\zeta} \clj u_{\zeta}^*=\clj, \;\; 
\;w_{\zeta} \clj w_{\zeta}^*=\clj$ and
$u_{\zeta} \Delta^{1 \over 2}u_{\zeta}^*=w_{\zeta} \Delta^{1 \over 2}w^*_{\zeta}=\Delta^{-{1 \over 2}}$.

Thus
\be
u_{\zeta} w_{\zeta} v_k w_{\zeta}^* u_{\zeta}^* = u_{\zeta} \clj \beta_{\zeta}(v_k) \clj u_{\zeta}^* 
= \clj \beta_{\zeta}(u_{\zeta}v_ku_{\zeta}^*) \clj = \clj \beta_{\zeta}(\beta_{g_0}(\beta_{\bar{\zeta}}(\tilde{v}_k)))\clj= \clj \beta_{g_0}(\tilde{v}_k) \clj
\ee

We also compute that
\be
w_{\zeta}u_{\zeta}v_ku_{\zeta}^*w_{\zeta}^* = w_{\zeta} \beta_{\bar{\zeta}}(\beta_{g_0}(\tilde{v}_k)) w_{\zeta}^* 
=\beta_{\bar{\zeta}}(\beta_{g_0}(w_{\zeta}\tilde{v}_k w_{\zeta}^*)) 
=\clj \beta_{\zeta} \beta_{\bar{g_0}}(\beta_{\bar{\zeta}}(\tilde{v}_k)) \clj 
= \clj \beta_{\bar{g_0}}(\tilde{v}_k) \clj 
\ee

Thus taking $v_{g_0}=w_{\zeta}u_{\zeta}$, as $g_0^2=g_0$ we also have 
\be 
v_{g_0}\beta_{g_0}(v_k)v^*_{g_0} = \clj \tilde{v}_k \clj
\ee 
for all $1 \le k \le d$ where $v_{g_0}$ is an unitary operator 
commuting with $\Delta^{1 \over 2}$ and $\clj$. Unlike the twist free case, self-adjoint property 
of $v_{g_0}$ is not guaranteed in general. In fact we get from the following computation  
$$v_{g_0}^2v_k(v^*_{g_0})^2=v_{g_0} \clj \beta_{\bar{g_0}}(\tilde{v}_k) \clj v^*_{g_0}$$
$$=\clj \beta_{\bar{g_0}} (\clj \beta_{\bar{g_0}}(v_k) \clj)\clj=\beta_{\bar{g_0}g_0}(v_k)$$   
Thus $v_{g_0}$ is self adjoint if and only if $g_0=\bar{g_0}$ as $g_0^2=1$.  

\vsp 
Such a $\omega$ is called reflection positive if $\omega(\bar{\beta_g(r(x))}x) \ge 0 $ for 
all $x \in \clb_R$ where $r$ is the reflection around the lattice point 
${ 1 \over 2 }$ so that $r(\clb_R)=\clb_L$ and $\bar{x}$ stands for complex conjugation 
i.e. $\bar{x}=\clj_0 x \clj_0$ where $\clj_0:(z_1,..z_d)=(\bar{z}_1,..\bar{z}_d)$ with respect to a basis. 
Such an involution are included within the abstract framework of positive reflection symmetric with twist 
introduced in [FILS].

\vsp
\begin{thm} 
Let $\omega$ be a translation invariant, reflection symmetric with twist $g_0$, pure state on $\clb$ 
and $\psi$ be an extremal point $K_{\omega'}$ and $\pi$ as described as in Theorem 4.4. Then the following statements are 
true:

\NI (a) $\psi$ is reflection positive with twist $g_0$ on $\pi(\tilde{\clo_d} \otimes \clo_d)$ if and only if $v_{g_0}$ in 
Theorem 3.10 is equal to $1$ i.e. we have $\clj \tilde{v}_k \clj = \beta_{g_0}(v_k) $ for all $1 \le k \le d$; 

\NI (b) $\psi$ is reflection positive with twist on $\clb$ if and only if $\clj \tilde{v}_I\tilde{v}_J^* \clj = \beta_{g_0}(v_Iv_J^*)$ 
for all $|I|=|J| < \infty $; In such a case $v_{g_0}$ commutes with $P_0$ and $v_{g_0}P_0=P_0$ 
and $\tilde{\tau}(y)=\clj \tau(\clj y \clj )\clj$ for $y \in \clm_0' \subseteq \clb(\clk_0)$ where 
we recall $P_0=[\clm_0\Omega]$ and $\clk_0$ is the Hilbert subspace for $P_0 \clk$. 

\NI (c) $\Delta=I$ if and only if $v_{g_0}\beta_{g_0}(v_k)v^*_{g_0} = v_k^*,\; 1 \le k \le d$. In such a case $H$ is trivial 
and $\clm$ is finite type-I and spacial correlation functions of $\omega$ decays exponentially. Further if $\omega$ is reflection 
positive with twist $g_0$, then $v_{g_0}=1$.      
\end{thm} 

\vsp
\NI {\bf PROOF:} We will prove for $g_0=1$ as a proof for $g_0 \neq 1$ needs no difference except 
involving a twist action $\bar{g}_0$ on $\clb_L$ to accomodate conjugate linearity on Popescu 
elements of the map $\clj_v$. We recall from Theorem 3.8 that $P\pi(\clo_d)''P=\clm$ 
and $P\pi(\tilde{\clo}_d)''P = \tilde{\clm} \subseteq \clm'$ 
( we do not need equality here ) and  $P=E\tilde{E}$. Thus for any $x \in \clo_d$ we may write 
$$\omega(\clj_v(x) x ) = <\Omega, \bar{\pi}(\clj_v(x)) \pi(x) \Omega> $$
$$= <\Omega, P \pi(\clj_v(x))P\pi(x)P \Omega> = <\Omega, \clj_v P\pi(x)P \clj_v P\pi(x)P \Omega>$$
where we have used equality $\bar{\pi}(\clj_v(x)) = \clj_v \pi(x) \clj_v$ from Theorem 3.10.  
If $v=1$ i.e. $\clj_v=\clj$ we have $\omega(\clj x \clj x) \ge 0$ by the self-dual property of 
Tomita's positive cone $\overline{ \{ \clj a \clj a \Omega: a \in \clm \} }$ [BR1] and being a 
pointed cone equality holds if and only $x=0$. Thus $\omega$ is a reflection positive map on 
$\pi(\tilde{\clo}_d \otimes \clo_d)$.  

\vsp 
Conversely if $\omega$ is reflection positive on $\pi(\tilde{\clo}_d \otimes  \clo_d)''$ we have 
$<\Omega,y \clj_v y \clj_v \Omega> \ge 0$ where $y \in \clm = P\pi(\clo_d)''P$. Since $v$ commutes 
with $\clj$ and $\Delta^{1 \over 2}$ we may rewrite $<\Omega, y v \clj y \Omega> = 
<y^*\Omega, v \Delta^{1 \over 2} y^*\Omega>  \ge 0$ i.e. $v\Delta^{1 \over 2}$ is a non-negative operator. 
Since $\Delta^{-{1 \over 2}}$ is also a non-negative operator commuting with $v\Delta^{1 \over 2}$, 
we conclude that $v$ is a non-negative operator. $v$ being unitary we conclude that $v=1$.  

\vsp 
Proof for (b) follows the same route that of (a) replacing the role of $\clm$ and $\tilde{\clm}$ by $\clm_0$ and 
$\tilde{\clm}_0$ respectively.  

\vsp 
We will deal with the non-trivial part of (c). Assume $v_{g_0}\beta_{g_0}(v_k)v^*_{g_0}=v_k^*$ for all $1 \le k \le d$. 
So $\Delta$ is affiliated to $\clm'$. As $\clj \Delta{1 \over 2} \clj =\Delta^{-{1 \over 2}}$, $\Delta$ is also affiliated 
to $\clm$. Hence $\Delta=I$ as $\clm$ is a factor and $\Delta \Omega=\Omega$. 

\vsp 
In general $\omega$ being a pure state $\clm$ is either a type-I or type-III factor [Mo3, Theorem 3.4]. Thus 
we conclude that $\clm$ is a finite type-I factor if $\Delta=1$ ( i.e. $\phi_0$ is a tracial state 
on $\clm$  ). This completes the proof of the first part of (c). 

\vsp
The last part of (c) is rather elementary. We note that purity of $\omega$ ensures that the point spectrum of 
the self-adjoint contractive operator $T$, defined by $Tx\Omega=\tau(x)\Omega$ on the KMS Hilbert space, in the unit 
circle is trivial i.e. $\{ z \in S^1: Tf=zf \mbox{for some non zero } f \in \clk \}$ is the trivial set $\{ 1 \}$ 
( as an consequence of strong mixing property ). Thus $T$ being a contractive matrix on a finite dimensional 
Hilbert space, the spectral radius of $T-|\Omega><\Omega|$ is $\alpha$ for some $\alpha < 1$. Now we use 
Proposition 3.1 for any $X_l \in \clb_L $ and $X_r \in \clb_R$ to verify the following
$$e^{\delta k }|\omega(X_l \theta_k(X_r)) -\omega(X_l)\omega(X_r)| $$
$$= e^{\delta k}|\phi_0(\clj_v x_l \clj_v \tau_k(x_r))- \phi_0(x_l)\phi_0(x_r)| \raro 0$$
as $k \raro \infty $ for any $\delta >0$ so that $e^{\delta}\alpha < 1$ where $\clj x_l \clj = PX_lP$ and $x_r=PX_rP$ for
some $x_l,x_r \in \clm$. As $\alpha < 1$ such a $\delta >0$ exists. This completes the proof for (c) as last 
statement follows from (b). \qed

\section{ Translation invariant twisted reflection positive pure state and it's split property: }

\vsp
Let $\omega$ be a translation invariant real lattice symmetric with a twist $g_0$ pure state on $\clb$ as in Theorem 4.5. We fix an extremal element $\psi \in K_{\omega'}$ so that $\bar{\psi}=\tilde{\psi}= \psi \beta_{\zeta}$ and consider the Popescu elements $(\clk,\clm,v_i,\Omega)$ as in Theorem 4.5. $P$ being the support projection of a factor state $\psi$ we have $\clm = P\pi(\clo_d)''P = \{ v_k,v^*_k: 1 \le k \le d \}''$ ( Proposition 2.4 ). So the dual Popescu elements $(\clk,\clm',\tilde{v}_k,1 \le k \le d, \Omega)$ satisfy the relation $v\tilde{v}_kv= \clj \beta_{g_0}(v_k) \clj$ ( recall that the factor $\zeta$ won't show up as two symmetry will kill each other as 
given in Theorem 4.5 ) for all $1 \le k \le d$. 

\vsp
We quickly recall as $\clm_0$ is the $\{\beta_z:\;z \in H \}$ invariant elements of $\clm ( =P\pi(\clo_d)''P )$, 
the norm one projection $x \raro \int_{z \in H}\beta_z(x)dz$ from $\clm$ onto $\clm_0$ preserves the 
faithful normal state $\phi_0$. So by Takesaki's theorem modular group associated with $\phi_0$ preserves 
$\clb_0$. Further since $\beta_z(\tau(x))=\tau(\beta_z(x))$ for all $x \in \clm$, the restriction of the 
completely positive map $\tau(x)=\sum_k v_kxv_k^*$ to $\clm_0$ is a well defined map on $\clm_0$. Hence the 
completely positive map $\tau(x)=\sum_k v_kxv_k^*$ on $\clm_0$ is also KMS symmetric modulo an unitaryconjugation by $v$ 
i.e. $$<<x,\tau(y)>>=<<\tau_v(x),y>>$$
where $x,y \in \clm_0$ and $<<x,y>>=\phi_0(x^*\sigma_{i \over 2}(y))$ and $(\sigma_t)$ is the modular automorphism
group on $\clm_0$ associated with $\phi_0$ and $[\clm_0\Omega]=P_0$ where $P_0=PF_0$ and $\tau_v(x)=v^*\tau(vxv^*)v$
for all $x \in \clm_0$. Thus $\tau_v=\tau$ if and only if $\omega$ is reflection positive on $\clb$ with 
twist $g_0$ (Theorem 4.5).    

\vsp 
However the inclusion $\clm_0 \subseteq \clm$ need not be an equality in general unless $H$ is trivial. 
The unique ground state of $XY$ model in absence of magnetic field give rise to a non-split 
translation invariant real lattice symmetric pure state $\omega$ and further $H=\{1,-1\}$ (see next section).

\vsp
We now fix a translation invariant real lattice symmetric pure state $\omega$ which is also reflection positive with a twist $g_0$ on $\clb$ and explore KMS-symmetric property of $(\clm_0,\tau,\phi_0)$ and the extended Tomita's conjugation operator $\clj$ on $\clh \otimes_{\clk} \tilde{\clh}$ defined in Theorem 4.5 to study the 
relation between split property and exponential decaying property of spacial correlation functions of $\omega$. 

\vsp
For any fix $n \ge 1$ let $Q \in \pi(\clb_{[-k+1,k]})$. We write
$$Q = \sum_{|I|=|J|=|I'|=|J'|=n} q(I',J'|I,J)\beta_{g_0}(\tilde{S}_{I'}\tilde{S}^*_{J'})S_IS^*_J $$
and $q$ be the matrix $q=((q(I',J'|I,J) ))$ of order $d^{2n} \times d^{2n}$. \qed 

\vsp
\begin{pro}
The matrix norm of $q$ is equal to operator norm of $Q$ in $\pi(\clb_{[-n+1,n]})$.
\end{pro}

\vsp
\NI {\bf PROOF: } $C^*$ completions of $\pi(\tilde{\mbox{UHF}}_d \otimes \mbox{UHF}_d)$ is isomorphic 
to $\clb$. Thus we note that the operator norm of $Q$ is equal to the matrix norm of $\hat{q}$
where $\hat{q}=((\hat{q}(I',I|J',J) ))$ is a $d^{2n} \times d^{2n}$ matrix with $\hat{q}(I',I|J',J)=
q(I',J'|I,J)$. However the map $L(q)=\hat{q}$ is linear and identity preserving. Moreover $L^2(q)=q$.
Thus $||L||=1$. Hence $||q||=||\hat{q}||$. This completes the proof. \qed

\vsp
\begin{pro}
Let $\omega$ be a translation invariant real lattice symmetric pure state on UHF$_d$
$\otimes_{\IZ}M_d$ . Then there exists an extremal point $\psi \in K_{\omega'}$ so that
$\psi \beta_{\zeta}=\tilde{\psi}=\bar{\psi}$ where $\zeta \in \{1, \mbox{exp}^{i \pi \over 2} \}$ and 
the associated Popescu systems $(\clh,S_k,\;1 \le k \le d,\Omega)$ and $(\clh,\tilde{S}_k,\;1 \le k \le d,\Omega)$ described in 
Proposition 3.2 and Theorem 4.5 satisfies the following:

\NI (a) For any $n \ge 1$ and $Q \in \pi(\clb_{[-n+1,n]})$ we write
$$Q=\sum_{|I'|=|J'|=|I|=|J|=n} q(I',J'|I,J)\beta_{g_0}(\tilde{S}^*_{I'}\tilde{S}^*_{J'})S^*_IS_J$$ and
set a notation for simplicity as $$\hat{\theta}_k(Q)= \sum_{|I|=|J|=|I'|=|J'|=n} q(I',J'|I,J)
\beta_{g_0}(\tilde{S}_{I'}\tilde{S}^*_{J'})\Lambda^{2k}(S_IS^*_J).$$
Then $\hat{\theta}_k(Q) \in \clb_{(-\infty,-k] \bigcup [k+1,\infty)}$.

\NI (b) $Q=\clj_{g_0} Q \clj_{g_0}$ if and only if $q(I',J'|I,J) = \overline{q(I,J|I',J')};$

\NI (c) If the matrix $q=(( q(I',J'|I,J) ))$ is non-negative then there exists a matrix $b=(( b(I',J'|I,J) ))$
so that $q=b^*b$ and then
$$q=PQP=\sum_{|K|=|K'|=n} \clj_v x_{K,K'} \clj_v x_{K,K'}$$
where $x_{K,K'} = \sum_{I,J:\;|I|=|J|=n} b(K,K'|I,J)v_Iv^*_J \in \clm_0$
\vsp
\NI (d) In such a case i.e. if $Q=\clj Q \clj$ the following holds:

\NI (i) $\omega(Q)=\sum_{|K|=|K'|=n} \phi_0(\clj_v x_{K,K'} \clj_v x_{K,K'})$

\NI (ii) $\omega(\hat{\theta}_{2k}(Q))=\sum_{|K|=|K'|=n}\phi_0(\clj_v x_{K,K'} \clj_v \tau_{2k}(x_{K,K'})).$

\end{pro} 

\vsp
\NI {\bf PROOF: } Since the elements $\beta_{g_0}(\tilde{S}_{I'}\tilde{S}^*_{J'})S^*_IS_J: |I|=|J|=|I'|=|J'|=n$ form
an linear independent basis for $\pi(\clb_{[-n+1,n]})$, (a) follows. (b) is also a simple consequence of linear 
independence of the basis elements and the relation $\clj \beta_{g_0}(\tilde{S}_{I'}\tilde{S}^*_{J'})S_IS^*_J \clj = 
S_{I'}S^*_{J'} \beta_{g_0}(\tilde{S}_I\tilde{S}^*_J) $ as described in Theorem 4.5.

\vsp 
For (c) we write 
$$Q=\sum_{|K|=|K'|=n} \clj \beta_{g_0}(Q_{K,K'}) \clj Q_{K,K'}$$ 
where $Q_{K.K'}= \sum_{I,J:\;|I|=|J|=n} b(K,K'|I,J)S_IS^*_J$. $\omega$ being pure we have ( Theorem 3.6)
$P=E \tilde{E}$ where $E$ and $\tilde{E}$ are support projection of $\psi$ in $\pi(\clo_d)''$
and $\pi(\tilde{\clo}_d)''$ respectively. So for any $X \in \pi(\clo_d)''$ and $Y \in \pi(\tilde{\clo}_d)''$
we have $PXYP=\tilde{E}EXY\tilde{E}E= \tilde{E}EYE\tilde{E}X\tilde{E}E=PXPYP$. Thus (c) follows as $\omega(Q)
= \phi_0(q)$ by Proposition 3.1 (b) and Theorem 4.5 as $\omega$ is reflection symmetry with twist $g_0$. 
For (d) we use (a) and (c). This completes the proof. \qed

\vsp
\begin{pro} 
Let $\omega$, a translation invariant pure state on $\clb$, be in detailed balance and 
reflection positive with a twist $g_0$. Then the following are equivalent:

\NI (a) $\omega$ is decaying exponentially.

\NI (b) The spectrum of $T-|\Omega><\Omega|$ is a subset of $[-\alpha,\alpha]$ for some $0 \le \alpha < 1$
where $T$ is the self-adjoint contractive operator defined by
$$Tx\Omega= \tau(x)\Omega,\;\;x \in \clm_0$$
on the KMS-Hilbert space $<<x,y>>=\phi_0(x^*\sigma_{i \over 2}(y)>>$.
\end{pro} 

\vsp
\NI {\bf PROOF: } Since $T^k x \Omega = \tau_k(x)\Omega$ for $x \in \clm_0$ and
for any $L \in \clb_L$ and $R \in \clb_R$ we have
$\omega(L \theta_k(R)) = \phi_0(\clj y \clj \tau_k(x)) = << y,T^kx >>$
where $x=P\pi(R)P$ and $y=\clj P\pi(L)P \clj$ are elements in $\clm_0$. Since $P\pi(\clb_R)''P = \clm_0$ and
$P\pi(\clb_L)''P = \tilde{\clm}_0 = \clm_0'$ as $\tilde{\clm}=\clm'$ by Theorem 4.4, we conclude that (a) 
holds if and only if $e^{k \delta}|<f,T^kg>-<f,\Omega><\Omega,g>| \raro 0$ as $k \raro \infty$ for any vectors $f,g$ in a 
dense subset $\cld$ of the KMS Hilbert space.

\vsp
That (b) implies (a) is now obvious since $e^{k\delta}\alpha^k=(e^{\delta}\alpha)^k \raro 0$ whenever we choose
a $\delta > 0$ so that $e^{\delta}\alpha < 1$ where $\alpha < 1 $.
\vsp
For the converse suppose that (a) holds and $T^2-|\Omega><\Omega|$ is not bounded away from $1$. Since
$T^2-|\Omega><\Omega|$ is a positive self-adjoint contractive operator, for each $n \ge 1$, we find an unit 
vector $f_n$ in the Hilbert space so that $E_{[1-1/n,1]}f_n=f_n$ and $f_n \in \cld$, where $s \raro E_{[s,1]}$ 
is the spectral family of the positive self-adjoint operator $T^2-|\Omega><\Omega|$ and in order to ensure 
$f_n \in \cld$ we also note that $E_{[s,1]}\cld=\{E_{s,1]}f:\;\;f \in \cld \}$ is dense in $E_{[s,1]}$ 
for any $0 \le s \le 1$.  

\vsp
Thus by exponential decay there exists a $\delta > 0$ so that
$$e^{2k \delta} (1-{1 \over n})^k \le e^{2k \delta} \int_{[0,1]}s^k<f_n,dE_sf_n>=
e^{2k \delta }<f_n,[T^{2k}-|\Omega><\Omega|]f_n> \raro 0$$
as $k \raro \infty$ for each $n \ge 1$. Hence $e^{2\delta}(1-{1 \over n}) < 1$. Since $n$ is any
integer, we have $e^{2 \delta} \le 1$. This contradicts that $\delta > 0$.
This completes the proof. \qed

\vsp
For the simplicity of notation we take in the following $g_0=1$. Now we are set to state our main result in this section. 
For any $Q \in \pi(\clb)$ we set $\clj(Q)=\clj Q \clj $. 
Recall that $\clj^2=I$. Any element $Q={1 \over 2}(Q+\clj(Q))+{1 \over 2}(Q-\clj(Q))$ is a sum of an even element in
$\{ Q: \clj(Q)=Q \}$ and an odd element in $\{Q: \clj(Q)=-Q \}$. Moreover $iQ$ is an
even element if $Q$ is an odd element. Also note that $||Q_{even}|| \le ||Q||$ and $||Q_{odd}||
\le ||Q||$. Hence it is enough if we verify (1) for all even elements for split property.
We fix any $n \ge 1$ and an even element $Q \in \clb_{[-n+1,n]}$. We write as in Proposition 5.2
$Q=\sum_{|I'|=|J'|=|I|=|J|=n}q(I',J'|I,J)\tilde{S}^*_{I'}\tilde{S}_{J'}S^*_IS_J$. The matrix
$q=(q(I',J'|I,J)$ is symmetric and thus $q=q_+-q_-$ where $q_+$ and $q_-$ are the unique
non-negative matrix contributing it's positive and negative parts of $q$. Hence $||q_+|| \le ||q||$
and $||q_-|| \le ||q||$. We set a notation for simplicity that
$$\hat{\theta}_k(Q)= \sum_{|I|=|J|=|I'|=|J'|=n} q(J',I'|I,J)\tilde{\Lambda}^k(\tilde{S}_{I'}\tilde{S}^*_{J'})
\Lambda^k(S_IS^*_J)$$ which is an element in $\clb_{(-\infty,-k] \bigcup [k,\infty)}$ and by Proposition
5.2 (d) $$\omega(\hat{\theta}_k(Q))=\sum_{|K|=|K'|=n}\phi_0(\clj x_{K,K'} \clj \tau_{2k}(x_{K,K'}))$$ provided
$q=(q(I',J'|I,J)$ is positive, where $PQP=\sum_{|K|=|K'|=n}\clj x_{K,K'}\clj x_{K,K'}$ and $x_{K,K'}=
\sum_{I,J}b(K,K'|I,J)v_Iv^*_J$ and $q=b^*b$.
Thus in such a case we have by Proposition 5.2 (d) that
$$|\omega(\hat{\theta}_k(Q))-\omega_L\otimes \omega_R(\hat{\theta}_k(Q))| =
\sum_{|K|=|K'|=n}\phi_0(\clj x_{K,K'} \clj (\tau_{2k}-\phi_0)(x_{K,K'}))$$
$$=\sum_{|K|=|K'|=n}<< x_{K,K'}, ( T - |\Omega><\Omega| )^{2k} x_{K,K'}>> $$
$$\le \alpha^{2k} \sum_{|K|=|K'|=n}<< x_{K,K'}, x_{K,K'}>> $$
provided $||T-|\Omega><\Omega||| \le \alpha$ and so 
$$\le \alpha^{2k}\omega(Q) \le \alpha^{2k}||\hat{q}|| = \alpha^{2k}||q|| $$
In the last identity we have used Proposition 5.1.

\vsp
Hence for an arbitrary $Q$ for which $\clj(Q)=Q$ we have
$$ |\omega(\hat{\theta}_k(Q))-\omega_L\otimes \omega_R(\hat{\theta}_k(Q))| \le \alpha^{2k}(||q_+||+||q_-||) \le 2 \alpha^{2k}||q|| = 2\alpha^{2k}||Q||$$
where in the last identity we have used once more Proposition 5.1.
Thus we have arrived at our main result by a well know criteria [BR1] on split property. 

\vsp
\begin{thm} 
Let $\omega$ be a translation invariant pure state. Let $\omega$ be also 
real ( with respect to a basis for $\!C^d$ ) and lattice symmetric with twist $g_0$.  If $\omega$ is 
reflection positive with twist a $g_0$ ( $g_0 \in U_d(\!C), g_0^2=1 )$ and the spatial  correlation 
function of $\omega$ decays exponentially then $\omega$ is split i.e. $\pi_{\omega}(\clb_R)''$ is a 
type-I factor.  
\end{thm}

\section{ Spontaneous symmetry breaking in quantum spin chain }

\vsp 
We aim to investigate the methodology developed in section 3 to study properties of the ground 
states of a translation invariant Hamiltonian for one lattice dimensional quantum spin chain 
$\clb=\otimes_{\IZ}M_d$. 

\vsp
Let $G$ be a compact group and $g \raro v(g)$ be a $d-$dimensional unitary representation of $G$. By $\gamma_g$ we denote
the product action of $G$ on the infinite tensor product $\clb$ induced by $v(g)$,
$$\gamma_g(Q)=(..\otimes v(g) \otimes v(g)\otimes v(g)...)Q(...\otimes v(g)^*\otimes v(g)^*\otimes v(g)^*...)$$
for any $Q \in \clb$. We recall now that the canonical action of the group $U(d)$ of $d \times d$ matrices on
$\clo_d$ is given by
$$\beta_{v(g)}(s_j)=\sum_{1 \le i \le d} s_i v(g)^i_j $$
and thus
$$\beta_{v(g)}(s^*_j) = \sum_{1 \le i \le d} \bar{v(g)^i_j} s^*_i $$

\vsp
Note that $v(g)|e_i><e_j|v(g)^*=|v(g)e_i><v(g)e_j| = \sum_{k,l} v(g)^l_i \bar{v(g)}^k_j|e_l><e_k|$, where $e_1,..,e_d$
are the standard basis for $\IC^d$. Identifying $|e_i><e_j|$ with $s_is^*_j$ we verify that on $\clb_R$ the gauge action
$\beta_{v(g)}$ of the Cuntz algebra $\clo_d$ and $\gamma_g$ coincide i.e. $\gamma_g(Q)=\beta_{v(g)}(Q)$
for all $Q \in \clb_R$.

\vsp
\begin{pro} 
Let $\omega$ be a translation invariant factor state on $\clb$. Suppose that
$\omega$ is $G-$invariant,
$$\omega(\gamma_g(Q))=\omega(Q) \; \mbox{for all } g \in G \mbox{ and any } Q \in \clb. $$
Let $\psi$ be an extremal point in $K_{\omega'}$ and $(\clk,\clm,v_k,\;1 \le k \le d,\phi_0)$ be the Popescu system
associated with $(\clh, S_i=\pi(s_i),\Omega)$ described as in Proposition 2.4. Then we have the following:

\NI (a) There exists an unitaryrepresentation $g \raro U(g)$ in $\clb(\clh)$ and a representation $g \raro \zeta(g)$ 
so that
\be
U(g)S_iU(g)^*= \zeta(g) \beta_{v(g)}(S_i),\;1 \le i \le d
\ee
for all $g \in G$ and

\NI (b) There exists an unitaryrepresentation $g \raro u(g)$ in $\clb(\clk)$ so that $u(g)\clm u(g)^*=\clm$ for
all $g \in G$ and $\phi_0(u(g)xu(g)^*)=\phi_0(x)$ for all $x \in \clm$. Furthermore the operator
$V^*=(v^*_1,..,v^*_d)^{tr}: \clk \raro \IC^d \otimes \clk $ is an isometry which intertwines the representation
of $G$,
\be
( \zeta(g) v(g) \otimes u(g) )V^*=V^*u(g)
\ee
for all $g \in G$, where $g \raro \zeta(g)$ is the representation of $G$ in $U(1)$.

\NI (c) $\clj u(g) \clj=u(g)$ and $\Delta^{it}u(g)\Delta^{-it}=u(g)$ for all $g \in G$.
\end{pro}

\vsp
\NI {\bf PROOF: } $\omega$ being a factor state by Proposition 3.3, $H$ is a closed subgroup of $S^1$. Thus $H$ is either
$S^1$ or a finite cyclic subgroup. We also recall that $\lambda \beta_g = \beta_g \lambda$ for all $g \in G$ and $\omega$
being $G$-invariant we have $\psi \beta_g \in K_{\omega'}$ for all $\psi \in K_{\omega'}$ and $g \in G$.

\vsp
In case $H=S^1$, by Proposition 2.6 $K_{\omega'}$ is having a unique   element and thus by our starting remark we have $\psi 
\beta_g = \psi$ for the unique extremal element $\psi \in K_{\omega'}$. In such a case we define unitary operator 
$U(g)\pi(x)\Omega=\pi(\beta_g(x))\Omega$ and verify (a) with $\zeta(g)=1$ for all $g \in G$.

\vsp
Now we are left to deal with the more delicate case. Let $H = \{z: z^n=1 \}''$ for some $n \ge 1$. In such a case by Proposition
2.5 and Proposition 3.2 (a) we have $\pi(\clo^H_d)''=\pi(\mbox{UHF}_d)''$ and $\pi(\tilde{\clo}^H_d)''=\pi(\tilde{\mbox{UHF}}_d)''$.
Thus for any $0 \le k \le n-1$ orthogonal projection $F_k$ is spanned by the vectors
$\{\tilde{S}_{I'}\tilde{S}^*_{J'}S_IS^*_JS^*_K\Omega: |I'|=|J'|,|I|=|J|,\; \mbox{and}\; |K|=k \}.$
We set unitary operator $U(g)'$ on $F_k: k \ge 0$ by 
$$U(g)'\pi(\tilde{s}_{I'}\tilde{s}^*_{J'}s_Is^*_Js^*_K)\Omega =
\pi(\beta_{v(g)}(\tilde{s}_{I'}\tilde{s}^*_{J'}s_Is^*_Js^*_K))\Omega$$
where $|I'|=|J'|,|I|=|J|$ and $|K|=k$. It is a routine work to check that $U(g)'$ is indeed an inner product preserving map
on the total vectors in $E_k$. The family $\{ F_k: 0 \le k \le n-1 \}$ being an orthogonal family projection with $\sum_k F_k= I$,
$U(g)'$ extends uniquely to an unitary operator on $\tilde{\clh} \otimes_{\clk} \clh$. It is obvious by our
construction that $g \raro U(g)'$ is a representation of $G$ in $\clh \otimes_{\clk} \tilde{\clh}$.

\vsp
For each $g \in G$ the Popescu element $(\clh,\beta_{v(g)}(S_k),1 \le k \le d,\;\Omega)$ determines an extremal point
$\psi_g \in K_{\omega'}$ and thus by Proposition 2.6 there exists a complex number $\zeta(g)$ with modulus $1$
so that $\psi_{g}=\psi \beta_{\zeta(g)}$. Note that for another such a choice $\zeta'(g)$, we have $\bar{\zeta(g)}\zeta(g)' \in H$.
As $H$ is a finite cyclic subgroup of $S^1$, we have a unique   choice once we take $\zeta(g)$ to be an element in the group
$S^1 / H$ which we identify with $S^1$. That $g \raro \zeta(g)$ is a representation of $G$ in $S^1=\{z \in \IC: |z|=1 \}$
follows as the choice in $S^1 / H$ of $\zeta(g)$ is unique. Hence there exists an unitary operator $U(g)$ and a representation
$g \raro \zeta(g)$ in $S^1$ so that
$$U(g)\Omega=\Omega,\;U(g)S_iU(g)^* = \zeta(g)\beta_{v(g)}(S_i)$$
for all $1 \le i \le d$. Thus 
$U(g) = \sum_k \zeta(g)^k U(g)'E_k$ for all $g \in G$ as their actions on any typical vector
$S_IS^*_JS_K\Omega,\;|I|=|J|,|K|=k < \infty $ 
are same. Both $g \raro U'(g)$ and $g \raro \zeta(g)$ being representations of $G$, we conclude that
$g \raro U(g)$ is an unitary representation of $G$. 

\vsp
The above covariance relation ensures that $U(g)\pi(\clo_d)''U(g)^* = \pi(\clo_d)''$ for all $g \in G$ and thus also
$U(g)\pi(\clo_d)'U(g)^*=\pi(\clo_d)'$ for all $g \in G$. Now it is also routine work to check that $U(g) \clf U(g)^*=\clf$, where we call 
$\clf=\{\pi(\clo_d)''\Omega]$ and $U(g)\cle U(g)^*=\cle$ where $\cle=[\pi(\clo_d)'\Omega]$ is the 
support projection of the state $\psi$ in $\pi(\clo_d)''$. Hence the support projection $P=\cle \clf$ of the state $\psi$ in the cyclic subspace $\clf$ is
also $G$ invariant i.e. $U(g)PU(g)^*=P$ for all $g \in G$. Thus we define $g \raro u(g)=PU(g)P$ an unitary representation of $g$ in $\clk$. Hence we have 
\be 
u(g)v_ju(g)^*=\zeta(g) \beta_{v(g)}(v_j) = \zeta(g)v_i v(g)^i_j
\ee 
for all $1 \le k \le d$. By taking adjoint we get
$u(g)v^*_ju(g)^*=\bar{\zeta(g)}\bar{v(g)^i_j} v^*_i$ for all $1 \le j \le d$.

\vsp
We are now left to prove (c). To that end we first verify that $S_0u(g)=u(g)S_0$ as their actions on
any typical vector $v_Iv^*_J\Omega$ are same, where $S_0x\Omega=x^*\Omega$ for $x \in \clm$. Hence by
uniqueness of the polar decomposition we conclude that (c) holds.

\vsp 
\begin{pro} 
Let $\omega$ be a $G$-invariant translation invariant factor state on $\clb$ as in Proposition 4.1 and 
$\omega$ be also pure real (with respect to an orthonormal  basis for $\!C^d$ ) and lattice symmetric with a twist $g_0$. We fix an extremal 
element $\psi \in K_{\omega}$ so that $\tilde{\psi}^{g_0} = \bar{\psi} = \psi \beta_{\zeta}$ where $\zeta=1$ or $exp^{i \pi\over n }$ as in 
Theorem 4.5. Let the associated family $\{v_k: 1 \le k \le d \}$ of operators defined in Theorem 3.10 be linearly independent (i.e. $\sum_k c_kv_k=0$ 
if and only if $c_k=0$ for all $1 \le k \le d$ ). Then $g_0$ intertwines the representation $g \raro \zeta(g)v(g)$ with it's complex conjugate matrix representation with respect to the orthonormal basis $(e_i)$ in $\IC^d$ if and only if $v_{g_0}$ commutes with $u(g)$ for all $g \in G$;

\end{pro}

\vsp
\NI {\bf PROOF: } For simplicity of notation in the proof we use $v$ for $v_{g_0}$. For 
(a) we fix an extremal element $\psi \in K_{\omega'}$ as described in Theorem 4.5 and consider the Popescu 
elements $\{v_j,1 \le j \le d \}$ satisfying 
$$v\beta_{g_0}(v_j)v^*= \clj \tilde{v}_j \clj$$ 
for all $1 \le j \le d$. 

\vsp 
We set temporary notation $$v^g=u(g)vu(g)^*$$ 
and $\overline{v(g)}$ for complex conjugation matrix of $v(g)$ i.e. $\overline{v(g)}^i_j = \overline{v^i_j(g)}$. 
We compute the following simple steps.

$$u(g)v\beta_{g_0}(v_j)v^*u(g)^*=\clj u(g)\tilde{v}_ju(g)^* \clj=\overline{\zeta(g)}\beta_{\overline{v(g)}}(\clj \tilde{v}_j \clj)$$
and so 
$$v^g \zeta(g) \beta_{v(g)}(\beta_{g_0}(v_j))(v^g)^* $$
$$=u(g)vu(g)^*u(g)\beta_{g_0}(v_j)u(g)^*u(g)v^*u(g)^*$$
$$=\overline{\zeta}(g) \beta_{\overline{v(g)}} (\clj \tilde{v}_j \clj ) $$
$$=\overline{\zeta}(g) \beta_{\overline{v(g)}}(v\beta_{g_0}(v_j)v^*)$$  
i.e. 
$$v^*v^g \zeta(g) \beta_{v(g)}(\beta_{g_0}(v_j))(v^g)^*v = 
\overline{\zeta(g)} \beta_{ \overline{v(g)}}(\beta_{g_0}(v_j))$$ 
 
\vsp 
Explicitly if $v$ commutes with $u(g)$, by the covariance relation we have
$$\sum_{1 \le j \le d}[\zeta(g)v^i_j(g)-\overline{\zeta(g)v^i_j(g)}]v_j=0$$ 
for all $1 \le i \le d$. By linear independence we conclude that there exists 
an orthonormal basis for $\IC^d$ so that the matrix representation $(\zeta(g)v^i_j(g))$ 
of $\zeta(g)v(g)$ with respect to the basis having real entries. Conversely we have $v^*v_g$ 
commuting with all $v_i$ i.e. $v^*v_g \in \clm'$. But $v,u(g)$ and so $v^*v^g$ commutes with $\clj$ 
and thus $v^*v^g \in \clm$ and as $v^*v^g\Omega=\Omega$ we have $v^*v^g=1$ by factor property of $\clm$. 
For $g_0 \neq 1$ we get the relation $g_0 v(g) g_0^* = \overline{v}(g)$ where conjugation with 
respect to the basis $(e_i)$ as a necessary and sufficient condition for commuting property of $u(g)$ and $v$.  
\qed 

\vsp
We are left to deal with another class of example. Let $\omega$ be a translation invariant pure state on 
$\clb=\otimes_{\IZ}M_2$. If $\omega$ is $G=U(1) \subseteq SU(2)$ invariant then by a Theorem [Ma2] $\omega$ 
is either a product state or a non-split state. The following theorem says more when $\omega$ is also real
and lattice symmetric. The following theorem is originated from reviewer's remark on an earlier version.    

\vsp
\begin{thm} 
Let $\omega$ be a translation invariant pure state on $\clb=\otimes_{\IZ}M_2$ which 
also be reflection symmetric with a twist $g_0$ and real with respect to a basis $(e_i)$. If $\omega$ is also 
$U(1) \subseteq SU(2)$ invariant and reflection positive on $\clb$ with twist then following holds: 

\NI (a) $\omega$ is a non-split state if and only if $H=\{1, -1\}$;  

\NI (b) $\omega$ is a split state ( hence product state ) if and only if $H=\{1 \}$.
\end{thm}

\vsp
\NI {\bf PROOF: } Assume for the time being that $\{v_0,v_1 \}$ are linearly independent. In such a case by Proposition 6.1 
there exists a representation $z \raro u_z$ in $\clk$ of $U(1)$ so that 
$u_zv_0u_z^*=\zeta(z)zv_0$ and $u_zv_1u_z^*=\zeta(z)\overline{z}v_1$  
where $z \raro \zeta(z)=z^k$ for some $k \in \IZ$. $\omega$ being lattice symmetric 
and real with respect to a basis, by Theorem 5.4 we also have $H \subseteq \{1,-1\}$.  

\vsp 
Once we also assume $\omega$ to reflection positive with twist, $z \raro \zeta(z)v(z)$ has a real representation of 
$U(1) \subseteq SU(2)$ in a basis and $v$ commutes with $u_z$ for all $z \in H$ by Proposition 5.2. Thus $((z^kv^i_j(z)))=((z^{-k}\overline{v^i_j(z)}))$ for all $z \in S^1$. Hence $z^{2k}I = ((v^i_j(z))^{-1} ((\overline{v^i_j(z)}))$ where 
the right hand side is a multiple of two elements in $SU(2)$. By taking determinants of both the side we conclude that $z^{4k}=1$ 
for all $z \in S^1$. Hence $k=0$. Now by $U(1) \subseteq SU(2)$ invariance we also check that $\psi(s_Is^*_J)=0$ if $|I|+|J|$ is 
not a multiple of $2$. Further as $\psi(\beta_z(s_Is^*_J)) = z^{|I|-|J|}\psi(s_Is_J^*)$ and $|I|-|J|=|I|+|J|-2|J|$ is a multiple 
of $2$ whenever $|I|+|J|$ is an even number, we check that $-1 \in H$. Thus we conclude 
that $\{-1,1\} \subseteq H$ whenever $\{v_0,v_1 \}$ are linearly independent. In such a case 
$\omega$ is a non-split state.     

\vsp
Now we consider the case where $\{v_0, v_1\}$ are linearly dependent. We assume without loss of generality that 
$v_0 \ne 0$ and $v_1=\alpha v_0$. $v_0v_0^* + v_1v_1^*=1$ ensures that $(1+|\alpha|^2)v_0v_0^*=1$. 
Thus $v_iv_j^*=\phi_0(v_iv_j^*)$ for all $0 \le i,j \le 1$. Hence $\omega(|e_{i_1}><e_{j_1}| \otimes ...\otimes |e_{i_n}><e_{j_n}|)
=\phi_0(v_Iv^*_J) = \phi_0(v_{i_1}v^*_{j_1})...\phi_0(v_{i_n}v^*_{j_n})$. This clearly shows that $\omega$ is a product state. 
$\omega$ being pure, we conclude that there is an extremal point $\psi$ so that the associated Popescu elements are given by 
$v_0=1$ and $v_1=0$. In such a case $\omega$ is a split state and $H=\{ 1 \}$. 

\vsp
Thus we can sum up from the above argument that such a state $\omega$ is split ( non-split ) if and only if $H = \{ 1 \}$ ($H=\{1,-1\}$).  \qed 

\vsp
\begin{thm} 
Let $G$ be a simply connected compact Lie group and $g \raro v(g)$ be an irreducible representation on 
$\IC^d$ so that the invariance vectors of $\IC^d \otimes \IC^d$ with respect to the representation $g \raro \overline{v(g)} \otimes v(g)$ is one dimensional. Let $\omega$ be a translation and $\{\beta_g: g \in G \}$ 
invariant factor state on $\clb = \otimes_{\IZ}M_d$. Then the following holds:

\NI (a) The family $(v_k)$ is linearly independent and  
$\{v^*_k\Omega: 1 \le k \le d \}$ is a set of orthogonal 
vectors in $\clk$;   

\NI (b) The following statements are equivalent:

\NI (i) $\sum_k v_k^*v_k=1$;

\NI (ii) $\Delta=I$; 

\NI (iii) The action $g:x \raro u(g)xu(g)^* $ on $\clm$ is ergodic.

\vsp 
Further in such a case $\clm$ is a finite type-I$_n$ for some $n \ge 1$ and the unique normalized trace $\phi_0$ on 
$\clm$ is strongly mixing for $\tau: \clm \raro \clm$, where $(\clm,\tau,\phi_0)$ is defined as in Proposition 2.4; 
Further there exists a unique  representation of $G$, $g \raro u'(g) \in \clm$ so that
$$u'(g)v_ku'(g)^*= \sum_{1 \le j \le d} v^k_j(g)v_j$$
for all $1 \le k \le d$. The unique representation $g \raro u'(g)$ is irreducible i.e. 
$(\clm,\alpha_g,\phi_0)$ is $G$-ergodic. 
\end{thm} 

\vsp
\NI {\bf PROOF: } As a first step we explore the hypothesis that the invariant vectors in $\IC^d \otimes \IC^d$ of the 
representation $g \raro \overline{v(g)} \otimes v(g)$ of $G$ is one dimensional. By appealing to the assumption with 
invariant vectors $((\phi_0(v_i^*v_j)))$ and 
$((<\Omega,S_i^*S_j\Omega>))$ we have
$$\phi_0(v^*_iv_j)= \delta^i_j {\lambda \over d}$$
for all $1 \le i,j \le d$ and scaler $\lambda > 0$. That the scaler $\lambda$ is indeed non-zero follows 
by separating property of $\Omega$ ( otherwise we will have $v_k=0$ for all $k$ ).  In particular we get the 
vectors $\{v_i\Omega: 1 \le i \le d \}$ are linearly independent and so by separating property of $\Omega$ for 
$\clm$ (a) follows. (b) is immediate from (a) and Proposition 6.2.   

\vsp
Vectors $((<\Omega,S^*_iS_j\Omega>))$ and $((\phi_0(v_iv^*_j) ))$ are also invariant for the product 
representation $g \raro \overline{v(g)} \otimes v(g)$ and thus we also have
$$\phi_0(v_iv^*_j) = {1 \over d}\delta^i_j$$
for all $1 \le i,j \le d$, where we have used Popescu's relation $\sum_kv_kv_k^*=1$. Similarly we also 
have $\phi_0(v^*_iv_j)= { \lambda \over d }\delta^i_j$ for some $\lambda > 0$. Same is true for
$\phi_0(v_i\Delta^sv_j^*)$ as $\Delta$ commutes with $g \raro u(g)$ as the induced automorphism preserves 
the state $\phi_0$.  

\vsp 
Now we will prove the equivalence of those three statements. (ii) implies (i) is obvious. Now we will prove 
(i) implies (ii): (i) ensures that $\lambda=1$ and $\sum_k \tilde{v}^*_k\tilde{v}_k=1$. Thus a simple computation shows that
$$||\Delta^{1 \over 2}v_k^*\Omega-\Delta^{-{1 \over 2}}v_k^*\Omega||^2 $$
$$=<v_k^*\Omega,\Delta v_k^*\Omega>+<v^*_k\Omega,\Delta^{-1}v_k^*\Omega>-2<\Omega,v_kv_k^*\Omega>$$
$$=<v_k\Omega,v_k\Omega>+<\clj \Delta^{-{1 \over 2}} \tilde{v}_k^*\Omega, \clj \Delta^{-{1 \over 2}}\tilde{v}_k^*\Omega> - 2<\Omega,v_kv_k^*\Omega>$$
$$=<v_k\Omega,v_k\Omega>+<\tilde{v}_k\Omega, \tilde{v}_k\Omega> - 2<\Omega,v_kv_k^*\Omega>=0 $$
 
\vsp 
Hence by separating property of $\Omega$ for $\clm$ we conclude that 
$\Delta v^*_k \Delta^{-1} =v^*_k$ for all $ 1 \le k \le d$. So $\Delta$
is affiliated to $\clm'$. As $\clj \Delta \clj =\Delta^{-1}$, $\Delta$ is
also affiliated to $\clm$. Hence $\Delta=I$ as $\clm$ is a factor. This completes 
the proof for (i) implies (ii).

\vsp 
(iii) implies (i) follows as ergodic action ensures that $\clm$ is a type-I finite factor [Wa] and $\phi_0$ is 
the normalized  trace. Since $\sum_k v_k^*v_k$ is a $G$-invariant element, we conclude that the element is a scaler and 
hence by trace property of $\phi_0$ it is equal to $1$.  

\vsp
That irreducibility $g \raro u'(g)$ is equivalent to $G$-ergodicity of $(\clm,\alpha_g,\phi_0)$ follows from a more general result 
[BR1], however here one can verify easily as $\clm$ is a type-I factor. For the non-trivial part of last statement we will use once 
more the property (b) i.e. self-adjointness of the Popescu elements $\{v_k,\;1 \le k \le d \}$. To that end let $E$ be a $G$-
irreducible projection in $\clm$ and set $E'= \clj E \clj $. We set von-Neumann algebra $\clm_{E}=E'\clm E'$. 
$\clm_{E}$ is a type-I finite factor as $\clm$ is so and $u(g) \clm_E u(g)^* = \clm_E$ for all $g$ as $E'$ is $G$-invariant. 
Further $\clm_{E}$ is also $G$-irreducible as $E$ is an $G$-irreducible projection in $\clm$. The vector state $\phi_{E}(X) = 
<\Omega,X\Omega>$ on $\clm_E$ being $G$-invariant and irreducible, we conclude that $\phi_{E}$ is a scaler multiple of the unique 
normalized trace on $\clm_{E}$ ( see [HLS] for details ). 

\vsp 
That the non-normalized state $\psi_{E}$ on $\clo_d$ defined by 
$$\psi_{E}(s_Is_J^*)=<\Omega,E'v_Iv^*_JE'\Omega>\;\;|I|,|J| < \infty $$ 
is $\lambda$-invariant follows by the tracial property of $\phi_E$ on $\clm_{E}$ 
and by (i) we make the following computation as follows:

$$\psi_{E}(\lambda(s_Is_J^*))=\sum_{1 \le k \le d} <\Omega,E'v_k v_Iv^*_Jv^*_k E'\Omega>\;\;|I|,|J| < \infty $$ 
$$=\sum_{1 \le k \le d} <\Omega, E'v^*_k E'v_kE'v_Iv_J^*E'\Omega>$$
$$=\psi_{E}(s_Is_J^*)$$ 
where $\sum_k E'v^*_kE'v_kE'=E'$ by (i). Further $\lambda$-invariance of $\psi_{E}$ ensures that
$$<\Omega,E'xE'\Omega>=<\Omega,E'\tau(x)E'\Omega>$$
for all $x \in \clm$. Tomita's modular operator being trivial we also have $E'\Omega=E\Omega$ and so by duality relation  we have
$$<\Omega,\clj E \clj x\Omega> = <\Omega,\tilde{\tau}(\clj E \clj)x\Omega>$$ 
for all $x \in \clm$. Thus by cyclic and separating property of $\Omega$ for $\clm$, we get $\tilde{\tau}(\clj E \clj )
= \clj E \clj $. Now by ergodic property of $(\clm,\tau,\phi_0)$ ( See Proposition 2.4 (e) ) and so the property for the dual Markov map, we conclude that $\clj E \clj $ is either $0$ or $1$. This completes the proof for (i) implies (iii). 

\vsp 
Thus in such a case i.e. if any of the statement (i)-(iii) is true, $\phi_0$ is a tracial state on $\clm$. That it is the 
unique trace follows as $\clm$ is either a type-I finite factor or a type-II$_1$ factor. However $\omega$ being pure, $\clm$ can not be a type-II$_1$ factor [Mo3]. $\clm$ being a type-I factor and $G$ being a simply connected group, by a general result 
[Ki] any continuous action of the $G$ is implemented by an inner conjugation, i.e. there exists an unitary operator 
$g \raro u'(g) \in \clm$ so that $u(g)'x{u(g)'}^*=u(g)xu(g)^*$ for all $g \in G$. Thus we have $u'(g)v_ku'(g)^* = 
\sum_j v^k_j(g)v_j$ for all $1 \le k \le d$. Now for uniqueness let $g \raro u''(g) \in \clm$ be another such 
representation. Then $u''(g)u'(g)^* \in \clm \bigcap \clm'$ and $\clm$ being a factor, $\lambda(g)=u''(g)u'(g)^*$ 
is a scaler and as $u'(g) = \lambda(g) u''(g)$ and each one being a representation we also get $g \raro \lambda(g)$ 
is a representation and so $\lambda(g)=1$ as $G$ is simply connected. Hence uniqueness follows. \qed 

\vsp 
We are left to discuss few motivating examples for this abstract framework to study symmetry. For basic facts 
about thermodynamic limit, KMS states, ground states on quantum spin chain we refer to [Ru,BR2,EK].   
 
\vsp 
We consider the following standard ( irreducible ) representation of Lie algebra $su(2)$ in $\!C^2$

\ben
\sigma_x = \left (\begin{array}{llll} 0&,&\; 1 \\ 1&,&\;\;0 \\ 
\end{array} \right ),
\een
\ben
\sigma_y = \left (\begin{array}{llll} \;0&,&\;\; i \\ -i&,&\;\;0 \\
\end{array} \right ),
\een
\ben
\sigma_z = \left (\begin{array}{llll} 1&,&\;\; 0 \\  0&,& \;-1\\
\end{array} \right ).
\een
\ben
g_0 = \left (\begin{array}{llll} \;0&,&\;\; i \\ -i&,&\;\;0 \\
\end{array} \right ),
\een

Note also that $g_0^2=1$ and $ig_0 \in SU(2)$ and $g_0\sigma_xg^*_0=-\sigma_x$, $g_0\sigma_y\sigma_0^*=\sigma_y$
and $g_0 \sigma_z g_0^*=-\sigma_z$. By taking exponential of $i\sigma_x,i\sigma_y,i\sigma_z$, we conclude that 
$g_0$ inter-twins $g \raro v^i_j(g)$ with it's complex conjugate representation. The representation being irreducible, such an intertwining $g_0$ is unique modulo a phase factor. Since we also need $g_0^2=1$, we only have the choice given above or $-g_0$. Same is true if we work with any irreducible representation of $su(2)$ in $d=2s+1$ dimension and one such $g_0 \in U_d(\!C)$ exists which inter-twins the representation with it's complex conjugate. Now we aim to deal with two instructive examples. To that end we will now study the relation (33) for ${1 \over 2}$-integer quantum spin chain. 

\vsp 
\begin{thm} Let $\omega$ be reflection positive with twist $g_0$ as in Theorem 4.5 (b) and $d$ be $2$ with $g_0=\sigma_y$ given above. Then 
\be 
v_1 =  \epsilon \clj \tilde{v}_2 \clj 
\ee 
where $\epsilon$ is either $1$ or $-1$. Further if the family $\{v_k: 1 \le k \le 2 \}$ are linearly independent and 
$\omega$ is $G$-invariant where $g \raro v^i_j(g)$ is a representation taking values in $SU_2(\!C)$ then the 
range $(v^i_j(g)) \in S^1 \subseteq SU(2)$. There exists no $SU(2)$-invariant state $\omega$ on $\clb$ 
satisfying conditions for Theorem 4.5 (b) for $d=2$ with $g \raro v^i_j(g)$ irreducible. 
\end{thm}

\vsp 
\NI{ \bf PROOF: } By Theorem 4.5 (b) we have $v_{g_0} X v_{g_0}^*=X$ for all $x \in \clm_0 \vee \tilde{\clm}_0$  
\be 
v_{g_0}v_1v_{g_0}^*= i \clj \tilde{v}_2 \clj
\ee 
and
\be 
v_{g_0}v_2v_{g_0}^*=-i \clj \tilde{v}_1 \clj
\ee 
Since $v_{g_0}\clm v^*_{g_0}=\clm$ by separating property of $\Omega$ for $\clm$, we also claim that 
\be 
v_{g_0}\tilde{v}_1v_{g_0}^*= i \clj v_2 \clj
\ee 
which we verify as
$$v_{g_0}\tilde{v}^*_1v_{g_0}^*\Omega=v_{g_0}\tilde{v}^*_1\Omega= v_{g_0}v_1^*\Omega$$
$$=v_{g_0}v^*_1v_{g_0}^*\Omega =-i \clj \tilde{v}^*_2 \clj\Omega=-i\clj v_2^*\Omega= -i \clj v_2^* \clj \Omega$$
So by separating property we have 
$v_{g_0}\tilde{v}^*_1v_{g_0}^*=-i \clj v_2^* \clj$. 
Similarly we have the other relation 
$$v_{g_0}\tilde{v}_2v_{g_0}^*=-i \clj v_1 \clj$$
Further we recall $v_{g_0}$ commutes with $\clj$ and so 
\be 
v_{g_0}\clj \tilde{v}_2 \clj v_{g_0}^*= i v_1 
\ee
We also have $v_{g_0}^2v_k(v_{g_0}^2)^* = \beta_{-1}(v_k)$ for $1 \le k \le 2$. 
 
\vsp 
A simple computation now shows that 
$$v_{g_0}(v_1 + \clj \tilde{v}_2 \clj )v_{g_0}^*=i(v_1+\clj \tilde{v}_2 \clj)$$
and 
$$v_{g_0}(v^*_1-\clj \tilde{v}^*_2 \clj)v_{g_0}^*= -i \clj \tilde{v}^*_2\clj + i v_1^*$$ 
So we have 
$$v_{g_0}(v_1 + \clj \tilde{v}_2 \clj ) X (v^*_1-\clj \tilde{v}^*_2 \clj)v_{g_0}^*$$
$$= - (v_1 + \clj \tilde{v}_2\clj ) X (v^*_1-\clj \tilde{v}^*_2\clj) $$ 
where 
$X \in \clm_0 \vee \tilde{\clm}_0$. The element $(v_1 + \clj \tilde{v}_2\clj ) X (v^*_1-\clj \tilde{v}^*_2\clj)$
being $\{ \beta_z:z \in H \}$ invariant, is an element in $\clm_0 \vee \tilde{\clm}_0$ and thus  
$(v_1 + \clj \tilde{v}_2 \clj ) X (v^*_1-\clj \tilde{v}^*_2 \clj)=0$ 
for all $X \in \clm_0 \vee \tilde{\clm}_0$. 

\vsp 
Assume now that $v_1 - \clj \tilde{v}_2 \clj \neq 0$ and then action of $(v_1+\clj \tilde{v}_2 \clj)^*(v_1 + \clj \tilde{v}_2 \clj )$ on $[X (v^*_1-\clj \tilde{v}^*_2 \clj)(v^*_1-\clj \tilde{v}^*_2 \clj)^*\Omega]$ is equal to $0$ vector. By purity of $\omega$ we also have $[\clm_0 \vee \tilde{\clm}_0 f]=P_0$ for any $f \neq 0$ and $P_0f=f$. Note that $f=(v^*_1-\clj \tilde{v}^*_2 \clj)(v^*_1-\clj \tilde{v}^*_2 \clj)^*\Omega$ is $u_z$-invariant where 
$\beta_z(x)=u_zxu_z^*$ defined as in Proposition 3.2 and so $P_0f=f$. Thus by our assumption and separating 
property of $\Omega$, we have $f \neq 0$. Hence we conclude that 
$(v_1+\clj \tilde{v}_2 \clj )^*(v_1 + \clj \tilde{v}_2 \clj )P_0=0$. Since $P_0\Omega=\Omega$, 
by separating property for $\clm$ we conclude that $v_1+\clj \tilde{v}_2 \clj=0$. Interchanging 
the role of elements involve we conclude the first part of the result.  

\vsp 
By covariance relation (33) and linear independence of $\{v_k:1 \le k \le 2 \}$ we have now 
$v^1_1(g)= \overline{v^2_2(g)}$ and $\overline{v^1_2(g)}=v^2_1(g)$. $((v^i_j(g)))$ being an element 
in $SU(2)$ we also have $v^1_2(g)= - \overline{v^2_1(g)}$. Hence we conclude that $v^1_2(g)=v^2_1(g)=0$ 
for all $g \in G$. Thus $(v^i_j(g)) \in S^1$ as a subgroup of $SU(2)$. The last statement is now obvious 
once we use Proposition 6.3 to ensure linear independence hypothesis on the family $\{v_k: 1 \le k \le d \}$ 
since $G=SU(2)$ satisfies the hypothesis by Clebsch Gordan theory and so brings a contradiction 
if we claim to exist such a state which is also $SU(2)$ invariant with $g \raro v(g)$ irreducible. This
completes the proof. \qed  

\vsp
\NI {\bf XY model: } We consider the exactly solvable XY model. The Hamiltonian $H_{XY}$ of the XY model is
determined by the following prescription:

$$H_{XY}= J ( \sum_{j \in \IZ} \{ \sigma_x^{(j)}\sigma_x^{(j+1)}+\sigma_y^{(j)} \sigma_y^{(j+1)} \} - 
2 \lambda \sum_{j \in \IZ} \sigma_z^{(j)}),$$
where $\lambda$ is a real parameter stand for external magnetic field and $J$ is a non-zero real number, 
$\sigma_x^{(j)},\sigma_y^{(j)}$ and $\sigma_z^{(j)}$ are Pauli spin matrices at site $j$. It is well known 
[AMa] that ground state exists and unique. It is simple to verify that $\tilde{H}=H$ since we can rewrite 
$H_{XY}$ as sum over element of the form $\sigma_x^{(j-1)} \sigma_x^{(j)} + \sigma_y^{(j-1)} \sigma_y^{(j)}$. 
Since the transpose of $\sigma_x$ is itself, transpose of $\sigma_y$ is $-\sigma_y$ and transpose of 
$\sigma_z$ is itself, we also verify that $H_{XY}^t=H_{XY}$. Hence $H_{XY}$ is in detailed balance.  

\vsp 
For $J < 0$, it is also well known that for $|\lambda| \ge 1$ the unique ground state is a product state thus split state. 
On the other hand for $|\lambda| < 1$ the unique ground state is not a split state [Ma2 Theorem 4.3]. For $J > 0$ 
$H_{XY}$ is reflection symmetric with a twist $g_0$ which rotates an angle $\pi$ with respect $Y$-axis. Further 
by a general theorem [FILS] $\omega$ is also reflection positive with a twist $g_0$ when $J > 0$ and $\lambda =0$.    
Thus by Theorem 6.3 the unique ground state is a non split state and $H= \{1,-1 \}$. In such a case a simple application of Theorem 5.4 says that the correlation functions of the ground state does not decay exponentially. Theorem 6.5 gives 
functional relation between $v_1$ and $v_2$ and $\Delta$ is non-trivial. \qed       

\vsp
\NI {\bf XXX MODEL: } Here we consider the prime example where very little exact results were known. 
The Hamiltonian $H_{XXX}$ of the spin $s$ anti-ferromagnetic chain i.e. the Heisenberg's XXX model is 
determined by the following formula:
$$H_{XXX} = J \sum_{j \in \IZ} \{ S_x^{(j)}S_x^{(j+1)}+S_y^{(j)}S_y^{(j+1)} + 
S_z^{(j)} S_z^{(j+1)} \}$$
where $S_x^{(j)},S_y^{(j)}$ and $S_z^{(j)}$ are representation in $d=2s+1$ dimensional of Pauli spin matrices 
$\sigma_x,\sigma_y$ and $\sigma_z$ respectively at site $j$. Existence of ground 
state for XXX model follows from more general theory [BR vol-2]. Since $H_{XXX}$ can be rewritten 
as sum of elements of the form  $$\{ S_x^{(j-1)}S_x^{(j)}+S_y^{(j-1)}S_y^{(j)} +  
S_z^{(j-1)} S_z^{(j)}\}$$, it is simple to check that $\overline{H}_{XXX}=H_{XXX}$. We also claim that 
$H_{XXX}^t=H_{XXX}$. To that end we consider the space $V_d$ of homogeneous polynomials in two
complex variable with degree $m,\;m \ge 0$ i.e. $V_d$ is the space of functions of the form
$$f(z_1,z_2)=a_0z^d_1+a_1z_1^{d-1}z_2+...+a_dz_2^d $$ 
with $z_1,z_2 \in \IC$ and $a_i's $ are arbitrary complex constants. Thus $V_d$ is a $d$-dimensional complex
vector space. The $d-$dimensional irreducible representation 
$\pi_d$ of the Lie-algebra $su(2)$ is given by 
$$\pi_d(X)f= - { \partial f \over \partial z_1 } ( X_{11}z_1+X_{12}z_2 ) +  
{ \partial f \over \partial z_2 } ( X_{21}z_1+X_{22}z_2 )$$  
where $X$ in any element in Lie-algebra $su(2)$. It is simple to verify that the transpose of 
$S_x=\pi_d(\sigma_x)$ is itself, transpose of $S_y=\pi_d(\sigma_y)$ is $-S_y$ and transpose of 
$S_z=\pi_d(\sigma_z)$ is itself. Thus $H_{XXX}^t=H_{XXX}$ for any $d$. Further for $J > 0$, the 
unique positive temperature state (KMS state ) for $H_{XXX}$ is also reflection positive [FILS] with twist $g_0$ where $g_0$ is as in $XY$ model. Thus any limit point of KMS states as temperature goes to zero is also reflection positive. Thus for 
$J > 0$, if ground state is unique then it is also pure and reflection positive with twist $g_0$.  This brings a contradiction 
to our hypothesis on $\omega$ ( by the last statement in Theorem 6.5 ) that is real lattice symmetric positive with twist $g_0$. We can expect same result to hold for any $s$ with half-odd integer spin as 
$g_0$ acts non trivially and $\overline{g}_0g_0=-I$.       

\vsp 
Thus in particular for $J > 0$ ( anti-ferromagnet ) and $s$ is half-integer, there is a spontaneous symmetry breaking 
in the ground state. So we conclude that there are solutions other then the well known Bethe's ansatz [Be] solution. As an indirect consequence we conclude that Bethe's solution is not pure as it has all other property as infinite volume limit and low temperature limit preserves those property. It is not hard to prove now that the set of ground states does not even form a simplex. A valid question that we can ask at this stage whether the symmetry that we have got so far for $H_{XXX}$ is enough to make action transitive on the set of it's extremal points? 

\vsp 
Now we briefly discuss the situation when $d=3$ i.e. $s=1$. In such a case Pauli spin matrices are given by  
\ben
\sigma_x= 2^{-{1 \over 2}} \left (\begin{array}{llll} 0&,&\;\; 1,\;\; 0 \\ 1&,&\;\;0,\;\;1 \\ 0&,&\;\;1,\;\;0
\end{array} \right ),
\een
\ben
\sigma_y = 2^{-{1 \over 2}} \left (\begin{array}{llll} 0&,&\;\; -i,\;\; 0 \\ i&,&\;\;0,\;\;-i \\ 0&,&\;\;i,\;\;0
\end{array} \right ),
\een
\ben
\sigma_z= \left (\begin{array}{llll} 1&,&\;\; 0,\;\; 0 \\ 0&,&\;\;0,\;\;0 \\ 0&,&\;\;0,\;\;-1
\end{array} \right ).
\een
A direct calculation shows that the intertwiner $g_0$ is a matrix with real entries given below 
\ben
g_0 = \left (\begin{array}{llll} 0&,&\;\; 0,\;\; -1 \\ 0&,&\;\;1,\;\;0 \\ -1&,&\;\;0,\;\;0
\end{array} \right ).
\een
Thus we have $\overline{g_0}=g_0$ and unitary operator $v_{g_0}$ given in relation (31) is self adjoint. 
So far we have not used reflection positivity. In such case we will investigate now $v_{g_0}$. We write relation 
(31) now below:
$$v_{g_0}v_1v_{g_0}^*=-\clj \tilde{v}_3 \clj$$
$$v_{g_0}v_2v_{g_0}^*=\;\clj \tilde{v}_2 \clj$$
$$v_{g_0}v_3v_{g_0}^*=-\clj \tilde{v}_1 \clj$$

Going along the line of the proof for $d=2$, using reflection positivity and purity of the state $\omega$ we show also now that either 
$v_1 = \clj\tilde{v}_3 \clj$ or $v_1=-\clj \tilde{v}_3 \clj$. Thus $v_1\Omega$ and $v_2\Omega$ are eigen vectors of $v_{g_0}$ with both 
eigenvalues are equal to $1$ or $-1$ simultaneously. By the same reason $v_2 = \clj\tilde{v}_2 \clj$ or $v_2=-\clj \tilde{v}_2 \clj$ i.e. 
$v_2\Omega$ is an eigenvector with eigenvalue either $1$ or $-1$. Now we use $SU(2)$-invariance property of $\omega$ to conclude from 
consistency of covariance relation that eigen value for $v_2\Omega$ is equal to that of $v_1\Omega$. Thus $v_{g_0}$ is either $1$ or 
$v_{g_0}=\beta_{-1}$. Thus we find 
$$\epsilon v_1=-\clj \tilde{v}_3 \clj$$
$$\epsilon v_2=\;\clj \tilde{v}_2 \clj$$
$$\epsilon v_3=-\clj \tilde{v}_1 \clj$$
where $\epsilon$ is either $1$ or $-1$. In case $H$ is trivial, then $\epsilon = 1$ 
and at this point it is not clear whether converse is true.  

\vsp 
Our result says very little about integer spin $H_{XXX}$ model as $G-$ ergodic property in Theorem 6.4 is not evident as of now 
for group $G=SU(2)$.  \qed

\bigskip
{\centerline {\bf REFERENCES}}

\begin{itemize} 
\bigskip
\item{[Ac]} Accardi, L. : A non-commutative Markov property, (in Russian), Functional.  anal. i Prilozen 9, 1-8 (1975).

\item{[AcC]} Accardi, Luigi; Cecchini, Carlo: Conditional expectations in von Neumann algebras and a theorem of Takesaki.
J. Funct. Anal. 45 (1982), no. 2, 245–273. 

\item{[AM]} Accardi, L., Mohari, A.: Time reflected
Markov processes. Infin. Dimens. Anal. Quantum Probab. Relat. Top., vol-2
,no-3, 397-425 (1999).

\item{[AKLT]} Affleck, L.,Kennedy, T., Lieb, E.H., Tasaki, H.: Valence Bond States in Isotropic Quantum Antiferromagnets, Commun. Math. Phys. 115, 477-528 (1988). 

\item{[AL]} Affleck, L., Lieb, E.H.: A Proof of Part of Haldane's Conjecture on Spin Chains, Lett. Math. Phys, 12,
57-69 (1986). 

\item{[AMa]} Araki, H., Matsui, T.: Ground states of the XY model, Commun. Math. Phys. 101, 213-245 (1985).

\item{[Ar1]} Arveson, W.: On groups of automorphisms of operator algebras, J. Func. Anal. 15, 217-243 (1974).

\item{[Ar2]} Arveson, W.: Pure $E_0$-semigroups and absorbing states. 
Comm. Math. Phys 187 , no.1, 19-43, (1997)

\item{[Be]} Bethe, H.: Zur Theorie der Metalle. I. Eigenwerte und Eigenfunktionen der linearen Atomkette. (On the theory of metals. I. Eigenvalues and eigenfunctions of the linear atom chain), Zeitschrift für Physik A, Vol. 71, pp. 205–226 (1931). 

\item{[BR]} Bratteli, Ola,: Robinson, D.W. : Operator algebras
and quantum statistical mechanics, I,II, Springer 1981.

\item{[BJP]} Bratteli, Ola,: Jorgensen, Palle E.T. and Price, G.L.: 
Endomorphism of $\clb(\clh)$, Quantisation, nonlinear partial differential 
equations, Operator algebras, ( Cambridge, MA, 1994), 93-138, Proc. Sympos.
Pure Math 59, Amer. Math. Soc. Providence, RT 1996.

\item{[BJKW]} Bratteli, Ola,: Jorgensen, Palle E.T., Kishimoto, Akitaka and
Werner Reinhard F.: Pure states on $\clo_d$, J.Operator Theory 43 (2000),
no-1, 97-143.    

\item{[BJ]} Bratteli, Ola: Jorgensen, Palle E.T. Endomorphism of $\clb(\clh)$, II, 
Finitely correlated states on $\clo_N$, J. Functional Analysis 145, 323-373 (1997). 

\item{[Cu]} Cuntz, J.: Simple $C\sp*$-algebras generated by isometries. Comm. Math. Phys. 57, 
no. 2, 173--185 (1977).  

\item{[Ev]} Evans, D.E.: Irreducible quantum dynamical
semigroups, Commun. Math. Phys. 54, 293-297 (1977).

\item{[EvK]} Evans, David E.; Kawahigashi, Yasuyuki: Quantum symmetries on operator algebras. Oxford Mathematical Monographs. 
Oxford Science Publications. The Clarendon Press, Oxford University Press, New York, 1998. xvi+829 pp. 

\item{[Ex]} Exel, Ruy : A new look at the crossed-product of a $C^*$-algebra by an endomorphism. (English summary)
Ergodic Theory Dynam. Systems 23 (2003), no. 6, 1733–1750.

\item{[FNW1]} Fannes, M., Nachtergaele, B., Werner, R.: Finitely correlated states on quantum spin chains,
Commun. Math. Phys. 144, 443-490(1992).

\item{[FNW2]} Fannes, M., Nachtergaele, B., Werner, R.: Finitely correlated pure states, J. Funct. Anal. 120, 511-
534 (1994). 

\item{[Fr]} Frigerio, A.: Stationary states of quantum dynamical semigroups, Comm. Math. Phys. 63 (1978) 269-276.

\item{[FILS]} Fr\"{o}hlich, J., Israel, R., Lieb, E.H., Simon, B.: Phase Transitions and Reflection Positivity-I 
general theory and long range lattice models, Comm. Math. Phys. 62 (1978), 1-34 

\item{[Ha]} Hall, B.C. : Lie groups, Lie algebras and representations: an elementary introduction, Springer 2003. 

\item{[Li]} Liggett, T.M. : Interacting Particle Systems, Springer 1985. 

\item{[Mo1]} Mohari, A.: Markov shift in non-commutative probability, J. Funct. Anal. vol- 199 
, no-1, 190-210 (2003) Elsevier Sciences. 

\item{[Mo2]} Mohari, A.: Pure inductive limit state and Kolmogorov's property, J. Funct. Anal. vol 253, no-2, 584-604 (2007)
Elsevier Sciences.

\item{[Mo3]} Mohari, A.: Pure inductive limit state and Kolmogorov's property-II,  http://arxiv.org/abs/1101.5961

\item{[Mo4]} Mohari, A: Jones index of a completely positive map, Acta Applicandae Mathematicae. Vol 108, Number 3, 665-677 
(2009).  

\item{[Mo5]} Mohari, A.: Translation invariant pure states in quantum spin chain and it's split property -II under preparation.  

\item{[Ma1]} Matsui, T.: A characterization of pure finitely correlated states. 
Infin. Dimens. Anal. Quantum Probab. Relat. Top. 1, no. 4, 647--661 (1998).

\item{[Ma2]} Matsui, T.: The split property and the symmetry breaking of the quantum spin chain, Comm. 
Maths. Phys vol-218, 293-416 (2001) 

\item{[Po]} Popescu, G.: Isometric dilations for infinite sequences of non-commutating operators, Trans. Amer. Math.
Soc. 316 no-2, 523-536 (1989)

\item{[Pow1]} Powers, Robert T.: Representations of uniformly hyper-finite algebras and their associated von Neumann. rings, Annals of Math. 86 (1967), 138-171.

\item{[Pow2]} Powers, Robert T.: An index theory for semigroups of $*$-endomorphisms of
$\clb(\clh)$ and type II$_1$ factors.  Canad. J. Math. 40 (1988), no. 1, 86--114.

\item{[Ru]} Ruelle, D. : Statistical Mechanics, Benjamin, New York-Amsterdam (1969) . 

\item{[Sa]} Sakai, S. : Operator algebras in dynamical systems. The theory of unbounded derivations in $C\sp *$-algebras. 
Encyclopedia of Mathematics and its Applications, 41. Cambridge University Press, Cambridge, 1991. 

\item{[Ta]} Takesaki, M. : Theory of Operator algebras II, Springer, 2001.

\item{[Wa]} Wassermann, Antony : Ergodic actions of compact groups on operator algebras. I. General theory. Ann. of Math. (2) 130 (1989), no. 2, 273–319

\end{itemize}

\end{document}